\documentclass{amsart}

\usepackage[a4paper,left=40mm,right=35mm,top=30mm,bottom=25mm]{geometry}

\usepackage{tikz}
\usepackage{xcolor}
\usepackage{amssymb,latexsym,amsmath,extarrows}

\usepackage{geometry}
\usepackage{graphicx,mathrsfs,comment}
\usepackage{hyperref,url}

\usepackage{enumitem}

\usepackage{amstext}
\usepackage{bbm}

\usepackage{tabularx}

\numberwithin{equation}{section}

\makeatletter
\newcommand\@avprod[2]{%
  {\sbox0{$\m@th#1\prod$}%
   \vphantom{\usebox0}%
   \ooalign{%
     \hidewidth
     \smash{\vrule height\dimexpr\ht0+1pt\relax depth\dimexpr\dp0+1pt\relax}%
     \hidewidth\cr
     $\m@th#1\prod$\cr
   }%
  }%
}
\newcommand{\avprod}{\mathop{\mathpalette\@avprod\relax}\displaylimits}

\usepackage{esint}

\newtheorem{theorem}{Theorem}[section]
\newtheorem{lemma}[theorem]{Lemma}

\newtheorem{proposition}[theorem]{Proposition}
\newtheorem{remark}[theorem]{Remark}

\newtheorem{definition}[theorem]{Definition}
\newtheorem{corollary}[theorem]{Corollary}

\newtheorem{conjecture}[theorem]{Conjecture}

\newcommand{\al}{\alpha}

\newcommand{\ga}{\gamma}
\newcommand{\Ga}{\Gamma}
\newcommand{\de}{\delta}

\newcommand{\e}{\epsilon}

\newcommand{\cj}{\mathcal J}

\newcommand{\ci}{\mathcal I}

\newcommand{\p}{\partial}

\newcommand{\wt}{\widetilde}

\newcommand{\bQ}{{\bf Q}}

\newcommand{\Tau}{\mathcal{T}}
\newcommand{\bv}{{\bf v}}
\newcommand{\bw}{{\bf w}}

\newcommand{\cB}{{\mathcal B}}

\newcommand{\R}{\mathbb{R}}

\newcommand{\T}{\mathbb{T}  }

\newcommand{\N}{\mathbbm{N}}

\newcommand{\mc}{\mathcal}

\newcommand{\dist}{{\rm dist}}

\newcommand{\fd}{\mathfrak d}

\newcommand{\rank}{\textup{rank}}

\newcommand{\spn}{\textup{span}}
\newcommand{\Sym}{\textup{Sym}}
\newcommand{\Rrank}{\textup{Row-rank}}

\newcommand{\cM}{\mathcal{M}}

\begin{document}

\title[]{Restriction estimates for quadratic manifolds of arbitrary codimensions}

\author{Shengwen Gan} \address{Shengwen Gan\\  Deparment of Mathematics, Massachusetts Institute of Technology, USA}\email{shengwen@mit.edu}

\author{Larry Guth} \address{ Larry Guth\\  Deparment of Mathematics, Massachusetts Institute of Technology, USA}\email{}

\author{Changkeun Oh}\address{ Changkeun Oh\\  Deparment of Mathematics, Massachusetts Institute of Technology, USA} \email{coh@mit.edu}

\maketitle

\begin{abstract}
The restriction conjecture is one of the famous problems in harmonic analysis. There have been many methods developed in the study of the conjecture for the paraboloid.
In this paper, we generalize the multilinear method of Bourgain and Guth for the paraboloid, and obtain restriction estimates for all quadratic manifolds of arbitrary codimensions. In particular, our theorem recovers the main theorem of Bourgain and Guth for the paraboloid. A new ingredient is a covering lemma for varieties whose proof relies on Tarski's projection theorem in real algebraic geometry. We also provide algorithms to compute several algebraic quantities that naturally appear in the argument. These algorithms rely on a cylindrical decomposition in real algebraic geometry.
\end{abstract}



\section{Introduction}

\subsection{Backgrounds}

Introduced by Elias Stein in the 1970s, the restriction problem is a key model problem for understanding  general oscillatory integral operators. 
Consider the $d$-dimensional manifold in $\mathbb{R}^{d+n}$
\begin{equation}
    \mathcal{M}=\{ (\xi,\Phi(\xi)): \xi\in [0,1]^d    \}
\end{equation}
for a measurable function $\Phi$. Introduce an extension operator associated to $\mathcal{M}$
\begin{equation}
    E^{\mathcal{M}}f(x)=\int_{[0,1]^d} f(\xi)e^{2\pi i x \cdot (\xi,\Phi(\xi))}\,d\xi.
\end{equation}
The $L^p \rightarrow L^p$ restriction problem is to find the optimal range of $p$ satisfying
\begin{equation}\label{0802.13}
    \|E^{\mathcal{M}}f\|_{L^p(\R^{d+n})} \leq C_{p}\|f\|_{L^p([0,1]^d)}
\end{equation}
for all measurable functions $f$.

There has been a lot of work studying restriction estimates for manifolds of codimension one, in particular, the paraboloid. After a series of significant progress on a restriction estimate for the paraboloid (for example, \cite{MR1097257, MR1625056,  MR2033842}), a breakthrough has been made by \cite{MR2860188}. They introduced a broad-narrow analysis, which is later used to solve major open problems, for example, the main conjecture in Vinogradov's mean value theorem \cite{MR3548534} and pointwise convergence for  Schrodinger equation \cite{MR3702674, MR3961084}. We refer to  \cite{MR3966819, MR3939285, guth2022decoupling} for the survey for Vinogradov's mean value theorem and \cite{hickman2023pointwise} for the pointwise convergence problem.
\\

Compared to the paraboloid, restriction estimates for manifolds of codimension larger than one are rarely known.
Michael Christ \cite{christthesis, MR766216}  and Mockenhaupt \cite{Mockenhaupt} initiated the study of restriction estimates for such manifolds. Notably, a sharp $L^2$ restriction theorem is obtained for well-curved manifolds (see Definition \ref{def15}) by \cite{Mockenhaupt}. While some understanding has been achieved for this special class of manifolds of codimension two, it remains open whether there is a systematic approach to study restriction estimates for quadratic manifolds of arbitrary codimensions. In this paper, we generalize a broad-narrow analysis to the study of restriction estimates for quadratic manifolds of arbitrary codimensions. In particular, our theorem recovers the main theorem of \cite{MR2860188} for the paraboloid (see Theorem \ref{0705.thm16}).
We refer to references in the introduction of \cite{MR3694011, MR4525760} for the previous results on restriction estimates for quadratic manifolds of codimension larger than one. 
\\

Our main theorem (Theorem \ref{restriction}) involves several definitions. To motivate our definitions and theorem,
let us first explain 
the main ideas of the broad-narrow analysis for the paraboloid. Recall the definition of transversality for the paraboloid: For
$k$ points $p_i \in \mathcal{M}$, we say that the points $p_i$ are transverse provided that
\begin{equation}\label{0726.14}
    n(p_1) \wedge \cdots \wedge n(p_{k}) \gtrsim 1,
\end{equation}
where $n(p_i)$ is the normal vector of $\mathcal{M}$ at the point $p_i$.
The broad-narrow analysis for the paraboloid goes as follows.
Let $K$ be a sufficiently large number. Decompose $[0,1]^d$ into squares $\tau$ with side length $K^{-1}$, and we write
\begin{equation}
    E^{\mathcal{M}}f(x)=\sum_{\tau }E^{\mathcal{M}}(1_{\tau}f)(x).
\end{equation}
Denote by $c_{\tau}$ the center of $\tau$.
For a given point $x$, we consider a collection of squares $\tau$ making a significant contribution to the value $|E^{\mathcal{M}}f(x)|$. The point $x$ is called a $k$-broad point provided that there exist $k$ many significant squares $\tau_i$ such that the vectors $n(c_{\tau_i},\Phi(c_{\tau_i}) )$ are transverse. In such a situation, we apply a $k$-linear restriction estimate by \cite{MR2275834}, and this takes care of $k$-broad points. If $x$ is not a $k$-broad point, then it is called a $k$-narrow point. The advantage of a $k$-narrow point is that all the significant squares lie in a small neighborhood of a $(k-1)$-dimensional plane $H$. Hence, this can be thought of as a lower dimensional problem, and we apply a decoupling inequality for the paraboloid restricted to the plane $H$. Finally, we optimize $k$ which gives the best estimate.

To generalize the broad-narrow analysis to higher codimensions, we first need to introduce a notion of transversality for manifolds of arbitrary codimensions. Compared to the paraboloid (see \eqref{0726.14}), it is not entirely clear how to define transversality and what it means geometrically. In the next subsection, we define transversality so that in our applications a narrow case can be thought of as a lower dimensional problem.

\subsection{Main results}

We need to introduce several notations to state our main theorem. Let $d,n \geq 1$. Denote by $\mathbf{Q}(\xi)=(Q_1(\xi),\ldots,Q_n(\xi))$ an $n$-tuple of real quadratic forms in $d$ variables. Consider the $d$-dimensional manifold in $\mathbb{R}^{d+n}$
\begin{equation}\label{0714.11}
    \mathcal{M}=\{ (\xi,\mathbf{Q}(\xi)): \xi\in [0,1]^d    \}.
\end{equation}
For each $\xi \in [0,1]^d$, we denote by $V_{\xi}$ the tangent space of $\mathcal{M}$ at the point $(\xi,\mathbf{Q}(\xi))$, and 
 by $\pi_{V_{\xi}}$ the orthogonal projection from $\R^{d+n}$ onto $ V_{\xi}$. Let $\mathcal{P}(\frac1K)$ be the partition of $[0,1]^d$ into dyadic cubes of the side length $\frac1K$. (For our convenience, we always assume $K\in 2^\N$.)
 
 We introduce the notation of transversality for manifolds of arbitrary codimensions. 
For a set $A \subset \R^d$, we denote by $\dim A$ the Hausdorff dimension of the set $A$.

\begin{definition}[Transversality: $\theta$-uniform]\label{0614.def11}
Let $\theta>0$ and $M, K \geq 1$.
Let $\{ X_m\}_{m=1}^{d+n}$ be nonnegative integers and $\{\tau_j\}_{j=1}^{M}$ be a subset of $\mathcal{P}(\frac1K)$. 

We say that $\{\tau_j \}_{j=1}^{M}$ is $\theta$-uniform with the controlling sequence $\{X_m \}_{m=1}^{d+n}$ if for each $1\le m\le d+n$ and any subspace $V \subset \mathbb{R}^{d+n}$ with $\dim V=m$ there are at most $\theta M$ many $\tau_j$ intersecting
\begin{equation}\label{lowdimset}
    \{ \xi \in \mathbb{R}^d: \mathrm{dim} (\pi_{V_\xi}(V)) < X_m  \}.
\end{equation}
\end{definition}

\begin{remark}
    {\rm
    The $\theta$-uniform condition roughly says that the squares $\{\tau_j\}_{j=1}^{M}$ cannot be clustered near a low dimensional set of the form \eqref{lowdimset}. This non-stacking condition is analogous to the transversality condition defined for the squares of the paraboloid. We also remark that if the number $X_m$ is larger, then the $\theta$-uniform condition is stronger.
    }
\end{remark}

\begin{definition}\label{0618.def12}
    Let $\cM$ be a manifold of the form \eqref{0714.11}. Let $2\le k\le d+1$ and $ 0\le m\le d+n$ be integers.
    Define $X(\mathcal{M},k,m)$ to be the biggest integer $X \leq d+1$ such that
\begin{equation}\label{ineqX}
\sup_{\dim V=m}
   \dim \{ \xi \in \mathbb{R}^d: \mathrm{dim} (\pi_{V_\xi}(V)) < X  \} \leq k-2.
\end{equation}
\end{definition}

\begin{remark}
    {\rm We have trivial inequalities
    \begin{equation}
        0 \leq X(\cM,k,m) \leq m.
    \end{equation}
    Let us give a proof.
    Since \eqref{ineqX} always holds with $X=0$, we have $X(\cM,k,m)\ge 0$. When $X=m+1$, the left hand side of \eqref{ineqX} equals $d$ which is bigger than the right hand side of \eqref{ineqX} (recall that $k \leq d+1$), so $X(\cM,k,m)\le m$.
    }
\end{remark}

\bigskip

Let us now introduce an extension operator, which is the main object of our paper. For any set $\square \subset [0,1]^d$, define the extension operator associated with $\mathcal{M}$ by
\begin{equation}
E^{\mathcal{M}}_{\square}f(x):=\int_{\square}f(\xi) e^{2 \pi i(x \cdot (\xi, \mathbf{Q}(\xi) ))}\,d\xi.
\end{equation}
For a ball $B(x,\rho^{-1})$, define the weight function $w_{B}:\R^{d+n} \rightarrow \R$ by
\begin{equation}\label{0728.111}
    w_{B(x,\rho^{-1})}(z):=\frac{1}{(1+\rho|{x-z}|)^{100(d+n)}}.
\end{equation}
Lastly, we define a decoupling constant.

\begin{definition}[Decoupling constant] Let $d,n \geq 1$ and $2 \leq k \leq d+1$. Let ${D}_p(\mathcal{M}|_{L_{k-2}},\delta)$ be the smallest constant $D$ such that
    \begin{equation}
    \begin{split}
        \|E^{\mathcal{M}}_{N_{\de}L}f \|_{L^p(B_{\delta^{-2}})} \leq D \Big( \sum_{ \theta \in \mathcal{P}(\delta) } \|E^{\mathcal{M}}_{\theta \cap N_{\de}L}f \|_{L^p(w_{\delta^{-2}})}^p \Big)^{\frac1p}
    \end{split}
    \end{equation}
for any $N_{\de}L$ being the $\de$-neighborhood of a $(k-2)$-dimensional linear subspace $L$, and for all measurable functions $f$. Define
\begin{equation}
    \mathrm{Dec}_p(\mathcal{M}|_{L_{k-2}}):=\limsup_{\delta \rightarrow 0}\frac{\log D_p(\mathcal{M}|_{L_{k-2}},\delta)}{\log(\delta^{-1})}.
\end{equation}
\end{definition}

We are now ready to state the main theorem.

\begin{theorem}\label{restriction}
Let $d,n \geq 1$.
For every $2 \leq p < \infty$ and $\epsilon>0$,
\begin{equation}\label{ineqdecoupling}
\begin{split}
\|  E_{[0,1]^d}^{\mathcal{M}}f \|_{L^p(B_{\delta^{-1}})}  \leq C_{p,\epsilon} \delta^{- \epsilon} \min_{2 \leq k \leq d+1} &\Big(\delta^{d-\frac{2d+2n}{p}}\delta^{-\mathrm{Dec}_p(\mathcal{M}|_{L_{k-2}} )}
\\&+
       \sup_{0 \leq m \leq d+n}(\delta^{-\frac{m}{p}+\frac12{X(\mathcal{M},k,m)}}) \Big) \|f\|_{L^p}.
\end{split}
\end{equation}
\end{theorem}
Let us digest the quantity on the right hand side of \eqref{ineqdecoupling}. We use the broad-narrow analysis to prove Theorem \ref{restriction}.
Inside the ``$\min_{2\le k\le d+1}$", there are two numbers. The first one comes from the estimate of the narrow part, where we apply the decoupling inequalities for low dimensional sets. The second one comes from the estimate of the broad part, where we apply $k$-dimensional restriction estimates (see Theorem \ref{klinearthm}). Lastly, we optimize $k$ to obtain the quantity on the right hand side of \eqref{ineqdecoupling}. The special case $k=d+1$ is implicitly proved in \cite{MR4541334} (see Corollary 2.3 therein). 

The restriction estimate \eqref{0802.13} is false for any $p>2$ provided that $\mathbf{Q}=(Q_1,\ldots,Q_n)$ is independent of some variable or linearly dependent. If $\mathbf{Q}$ does not miss any variable and linearly independent, as observed in \cite{MR4541334} (see Corollary 2.3 therein), 
our theorem gives a restriction estimate \eqref{0802.13} for sufficiently large $p$.  

There are two mysterious numbers on the right hand side of \eqref{ineqdecoupling}:
\begin{equation}\label{0615.16}
    \mathrm{Dec}_p(\mathcal{M}|_{L_{k-2}}),\;\; \mathrm{and} \;\; X(\mathcal{M},k,m).
\end{equation}
One natural question is whether it is possible to calculate these quantities for a given quadratic manifold $\mathcal{M}$.
The first quantity is the decoupling constant for quadratic forms, which is obtained 
in \cite{MR4541334}. In that paper, the decoupling constant is characterized by ``the minimal number of variables''. Let us give a definition of it.

\begin{definition}\label{0706.def15}
For a tuple ${\bf{Q}}=(Q_1(\xi),\ldots,Q_n(\xi))$ of quadratic forms in $d$ variables, denote
\begin{equation}
    NV({\bf Q}):=| \{ 1 \leq d' \leq d: \partial_{\xi_{d'} }Q_{n'} \not\equiv 0 \mathrm{ \; for \; some \; } 1 \leq n' \leq n \}  |.
\end{equation}
    For $0 \leq n' \leq n$ and $0 \leq d' \leq d$,
\begin{equation}
    \fd_{d',n'}({\bf Q}):= \inf_{\substack{M \in \R^{d \times d} \\ \mathrm{rank}(M)=d' } }
    \inf_{\substack{M' \in \R^{n \times n} \\ \mathrm{rank}(M')=n' } } NV(M' \cdot ({\bf Q} \circ M)).
\end{equation}
\end{definition}
By the decocupling theorem in \cite{MR4541334}, the quantity $\mathrm{Dec}_p(\mathcal{M}|_{L_{k-2}})$ can be expressed in terms of $\fd_{d',n'}({\bf Q})$.
In Appendix A, we provide an algorithm to compute $\fd_{d',n'}({\bf Q})$. In Appendix B, we provide an algorithm to compute $X(\mathcal{M},k,m)$. These algorithms rely on a cylindrical decomposition in real algebraic geometry.
\\

We use the broad-narrow analysis introduced in \cite{MR2860188} to prove Theorem \ref{restriction}. One ingredient is the $k$-linear restriction estimate for manifolds of codimension larger than one. For the case of the paraboloid, the $k$-linear restriction estimate was proved in \cite{MR2275834}. We generalize the method there by using a scale-dependent Brsacamp-Lieb inequality in \cite{MR4395082}.

\begin{theorem}[$k$-linear restriction estimate]\label{klinearthm} Let $2 \leq k \leq d+1$ and $\theta>0$. Suppose that $\{\tau_j\}_{j=1}^{M}$ is $\theta$-uniform with the controlling sequence $\{X(\mathcal{M},k,m) \}_{m=1}^{d+n}$. Then for $2 \leq p \leq 2M$ we have
\begin{equation*}
    \Big\| \prod_{j=1}^{M} |E_{\tau_j}^{\mathcal{M}}f|^{\frac1M} \Big\|_{L^p(B_{\delta^{-1}})} \leq C_{\epsilon,p,K,
    \theta} \delta^{-\epsilon-(d+n) \theta/2} \max_{0 \leq m \leq d+n} \delta^{-\frac{m}{p}+\frac12{X(\mathcal{M},k,m)} } \prod_{j=1}^{M}\|f\|_{L^2(\tau_j)}^{\frac1M}.
\end{equation*}
\end{theorem}

Theorem \ref{klinearthm} helps us to deal with the broad part. The narrow part will be dealt with by using a decoupling inequality as in \cite{MR2860188}. However, things become tricky when dealing with manifolds of higher codimension. We would like to compare the codimension one case and the high codimension case. In the codimension one case, say paraboloid, the narrow set is a small neighborhood of a lower dimensional plane, for which we can directly apply a decoupling inequality. However, in the higher codimension case, the narrow set is no longer as simple as before. Probably the only useful information we have about the narrow set is that it lies in a small neighborhood of a variety $Z$ with dimension $k$. To apply a decoupling inequality, we hope the variety $Z$ can be approximated by several $k$-dimensional planes. As far as we can tell, $Z$ could have very bad singularities, and the approximation does not look possible near the singularities. Even though a broad-narrow analysis is widely used in a restriction theory, by this fundamental obstacle related to real algebraic geometry, the broad-narrow analysis has not been fully used in the study of manifolds of arbitrary codimensions. The main contribution of this paper is to devise a covering lemma for varieties to overcome the obstacle.
\medskip

As mentioned, the main novelty of this paper is to prove a covering lemma for varieties (Theorem \ref{algthm}), which helps to deal with the problems arising from the narrow part. Suppose that there is a $k$-dimensional semi-algebraic set $Z\subset [0,1]^d$, and we need a decoupling inequality for a function supported in a small neighborhood of $Z$.
The covering lemma says that a small neighborhood of a $k$-dimensional semi-algebraic set can be covered by a controlled number of  small neighborhoods of $k$-dimensional ``regular" graphs. Furthermore, a bootstrapping argument, originated from the work of Pramanik-Seeger \cite{pramanik2007p}, allows us to approximate each $k$-dimensional ``regular" graph by $k$-planes. Hence, the decoupling problem for a function whose  support is in a small neighborhood of $Z$ is reduced to the decoupling problem for a function whose  support is in a small neighborhood of $k$-planes. This is the reason that $\mathrm{Dec}_p(\mathcal{M}|_{L_{k-2}} )$ appears on the right hand side of \eqref{ineqdecoupling}.  

Although the statement of the covering lemma for varieties might look natural, it took a while for us to figure out the proof.
One main ingredient of the proof of the lemma is the use of a semi-algebraic set and Tarski's theorem in real algebraic geometry. 
A semi-algebraic set has been used for Kakeya-type problems, for example,  \cite{MR3830894, MR3881832, MR3820441,   MR4205111, MR4521046}. We also mention \cite{MR3231483, zahl2023maximal} for the study of a maximal function for a curve and \cite{basu2021stationary} for the study of estimating an oscillatory integral.
\medskip

As a side note, recently, Gressman \cite{gressman2022testing, gressman2023local} made a breakthrough on the study of Radon-like transforms. He gave a complete geometric characterization in all dimensions and codimensions of those Radon-like transforms which satisfy the largest possible range of local $L^p \rightarrow L^q$ inequalities permitted by quadratic-type scaling. It will be interesting to see if his method can be adapted to the study of restriction estimates. We also mention a surprising work \cite{basu2021stationary} on a restriction problem for some family of two-dimensional monomial manifolds of degrees higher than two.

\subsection{Examples}

Let $\mathcal{M}$ be a manifold given by \eqref{0714.11}. The $L^q \rightarrow L^p$ restriction problem is to find the optimal range of $p$ and $q$ satisfying
\begin{equation}\label{0714.110}
    \|E^{\mathcal{M}}f\|_{L^p(\R^{d+n} )} \leq C_{p,q}\|f\|_{L^q([0,1]^d)}
\end{equation}
for all measurable functions $f$. The $L^2$-restriction problem is to find a sharp range of $p$ for $q=2$. The sharp range of $p$ is sometimes called the range of Tomas-Stein. The case $p=q$ draws a lot of attention, and it is called the $L^p \rightarrow L^p$ restriction problem. Note that by H\"{o}lder's inequality it is expected that the range of $p$ where the $L^p \rightarrow L^p$ restriction estimate holds is wider than the range of $p$ for which the $L^2$-restriction estimate holds. 

One interesting example of a manifold is
\begin{equation}
    \{(\xi_1,\xi_2,\xi_1^2,\xi_1\xi_2,\xi_2^2): \xi_1,\xi_2 \in [0,1] \}.
\end{equation}
A restriction estimate for the manifold is obtained for an optimal range of $p$ by  \cite{Mockenhaupt}. As observed in \cite[Claim 2.4]{MR4541334}, our theorem recovers his (see subsubsection \ref{0803.sub133}).
This says that our theorem does not necessarily improve previous results on restriction estimates for manifolds of codimension larger than one.
However, it is very likely that our theorem improves restriction estimates for most manifolds.
In this subsection, we consider concrete examples and obtain  $L^p \rightarrow L^p$ restriction esimates from Theorem \ref{restriction}.

\subsubsection{Codimension one}

Our theorem recovers the main theorem of Bourgain and Guth.

\begin{theorem}[Theorem 
1 of \cite{MR2860188}] \label{0705.thm16}
Let $\mathcal{M}$ be the paraboloid in $\R^{d+1}$. Then the $L^p \rightarrow L^p$ restriction estimate  is true for
\begin{equation}
    p> \min_{2 \leq k \leq d+1} \max{\Big(\frac{2k}{k-1}, \frac{2(2d-k+4)}{2d-k+2}  \Big)}.
\end{equation}
This range of $p$ is asymptotically the same as
  \begin{equation}
  p>2+\frac{3}{d}+O(\frac{1}{d^2}).\footnote{This asymptotics is obtained by using $k=\lfloor \frac{2(d+2)}{3} \rfloor$.}
\end{equation}
\end{theorem}
We give a proof of this theorem in Appendix C.

The restriction estimate of \cite{MR2860188} is further improved by \cite{MR3454378}, \cite{guth2018}, \cite{HR2019}, \cite{hickman2020note} and \cite{MR4484215}.
The best-known bounds of $p$ are obtained by
\cite{wang2022improved} for $d=2$ and by \cite{guo2023dichotomy} for $d \geq 3$. We would like to mention that all the aforementioned improvements over \cite{MR2860188} are built on the broad-narrow analysis of \cite{MR2860188}. For previous results, we refer to references therein.

\subsubsection{Codimension two}

In the case of  codimension two, there is a special type of the manifold known as the \textit{well-curved manifold}. We first give the definition. For a polynomial $P$, denote by $H(P)$ the Hessian matrix of $P$.

\begin{definition}[Well-curved manifold]\label{def15}
    Given a quadratic manifold $\{(\xi,P(\xi),Q(\xi)) \in \mathbb{R}^d \times \mathbb{R} \times \mathbb{R} \}$ of codimension two, define $F(x,y):=\det (x H(P)+yH(Q))$.  Note that $F$ is a homogeneous polynomial. The manifold is called well-curved if
    \begin{itemize}
        \item  $F$ is not identically zero 

        \item F does not have any linear factor of multiplicity larger than $d/2$. 
    \end{itemize}
    Here, a linear factor of $F \in \mathbb{R}[x,y]$ means a homogeneous linear divisor $ax+by$ in $\mathbb{C}[x,y]$, and the multiplicity of $ax+by$ is the largest number $m$ such that $(ax+by)^m$ is still a divisor of $F$.
\end{definition}

This notion of ``well-curved" first appeared in the work of Michael Christ \cite{christthesis, MR766216}  and Mockenhaupt \cite{Mockenhaupt}, where the sharp $L^2$-restriction estimate for well-curved manifolds is proved; the sharp range is $p>2+\frac{8}{d}$. Since the $L^2\rightarrow L^p$ estimate implies the $L^p\rightarrow L^p$, we also know the $L^p\rightarrow L^p$ estimate holds for $p>2+\frac{8}{d}$. It is believed that $L^p \rightarrow L^p$ restriction estimate holds true for a wider range of $p$. However, it looks challenging to improve this range for all well-curved manifolds. The following example explains an obstacle. 
\begin{equation}\label{040814}
    \big\{ (\xi_1,\ldots,\xi_8,\xi_1^2+\xi_2^2-\xi_3^2-\xi_4^2,\xi_5^2+\xi_6^2-\xi_7^2-\xi_8^2): (\xi_1,\ldots,\xi_8) \in [0,1]^8  \big\}.
\end{equation}
This manifold is well-curved, and by the work of Christ and Mockenhaupt, the $L^p \rightarrow L^p$ restriction estimate is known for $p>3$. This manifold is the 2-tensor of the following hyperbolic paraboloid:
\begin{equation}\label{040815}
    \{(\xi_1,\xi_2,\xi_3,\xi_4,\xi_1^2+\xi_2^2-\xi_3^2-\xi_4^2):(\xi_1,\ldots,\xi_4) \in [0,1]^4 \}.
\end{equation}
Since \eqref{040814} is the tensor product of \eqref{040815}, we expect them to have the same range of $p$ for the $L^p\rightarrow L^p$ estimate.
However, the best known range of $p$ for which $L^p\rightarrow L^p$ estimate for hyperbolic paraboloid \eqref{040815} holds is the range of Tomas-Stein ($p>3$).\footnote{We refer to \cite{MR3653943, MR4405679} for the study of restriction estimates for hyperbolic paraboloids in $\R^n$.} 
So, it is hard to improve the $L^p\rightarrow L^p$ estimate for \eqref{040814} over $p>3$ without any progress on a restriction theory for the hyperbolic paraboloid.

For these reasons, we consider a smaller class of manifolds which are more like a paraboloid than a hyperbolic paraboloid. We call them the \textit{good manifolds}.

\medskip

\begin{definition}[Good manifold]
Consider the manifold
\begin{equation}\label{0705.111}
    \mathcal{M}=\{(\xi_1,\ldots,\xi_d,a_1\xi_1^2+\cdots+a_d\xi_d^2,b_1\xi_1^2+\cdots+b_d\xi_d^2) \in [0,1]^d \times \mathbb{R} \times \mathbb{R}  \}.
\end{equation}
This is called a good manifold  if $a_i$'s are all positive and every two by two minor of the following matrix
\begin{equation}\label{0707.115}
    \begin{pmatrix}
        a_1 & a_2 & \cdots & a_n \\
        b_1 & b_2 & \cdots & b_n
    \end{pmatrix}
\end{equation}
has rank two.
\end{definition}
This type of manifolds appeared in the work of Heath-Brown and Pierce \cite{MR3652248}.  If our quadratic manifold $\mathcal{M}$ is given by  $\{(\xi,P(\xi),Q(\xi))\}$ and the Hessian matrices of $P,Q$ are positive definite, then $P,Q$ are simultaneously diagonalizable, and after linear transformation, the manifold $\mathcal{M}$ becomes the form of \eqref{0705.111}. Therefore, good manifolds can be thought of as a natural generalization of a paraboloid to manifolds of codimension two. Note that the rank condition \eqref{0707.115} holds true for a generic choice of $a_i$ and $b_i$. Note also that good manifolds are well-curved.

\begin{theorem}\label{0409thm16}
Let $\mathcal{M}$ be a good manifold. Then the $L^p \rightarrow L^p$ restriction estimate  is true for
\begin{equation}\label{0707.116}
    p> \min_{3 \leq k \leq d+1} \max \Big(\frac{2(k+1)}{k-1}, \frac{2(2d-k+6)}{2d-k+2}\Big).
\end{equation}
This range is asymptotically the same as
  \begin{equation}
  p>2+\frac{6}{d}+O(\frac{1}{d^2}).\footnote{This asymptotics is obtained by using $k=\lfloor \frac{2(d+2)}{3} \rfloor$.}
\end{equation}
\end{theorem}

\begin{figure}[ht]
\centering
\begin{minipage}[b]{0.85\linewidth}
\includegraphics[width=11cm]{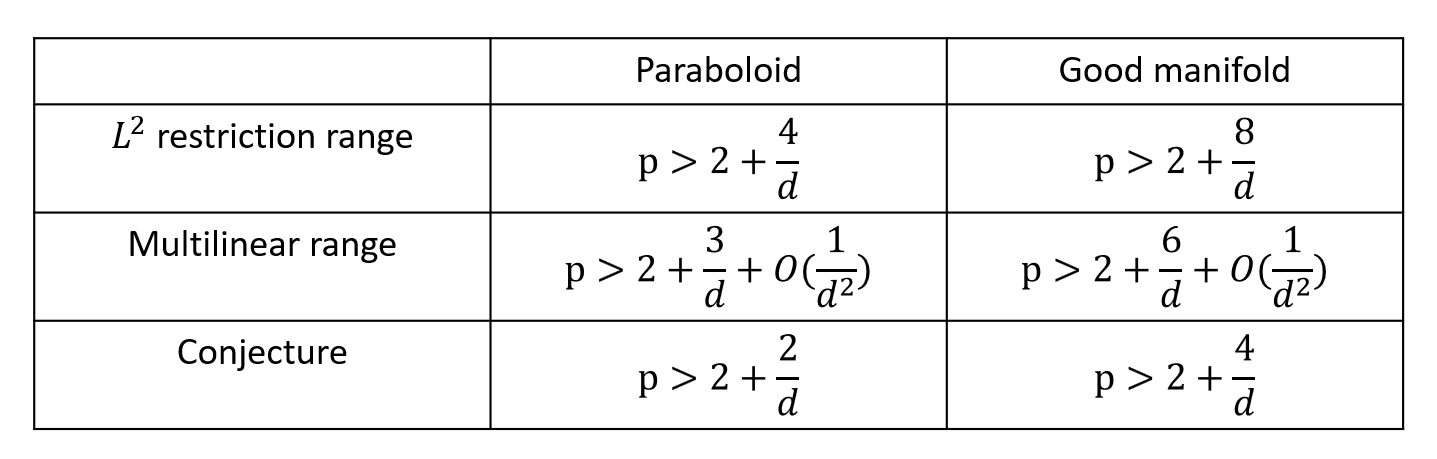}
\caption{Ranges of $p$}
\label{fig:minipage1}
\end{minipage}
\end{figure}

As mentioned before,
Mockenhaupt proved restriction estimates for well-curved manifolds for $p>2+\frac{8}{d}$. The conjectured range is $p>2+\frac{4}{d}$. See Figure \ref{fig:minipage1}.

\begin{remark}\label{0707.remark111}
   {\rm Actually, Theorem \ref{0409thm16} holds for manifolds that satisfy a weaker condition than the good manifold. We can replace the condition ``$a_i>0$" in the definition of good manifold by 
   the condition that there is no nonzero solution $(\xi_1,\ldots,\xi_n)$ to the system 
   \begin{equation}
       \begin{cases}
           a_1\xi_1+\cdots+a_n\xi_n=0 
           \\
           b_1\xi_1+\cdots+b_n\xi_n=0
       \end{cases}
   \end{equation}
   such that $\xi_i\ge 0$ for all $i$.
   Under this mild condition, the restriction estimate still holds for the range of $p$ described  in \eqref{0707.116}.}
\end{remark}

\subsubsection{Maximal codimension}\label{0803.sub133}

For $\alpha=(\alpha_1,\ldots,\alpha_d)\in \N^d$,
denote $\xi^{\alpha}:=\xi_1^{\alpha_1}\cdots \xi_d^{\alpha_d}$, $|\al|=\sum_{i=1}^d\al_i$. Consider 
\begin{equation}\label{0714.119}
 {\bf Q }:=(\xi^{\alpha})_{|\alpha|=2},   \;\;\; \mathcal{M}:=\{(\xi,{ \bf Q}(\xi)) \in \mathbb{R}^d \times \mathbb{R}^n \}.
\end{equation}
Since we are considering all possible monomials of degree two in $d$-variables, we have $n=d(d+1)/2$.  A restriction estimate for an optimal range of $p$ is obtained by \cite{Mockenhaupt} for $d=2$ and by  \cite{MR2064058} for $d \geq 3$. Let us state their theorem.

\begin{theorem}\label{0714.thm112}
    Let $\mathcal{M}$ be the manifold given by \eqref{0714.119}. Then the  $L^p \rightarrow L^p$ restriction estimate  is true for
   $      p>2d+2.$
    This range of $p$ is sharp.
\end{theorem}

Theorem \ref{0714.thm112} can be deduced from Theorem \ref{restriction} by using $k=d+1$. We refer to Claim 2.4 of \cite{MR4541334} and discussions there. As a side note, for this tuple $\bf Q$, the full range of $L^q \rightarrow L^p$ estimates is obtained by \cite{MR2064058, articleoberlin}.

\subsection{Necesary condition}
Denote by
$\mathbf{Q}(\xi)=(Q_1(\xi),\ldots,Q_n(\xi))$ an $n$-tuple of real quadratic forms in $d$ variables.
Consider the $d$-dimensional  manifold in $\mathbb{R}^{d+n}$
\begin{equation}
    \mathcal{M}=\{ (\xi,{\mathbf Q}(\xi)): \xi\in [0,1]^d    \}.
\end{equation}
The $L^q \rightarrow L^p$ restriction problem is to find a range of $p,q$ satisfying
\begin{equation}
    \|E^{\mathcal{M}}_{[0,1]^d}f\|_{L^p(\R^{d+n})} \leq C_{p}\|f\|_{L^q([0,1]^d)}
\end{equation}
for all measurable functions $f$. 
We would like to discuss for which $p,q$ the restriction estimate is false. Let us first introduce a notation. Let $\psi:\R^d \rightarrow \R$ be a non-negative smooth function so that 
\begin{enumerate}
    \item $\psi(\xi)=1$ on $\xi \in [\frac14,\frac34]^d$ \vspace{0.5mm}
    \item 
    $\psi(\xi)=0$ on $\xi \notin [\frac18,\frac78]^{d}$
\end{enumerate}
For a given rectangle $B \subset [0,1]^d$, we say that a function $f$ is adapted to $B$ if
\begin{equation}
    f(\xi)=\psi(L \xi),
\end{equation}
where $L$ is a transformation satisfying $[0,1]^d=L(B)$.

\begin{conjecture}
Given a quadratic manifold, the followings are equivalent.
\begin{enumerate}
    \item The $L^p \rightarrow L^p$ restriction estimate is true for all measurable functions $f$
    \item The $L^p \rightarrow L^p$ restriction estimate is true for a function adapted to $[0,1]^d$.
\end{enumerate}
\end{conjecture}

\begin{conjecture}
Given a quadratic manifold,
the followings are equivalent.
\begin{enumerate}
    \item The $L^q \rightarrow L^p$ restriction estimate is true for all measurable functions $f$
    \item The $L^q \rightarrow L^p$ restriction estimate is true for functions adapted to any rectangles.
\end{enumerate}
\end{conjecture}

We do not know a counterexample to these conjectures. We also wonder whether these conjectures are still true for any smooth manifolds.

\subsection{Notations}
For two non-negative numbers $A_1$ and $A_2$, we write $A_1 \lesssim A_2$ to mean that there exists a constant $C$, such that $A_1 \leq C A_2$. The notation $A_1 \lesssim_{d} A_2$ means that $A_1 \leq C(d)A_2$ for some constant $C(d)$ depending on $d$. Similarly, we use $O(A_1)$ to denote a number whose absolute value is smaller than $CA_1$ for some constant $C$.

For a set $A$, 
\begin{equation}
    \fint_A f(x)\,dx:=\frac{1}{|A|} \int_A f(x) \,dx.
\end{equation}
We use $B(x,\rho^{-1})$ to denote the ball of radius $\rho^{-1}$ centered at $x \in \R^{d+n}$. This is sometimes abbreviated as $B_{\rho^{-1}}$. For a ball $B(x,\rho^{-1})$, define the weight function $w_{B}:\R^{d+n} \rightarrow \R$ by
\begin{equation}
    w_{B(x,\rho^{-1})}(z):=\frac{1}{(1+\rho|{x-z}|)^{100(d+n)}}.
\end{equation}
For a function $F:\R^{d+n} \rightarrow \mathbb{C}$, define the weighted integral
\begin{equation}
\|F\|_{L^2_{\#}(w_{B(x,\rho^{-1}) } ) }:= \Big( \frac{1}{|B(x,\rho^{-1})|} \int_{\R^{d+n}} |F(z)|^2  w_{B(x,\rho^{-1}) }(z)\,dz \Big)^{1/2}.
\end{equation}

For a finite set $A$ and a sequence $\{a_n\}_{n \in A}$, denote
\begin{equation}
    \avprod_{n \in A } |a_n|: = \prod_{n \in A}|a_n|^{\frac{1}{|A|}}.
\end{equation}

For a polynomial $Q$ in $d$ variables, denote
\begin{equation}
    Z_Q:= \{ \xi \in \R^d: Q(\xi)=0 \}.
\end{equation}

\subsection{Acknowledgements}

 We acknowledge funding from the AIM Fourier restriction community.
 LG is supported by a Simons Investigator
award. We would like to thank Shaoming Guo for valuable comments.

\section{\texorpdfstring{$k$}{}-linaer restriction estimates}

In this section, we prove Theorem \ref{0614.thm21}. Theorem \ref{klinearthm} is a special case of it. The proof of Theorem \ref{0614.thm21} is an application of the scale-dependent Brascamp-Lieb inequality with a bootstrapping argument in \cite{MR3300318}.

\begin{theorem}\label{0614.thm21} Let  $\theta,M,K>0$ and $X_m \in \N_{\ge 0}$. Let $\{\tau_j\}_{j=1}^{M} \subset \mathcal{P}(\frac1K)$. Suppose that $\{\tau_j\}_{j=1}^{M}$ is $\theta$-uniform with the controlling sequence $\{X_m \}_{m=1}^{d+n}$. Then for $2 \leq p \leq 2M$ and $\epsilon>0$ we have
\begin{equation*}
    \Big\| \avprod_{j=1}^{M} |E_{\tau_j}^{\mathcal{M}}f| \Big\|_{L^p(B_{\delta^{-1}})} \leq C_{\epsilon,p,\theta,K} \delta^{-\epsilon-(d+n) \theta/2} \max_{0 \leq m \leq d+n} \delta^{-\frac{m}{p}+\frac12{X_m}} \avprod_{j=1}^{M}\|f\|_{L^2(\tau_j)}.
\end{equation*}
\end{theorem}

One ingredient in the proof is a $k$-linear Kakeya inequality.  Let us introduce some notations to state the inequality. Let $\mathcal{M}$ be a manifold, and $\tau_j$ be an element of $\mathcal{P}(\frac1K)$. Suppose that $\mathbb{T}_j$ is a collection of slabs $T_j$ with dimension  
\begin{equation}
 \underbrace{ \delta \times \cdots \times \delta   }_{d\text{-times}} \times
 \underbrace{ 1 \times \cdots \times 1   }_{n\text{-times}} .
\end{equation}
We say that $\T_j$ is associated with $\tau_j$ if for every $T_j \in \T_j$ there exists $\xi \in \tau_j$ such that the subspace spanned by the short directions of $T_j$ is parallel to $V_\xi$ (the tangent space of $\mathcal{M}$ at the point $(\xi,\mathbf{Q}(\xi))$).

\begin{proposition}\label{0304.Guth} Let $\theta,M,K>0$ and $X_m \geq 0$. Suppose that $\{\tau_j \}_{j=1}^{M}$ is $\theta$-uniform with the controlling sequence $\{X_m \}_{m=1}^{d+n}$ and $\T_j$ is associated with $\tau_j$. Then for $0 \leq p \leq M$ and $\epsilon>0$, we have
\begin{equation*}
    \Big\| \avprod_{j=1}^M \Big( \sum_{T_j \in \mathbb{T}_j } \chi_{T_j} \Big) \Big\|_{L^p([0,1]^{d+n})} \lesssim \delta^{-(d+n)\theta-\epsilon} \delta^{\frac{d+n}{p}} \sup_{0 \leq m \leq d+n}(\delta^{-\frac{m}{p}+X_m}) \avprod_{j=1}^{M} (\# \mathbb{T}_j).
\end{equation*}
\end{proposition}

This inequality will easily follow from a scale-dependent Brascamp-Lieb inequality. Let us state the inequality.

\begin{theorem}[Theorem 2 of \cite{MR4395082}]\label{0728.thm23} 
Fix $0 \leq p \leq M$. Let  $V_j^0$ be linear subspaces of $\R^{d+n}$.
Let $\{\pi_j \}_{j=1}^{M}$ be orthogonal projections from $\R^{d+n}$ onto $(V_j^0)^{\bot}$. Then there exists $\nu>0$ satisfying the following: For every $\epsilon>0$,
\begin{equation*}
    \int_{[0,1]^{d+n}} \avprod_{j=1}^M \big( \sum_{T_j \in \T_j} \chi_{T_j} \big)^p \leq C_{\epsilon} \delta^{-\epsilon+d+n}\sup_{V \leq \R^{d+n} }\big(\delta^{-\dim V+\frac{p}{M}\sum_{j=1}^M \dim \pi_j(V)} \big)\avprod_{j=1}^M(\# \T_j)^p
\end{equation*}
holds for all collections $\T_j$ of $\delta$-neighborhoods of $d$-dimensional affine subspaces of $\R^{d+n}$ which, modulo translations, are within a distance $\nu$ of the fixed subspace $V_j^0$.
\end{theorem}

The proof of Proposition \ref{0304.Guth} is rather straightforward. By the scale-dependent Brascamp-Lieb inequality (Theorem \ref{0728.thm23}) and compactness, for given $V$ with $\dim V= m$, it suffices to prove
\begin{equation}\label{0616.23}
    \frac{1}{M}\sum_{j=1}^M \dim \pi_{V_{\xi_j}}(V) \geq (1-\theta)X_m
\end{equation}
for all $\xi_j \in \tau_j$.
By Definition \ref{0614.def11}, there are at most $\theta M$ many $\tau_j$ intersecting
\begin{equation}
    \{ \xi \in \mathbb{R}^d: \mathrm{dim} (\pi_{V_\xi}(V)) < X_m  \}.
\end{equation}
In other words,
there are at least $(1-\theta)M$ many $\tau_j$ not intersecting the set. Hence, we have $\dim \pi_{V_{\xi_j}}(V) \geq X_m$ for all points $\xi_j \in \tau_j$ for $(1-\theta)M$ many $j$. This immediately gives \eqref{0616.23} and completes the proof of Proposition \ref{0304.Guth}.

\subsection{Proof of Theorem \ref{0614.thm21}}

The proof of Theorem \ref{0614.thm21} is identical to that for Corollary 9.1 of \cite{MR4541334}. This argument traces back to \cite{MR2275834}. Since the proof is standard nowadays, we only give a sketch here. We refer to \cite{MR4541334} for details. \medskip

We first obtain a ball inflation lemma from Proposition \ref{0304.Guth}. Once we have the ball inflation lemma (Proposition \ref{0617.prop24}), Theorem \ref{0614.thm21} will follow by applying the lemma repeatedly. The proof and numerology are identical to those of Corollary 9.1 in \cite{MR4541334}. To state the ball inflation lemma, recall the weight function $w_{B(x,\rho^{-1})}$ (see \eqref{0728.111}).
This function is adapted to the ball $B(x,\rho^{-1})$ of radius $\rho^{-1}$ centered at  $x \in \R^{d+n}$.
    
\begin{proposition}[cf. Proposition 4.9 of \cite{MR4541334}]\label{0617.prop24} Let  $\theta,M,K>0$ and $X_m \in\N_{\ge 0}$. Let $\{\tau_j\}_{j=1}^{M} \subset \mathcal{P}(\frac1K)$. Suppose that $\{\tau_j\}_{j=1}^{M}$ is $\theta$-uniform with the controlling sequence $\{X_m \}_{m=1}^{d+n}$. Then for  $2 \leq p \leq 2M$,
\begin{equation*}\label{0616.27}
    \begin{split}
        &\Big( \fint_{B(x_0,\rho^{-2})} \avprod_{j=1}^{M} \Big(\sum_{J_j \subset \tau_j: J_j \in \mathcal{P}(\rho) } \|E_{J_j}f\|_{L^2_{\#}(w_{B(x,\rho^{-1})}) }^2 \Big)^{\frac p2} \, dx\Big)^{\frac 1p} 
        \\& \lesssim_{\epsilon} \rho^{-\frac{(d+n)\theta}{2}-\epsilon}\rho^{\frac{d+n}{p}-\frac{d}{2}}\max_{0 \leq m \leq d+n} \big(\rho^{-\frac{m}{p}+\frac{ X_m}{2} }\big)\avprod_{j=1}^{M}\Big(\sum_{\substack{J_j \subset \tau_j: \\ J_j \in \mathcal{P}(\rho) }} \|E_{J_j}f\|_{L^2_{\#}(w_{B(x_0,\rho^{-2})} ) }^2 \Big)^{\frac12}.
    \end{split}
\end{equation*}
    
\end{proposition}

The proof is identical to Proposition 4.9 of \cite{MR4541334}, so we only give a sketch.

\begin{proof}[Sketch of the proof of Proposition \ref{0617.prop24}] 
Let $J_j \in \mathcal{P}(\rho)$. Then by the locally constant property, the function $E_{J_j}f$ is essentially constant on every slab with dimension 
    \begin{equation}
  \underbrace{ \rho^{-1} \times \cdots \times \rho^{-1}   }_{d\text{-times}} \times \underbrace{ \rho^{-2} \times \cdots \times \rho^{-2}   }_{n\text{-times}} 
\end{equation}
and the subspace spanned by the short directions of the slab is parallel to $V_{\xi}$ for $\xi \in J_j$. For each $J_j$, we tile $\R^{d+n}$ by translating the slabs. Denote by $\Tau_{J_j}$ the collection of slabs. We let $T_{J_j}(x) \in \Tau_{J_j}$ be the slab containing $x \in \R^{d+n}$. Define
\begin{equation}\label{0616.29}
    F_{J_j}(x):=\sup_{y \in T_{J_j}(x) }\|E_{J_j}f\|_{L^2_{\#}(w_{B(y,\rho^{-1})} )}^2.
\end{equation}
Then the left hand side of \eqref{0616.27} is bounded by
\begin{equation}\label{0616.210}
    \Big( \fint_{B(x_0,\rho^{-2})} \avprod_{j=1}^{M} \Big(\sum_{J_j \subset \tau_j: J_j \in \mathcal{P}(\rho) } F_{J_j}(x) \Big)^{\frac p2} \, dx\Big)^{\frac 1p}.
\end{equation}
Note that $F_{J_j}$ is constant on any $T_{J_j} \in \Tau_{J_j}$. By pigeonholing and homogeneity, we may assume that $F_{J_j}$ is either one or zero on each $T_{J_j}$. To apply Proposition \ref{0304.Guth}, we do scaling: $x \mapsto \rho^{-2}y$. Then \eqref{0616.210} is equal to 
\begin{equation}\label{0616.211}
    \Big( \fint_{B(x_0,1)} \avprod_{j=1}^{M} \Big(\sum_{J_j \subset \tau_j: J_j \in \mathcal{P}(\rho) } F_{J_j}(\rho^{-2}x) \Big)^{\frac p2} \, dx\Big)^{\frac 1p}.
\end{equation}
Now we can apply Proposition \ref{0304.Guth} with $\delta=\rho$ and $p$ replaced by $p/2$. Then this is bounded by
\begin{equation}\label{0616.212}
    \Big(\rho^{-(d+n)\theta-\epsilon}\rho^{\frac{2(d+n)}{p}}\max_{0 \leq m \leq d+n}\rho^{-\frac{2m}{p}+X_m }\avprod_{j=1}^{M}(\# \mathbb{T}_j)\Big)^{\frac12}.
\end{equation}
Note that
\begin{equation}\label{0616.213}
\begin{split}
    \#\T_j &= 
\rho^{-d} \fint_{B(x_0,1)} F_{J_j}(\rho^{-2}x)\,dx
\\&=\rho^{-d}\fint_{B(x_0,\rho^{-2})}F_{J_j}(x)\,dx \lesssim \rho^{-d} \|E_{J_j}f\|_{L^2_{\#}(w_{B(x_0,\rho^{-2})} )}^2.
\end{split}
\end{equation}
By \eqref{0616.29} and \eqref{0616.213}, we see that \eqref{0616.212} is further bounded by
\begin{equation*}
    \rho^{-\frac{(d+n)\theta}{2}-\epsilon}\rho^{\frac{d+n}{p}-\frac{d}{2}}\max_{0 \leq m \leq d+n} \big(\rho^{-\frac{m}{p}+\frac12 X_m }\big)\avprod_{j=1}^{M}\Big(\sum_{J_j \subset \tau_j: J_j \in \mathcal{P}(\rho) } \|E_{J_j}f\|_{L^2_{\#}(w_{B(x_0,\rho^{-2})} ) }^2 \Big)^{\frac12}.
\end{equation*}
This finishes the proof.
\end{proof}

Let us now give a proof of Theorem \ref{0614.thm21}.

\begin{proof}[Sktech of the proof of Theorem \ref{0614.thm21}]
    We will show that for any ball $B$ of radius $\delta^{-1}$,
    \begin{equation}\label{0618.214}
    \begin{split}
        &\Big\| \avprod_{j=1}^{M} |E_{\tau_j}^{\mathcal{M}}f| \Big\|_{L^p_{\#}(B)} \lesssim_{\epsilon} 
        \\&
        \delta^{-{(d+n)\theta/2}-\epsilon}\delta^{\frac{d+n}{p}-\frac{d}{2}}\max_{0 \leq m \leq d+n} \big(\delta^{-\frac{m}{p}+\frac12 X_m} \big)\avprod_{j=1}^{M}\big(\sum_{\substack{\square \in \mathcal{P}(\delta):  \\ \square \subset \tau_j } }\|E^{\mathcal{M}}_{\square}f\|_{L^2_{\#}(w_B)}^2\big)^{\frac12}.
        \end{split}
    \end{equation}
    Recall that we do not take an average over a ball on the integrals in Theorem \ref{0614.thm21}.
Theorem \ref{0614.thm21} follows from combining the above inequality with
\begin{equation}
\|E^{\mathcal{M}}_{\square}f\|_{L^2_{\#}(w_B)} \lesssim \delta^{\frac{d}{2}}\|f\|_{L^2(\square)}. 
\end{equation}
It remains to show \eqref{0618.214}. Note that the exponent of $\delta$ on the right hand side of \eqref{0618.214} is double the exponent of $\rho$ for Proposition \ref{0617.prop24}. The inequality \eqref{0618.214} follows by applying Proposition \ref{0617.prop24} repeatedly so that the frequency scales grow up as follows. 
\begin{equation}
    \delta^{\epsilon} \rightarrow \delta^{2\epsilon} \rightarrow \delta^{2^2\epsilon} \rightarrow \cdots \rightarrow \delta^{2^{N}\epsilon}:=\delta.
\end{equation}
This argument is identical to that of (9.5)--(9.7) in \cite{MR4541334}. We leave out the details.
\end{proof}

\section{Proof of Theorem \ref{restriction}}

In this section, we give the proof of Theorem \ref{restriction}. Our main method is the broad-narrow analysis in which there is the dichotomy: broad case and narrow case. Theorem \ref{0614.thm21} will take care of the broad case. To deal with the narrow case, we will use Theorem \ref{0619.thm32}. To state the theorem, let us first give the definition of a semi-algebraic set.

\begin{definition}\label{defsemialg}
    A set $Z\subset \R^d$ is called a semi-algebraic set if it can be written as a finite union of sets of the form
    \[ \{x\in \R^d: P_1(x)=0,\dots, P_l(x)=0, P_{l+1}(x)>0,\dots, P_{l+k}(x)>0  \}, \]
where $P_1(x),\dots,P_{l+k}(x)$ are polynomials. We call the polynomials appearing in the definition of $Z$ the defining polynomials of $Z$. Given a semi-algebraic set, the defining polynomials do not need to be unique. There is a way to represent $Z$ using appropriate polynomials so that the sum of degrees of these polynomials is minimal. We define the complexity of $Z$ to be this minimal number. Define the dimension of $Z$ to be the Hausdorff dimension of the set $Z$.
\end{definition}

Let us now state the theorem.

\begin{theorem}\label{0619.thm32}
Let $0 \leq k \leq d-1$.
Given $E>100$ and $0<\epsilon < \frac{1}{100}$, there exist  constants  $ M=M(d,E),\mu_0=\mu_0(d,E,\epsilon)$, and $c=c(d,E,\epsilon)$, such that the following is true.

Let $Z\subset \R^d$ be a $k$-dimensional semi-algebraic set with complexity $\le E$. For any $0<K^{-1}<\mu_0$, there exist dyadic numbers $\{K_i \}_{i=1}^{M}$ such that
\begin{equation}\label{0729.31}
    K^{{c} } \leq K_1 \leq K_2 \leq \cdots \leq K_M \leq K^{\frac{1}{2}} 
\end{equation}
and $\mathcal{P}_i(K_i^{-1}) \subset \mathcal{P}(K_i^{-1})$ depending on the choice of $Z$ such that
\begin{equation}
   {N}_{K^{-1}}(Z) \cap [0,1]^d \subset \bigcup_{i=1}^{M} \bigcup_{\tau_i \in \mathcal{P}_i(K_i^{-1}) } \tau_i 
\end{equation}
and
\begin{equation}\label{0620.33}
\begin{split}
    \Big\|\sum_{\substack{ \tau_i \in \mathcal{P}_i(K_i^{-1}) }} E^{\mathcal{M}}_{\tau_i}f \Big\|_{L^p} \leq
    C_{\epsilon,p,E} \cdot K_i^{\epsilon+\mathrm{Dec}_p(\mathcal{M}|_{L_{k}})}  \Big( \sum_{\tau_i \in \mathcal{P}(K_i^{-1}) } \|E^{\mathcal{M}}_{\tau_i}f\|_{L^p}^p \Big)^{\frac1p}.
    \end{split}
\end{equation}
The constant $C_{\epsilon,p,E}$ is independent of the choice of $Z$.
\end{theorem}

In the next subsection, we prove Theorem \ref{restriction} under assuming Theorem \ref{0619.thm32}. We postpone the proof of Theorem \ref{0619.thm32} to the next section.

\subsection{Proof of Theorem \ref{restriction}}
For simplicity, we use the notation $Ef$ for $E^{\mathcal{M}}f$.
Fix $2 \leq k \leq d+1$. Fix $\epsilon>0$ and $p$. Take sufficiently large $K$, which will be determined later. In particular, we require the number $K$ to satisfy
\begin{equation}
    p< K^{\epsilon}.
\end{equation}
This condition is required by a technical reason, which is related to the upper bound of $p$ in Theorem \ref{0614.thm21}. It will be clear when we use the theorem.

Recall Definition \ref{0618.def12}. We will show
\begin{equation*}
\begin{split}
\|  E_{[0,1]^d}f \|_{L^p(B_{\delta^{-1}})}  \leq C_{p,\epsilon} \delta^{-100(d+n) \epsilon}  &\Big(\delta^{d-\frac{2d+2n}{p}}\delta^{-\mathrm{Dec}_p(\mathcal{M}|_{L_{k-2}})}
\\&+
       \sup_{0 \leq m \leq d+n}(\delta^{-\frac{m}{p}+\frac12{X(\mathcal{M},k,m)}}) \Big) \|f\|_{p}.
\end{split}
\end{equation*}
Fix $\epsilon>0$. Take sufficiently large $K$, which will be determined later. Fix $B_{K^2}$. Consider
\begin{equation}\label{0620.34}
    \mathcal{C}_0:=\{ \tau \in \mathcal{P}(K^{-1}):\|E_{[0,1]^d}f\|_{L^p(B_{K^2})} \leq K^{10d} \|E_{\tau}f\|_{L^p(B_{K^2})}  \}.
\end{equation}
Recall Definition \ref{0614.def11}. There are two cases: whether $\mathcal{C}_0$ is $\epsilon$-uniform or not. 

If it is $\epsilon$-uniform with the controlling sequence $\{X(\mathcal{M},k,m) \}_{m=0}^{d+n}$, then we have
\begin{equation}\label{0620.35}
  \begin{split}  \|E_{[0,1]^d}f\|_{L^p(B_{K^2})} &\lesssim K^{10d} \avprod_{\tau \in \mathcal{C}_0 } \|E_{\tau}f\|_{L^p(B_{K^2})}
  \\& \lesssim
  K^{50d} \Big\| \avprod_{\tau \in \mathcal{C}_0} E_{\tau}f_{j(\tau)} \Big\|_{L^p(w_{B_{K^2}})}
  \end{split}
\end{equation}
where $f_{j(\tau)}$ is a modulation of $f$. The second inequality is by using random translation and locally constant property. For the readers who are not familiar with the proof of the second line, we refer to $(3.8)$ on page 850 of \cite{MR3961084}.

If it is not $\epsilon$-uniform, then there exists a subspace $V \subset \R^{d+n}$ whose dimension is $m$ such that there are more than $\epsilon |\mathcal{C}_0|$ many $\tau  \in \mathcal{C}_0$ intersecting
\begin{equation}
    Z_V:=\{ \xi \in \R^d: \mathrm{dim} (\pi_{V_\xi}(V)) < X(\mathcal{M},k,m) \}.
\end{equation}
By the definition of $X(\mathcal{M},k,m)$, this set has Hausdorff dimension at most $k-2$. Note that this set is a semi-algebraic set. We apply Theorem \ref{0619.thm32} to $Z_V$ with $k$ replaced by the dimension of $Z_V$. Then we have
\begin{equation}
   {N}_{K^{-1}}(Z_V) \cap [0,1]^d \subset \bigcup_{i=1}^{M} \bigcup_{\tau_i \in \mathcal{P}_i(K_i^{-1}) } \tau_i=:\Tau_0.
\end{equation}
Define $\mathcal{P}_{0,i}(K_i^{-1}):=\mathcal{P}_i(K_i^{-1})$ to indicate that this set appears in the first stage.
Define
\begin{equation}
    \mathcal{C}_1:=\big\{ \tau \in \mathcal{C}_0 : \tau \cap
 \big(\bigcup_{i=1}^{M} \bigcup_{\tau_i \in \mathcal{P}_{0,i}(K_i^{-1}) } \tau_i \big) = \emptyset \big\}.
\end{equation}
Note that the cardinality of $\mathcal{C}_1$ is at most $(1-\epsilon)|\mathcal{C}_0|$. If $\mathcal{C}_2$ is $\epsilon$-uniform, then we stop defining the set. If it is not, then we define $\mathcal{C}_3$ and continue this process. This process terminates after at most $O(\log{K})$ times. After this algorithm, we have
\begin{equation}
    \mathcal{C}_{m+1}=\big\{ \tau \in \mathcal{C}_m : \tau \cap \big(\bigcup_{i=1}^{M} \bigcup_{\tau_i \in \mathcal{P}_{m,i}(K_i^{-1}) } \tau_i \big) = \emptyset \big\}
\end{equation}
for $m=1, \ldots, L$ with $L \lesssim \log K$. We also have $\mathcal{P}_{m,i}(K_i^{-1})$. The collection $\mathcal{C}_{L+1}$ is either $\epsilon$-uniform or an empty set. Define
\begin{equation}
    \Tau:= \bigcup_{m=1}^{L} \bigcup_{i=1}^M \bigcup_{\tau_i \in \mathcal{P}_{m,i}(K_i^{-1}) } \tau_i.
\end{equation}
Recall that the contribution of $(\mathcal{C}_0
)^{c}$ is negligible (see the definition \eqref{0620.34}). By the triangle inequality, we have
\begin{equation}\label{0620.311}
    \|E_{[0,1]^d}f\|_{L^p(B_{K^2})} \lesssim 
    \Big\|\sum_{\tau \in \mathcal{C}_{L+1}} E_{\tau}f\Big\|_{L^p(B_{K^2})}+\Big\|\sum_{\tau \in \Tau} E_{\tau}f\Big\|_{L^p(B_{K^2})}.
\end{equation}
Recall that $\mathcal{C}_{L+1}$ is either $\epsilon$-uniform or empty.
If $\mathcal{C}_{L+1}$ is empty, then the first term on the right hand side is zero. If $\mathcal{C}_{L+1}$ is $\epsilon$-uniform, then we apply the arguments of \eqref{0620.35} to the first term.

Let us focus on the second term of \eqref{0620.311}. Since each $K_i$ is a dyadic number, by the triangle inequality, we have
\begin{equation}
    \Big\|\sum_{\tau \in \Tau} E_{\tau}f\Big\|_{L^p(B_{K^2})} \lesssim 
\sum_{m=1}^L \sum_{i=1}^{M} \Big\| \sum_{\tau_i \in \widetilde{\mathcal{P}}_{m,i}(K_i^{-1}) } E_{\tau_i}f \Big\|_{L^p(B_{K^2})}
\end{equation}
where $\widetilde{P}_{m,i}(K_i^{-1})$ is a subset of $\mathcal{P}_{m,i}(K_i^{-1})$. By \eqref{0620.33} (recall that we used $k-2$ instead of $k$) and the definition of $\mathrm{Dec}(\mathcal{M}|_{k-2})$, and replacing the summation over $m,i$ by maximum, this term is further bounded by
\begin{equation}
    \widetilde{C}_{\epsilon,p,E} \sup_{K^{c} \leq \widetilde{K} \leq K^{1/2 }}
    (\widetilde{K})^{\epsilon+\mathrm{Dec}_p(\mathcal{M}|_{L_{k-2}} )}
    \Big( \sum_{\tau_i \in \mathcal{P}(\widetilde{K}^{-1}) } \|E_{\tau_i}f\|_{L^p(w_{B_{K^2}})}^p \Big)^{\frac1p}.
\end{equation}

To summarize, we have
\begin{equation*}
    \begin{split}
        \|E_{[0,1]^d}f\|_{L^p(B_{K^2})} &\lesssim 
        K^{50d} \Big\| \avprod_{\tau \in \mathcal{C}} E_{\tau}f_{j(\tau)} \Big\|_{L^p(w_{B_{K^2}})}
        \\&
        + \sup_{K^{c} \leq \widetilde{K} \leq K^{1/2}}
    (\widetilde{K})^{\epsilon+\mathrm{Dec}_p(\mathcal{M}|_{L_{k-2}} )}
    \Big( \sum_{\tau_i \in \mathcal{P}(\widetilde{K}^{-1}) } \|E_{\tau_i}f\|_{L^p(w_{B_{K^2}})}^p \Big)^{\frac1p}.
    \end{split}
\end{equation*}
The collection $\mathcal{C}$ depends on a choice of $B_{K^2}$, but the number of the possibilities of the choice is at most $O(K^{d})$. The number of the possibilities of $j(\tau)$ is also $O(K^{C})$, so by summing over $B_{K^2} \subset B_{\delta^{-1}}$, we have
\begin{equation}\label{0621.316}
    \begin{split}
        \|E_{[0,1]^d}f\|_{L^p(B_{\delta^{-1}})} &\lesssim 
        K^{C} \Big\| \avprod_{\tau \in \mathcal{A}} E_{\tau}f_{j(\tau)} \Big\|_{L^p(w_{B_{\delta^{-1}}})}
        \\&
        + \sup_{K^{c} \leq \widetilde{K} \leq K^{1/2 }}
    (\widetilde{K})^{\epsilon+\mathrm{Dec}_p(\mathcal{M}|_{L_{k-2}} )}
    \Big( \sum_{\tau \in \mathcal{P}(\widetilde{K}^{-1}) } \|E_{\tau}f\|_{L^p(w_{B_{\delta^{-1}}})}^p \Big)^{\frac1p}.
    \end{split}
\end{equation}
for some $\epsilon$-uniform set $\mathcal{A}$ with the controlling sequence $\{X(\mathcal{M},k,m) \}_{m=0}^{d+n}$. We next apply Theorem \ref{0614.thm21} with $\{X(\mathcal{M},k,m) \}_{m=0}^{d+n}$ and $\theta=\epsilon$ and obtain
\begin{equation}
    \Big\| \avprod_{\tau \in \mathcal{A}} E_{\tau}f_{j(\tau)} \Big\|_{L^p(w_{B_{\delta^{-1}}})} \lesssim_{\epsilon}  \delta^{-\epsilon-(d+n) \epsilon} \max_{0 \leq m \leq d+n} \delta^{-\frac{m}{p}+\frac12{X_m} } \|f\|_{L^2}.
\end{equation}
This bound is already good enough to close the induction. Let us bound the second term on the right hand side of \eqref{0621.316}. Fix $\tau \in \mathcal{P}(\widetilde{K}^{-1})$. For convenience, assume that $\tau=[0,\widetilde{K}^{-1}]^d$.
We do rescaling; take a linear transform $L$ so that $L(\tau)=[0,1]^d$. Define the function $g$ so that $g(\xi):=f(L^{-1}\xi)$. Then we have
\begin{equation}
    \|E_{\tau}f\|_p \sim (\widetilde{K})^{-d+\frac{d+2n}{p}} \|Eg\|_{p}, \;\;\; \|f\|_{L^p(\tau)} \sim (\widetilde{K})^{-\frac{d}{p}}\|g\|_p.
\end{equation}
After applying the induction hypothesis to $Eg$, we rescale back, and obtain
\begin{equation*}
    \begin{split}
        \|E_{\tau}f\|_{L^p} \leq CC_{p,\epsilon} (\widetilde{K})^{-d+\frac{2d+2n}{p}}
        (\widetilde{K} \delta)^{-100(d+n) \epsilon}  &\Big((\widetilde{K}\delta)^{d-\frac{2d+2n}{p}}(\widetilde{K}\delta)^{-\mathrm{Dec}_p(\mathcal{M}|_{L_{k-2}})}
\\&+
       \sup_{0 \leq m \leq d+n}((\widetilde{K}\delta)^{-\frac{m}{p}+\frac12{X(\mathcal{M},k,m)}}) \Big) \|f\|_{p}.
    \end{split}
\end{equation*}
By plugging all these inequalities to \eqref{0621.316}, the right hand side of \eqref{0621.316} is bounded by
\begin{equation}\label{0715.319}
    \widetilde{C} \Big( A + C_{p,\epsilon}K^{-c\epsilon}( B + C) \Big)\|f\|_p  
\end{equation}
where
\begin{equation}
    \begin{split}
        & A:= K^C 
        \delta^{-\epsilon-(d+n) \epsilon} \max_{0 \leq m \leq d+n}( \delta^{-\frac{m}{p}+\frac12{X_m} }) 
        \\&
B:=\delta^{-100(d+n)\epsilon}\delta^{d-\frac{2d+2n}{p}-\mathrm{Dec}_p(\mathcal{M}|_{L_{k-2}} ) }
\\&
C:=\sup_{0 \leq m \leq d+n}\sup_{K^{c} \leq \widetilde{K} \leq K^{1/2}}
    (\widetilde{K})^{-d+\frac{2d+2n}{p}+\mathrm{Dec}_p(\mathcal{M}|_{L_{k-2}} )}(\widetilde{K}\delta)^{-\frac{m}{p}+\frac12{X(\mathcal{M},k,m)}}.
    \end{split}
\end{equation}
By monotonicity, we have
\begin{equation}
\begin{split}
    C \lesssim (\delta^{-1})^{-d+\frac{2d+2n}{p}+\mathrm{Dec}_p(\mathcal{M}|_{L_{k-2}} )}
    +\delta^{-\frac{m}{p}+\frac12{X(\mathcal{M},k,m)}}.
\end{split}
\end{equation}
Since we have a gain $K^{-c\epsilon}$ in front of $(B+C)$ in \eqref{0715.319}, by taking $K$ sufficiently large, we can close the induction. This finishes the proof of Theorem \ref{restriction}.

\section{A covering lemma for varieties}

The goal of this section is to prove Theorem \ref{0619.thm32}. This follows from a covering lemma for varieties (Theorem \ref{algthm}). One key ingredient of the proof of the lemma is Tarski's projection theorem from real algebraic geometry.
Recall the definition of a semi-algebraic set in Definition \ref{defsemialg}. Another important notion is a regular graph.

\begin{definition}[$(k,\rho,\Lambda,C)$-regular graph]
    Fix an integer $1\le k\le d$, numbers $0<\rho<1$, $C>1$, and $\Lambda \geq 2$. Let $B^k$ be a translation of $[0,\rho]^k$. 
    
    We call the set $\Ga\subset \R^d$ a $(k,\rho,\Lambda,C)$-regular graph over $B^k$ if it has the following form.
    \begin{equation}\label{ga}
        \Ga=\{ (x_{(k)},x_{k+1},\dots,x_d)\in\R^d:x_{(k)}\in B^k, \,x_j=\phi_j(x_{(k)}),\, j=k+1,\dots,d \}.
    \end{equation}
Here, each $\phi_j$ satisfies 
\begin{equation}\label{der}
    |\nabla^m\phi_j(x_{(k)})|\le C \rho^{-m}, \;\; \mathrm{for \; all}\; 1 \leq m \leq \Lambda
\end{equation}
for $x_{(k)}\in B^k$. If after permutation of coordinates, $\Ga$ has the form \eqref{ga}, then we still call $\Ga$ a $(k,\rho,\Lambda,C)$-regular graph. For the purpose of the paper, $\Lambda=2$ is enough, though we still prove for general $\Lambda$ as it may be interesting for other applications. Since $\Lambda$ is not important, we sometimes omit $\Lambda$ and just call $(k,\rho,C)$-regular graph. 

For convenience, we also consider the degenerate case $k=0$. We say $\Gamma \subset \R^d$ is a $(0,\rho,\Lambda,C)$-regular graph if it is a point.
\end{definition}

\begin{remark}
    {\rm
    The four indices $k,\rho,\Lambda,C$ have their respective meanings. $k$ is the dimension of the graph. $\rho$ is the size of the graph. $\Lambda$ is the order of regularity. $C$ is the regularity constant.
    }
\end{remark}

Let us state a basic property of a regular graph. As an application of Taylor's theorem, one can prove the following. We skip the proof.

\begin{lemma}\label{inbox}
    Suppose $\Ga$ is a $(k,\rho,\Lambda,C)$-regular graph over $B^k$, then there exists a cube $B^{d-k}\subset \R^{d-k}$ of side length $Cd\rho$ such that
    \[ \Ga\subset B^k\times B^{d-k}. \]
\end{lemma}

The bottom line of the lemma is that the side length of $B^{d-k}$ is comparable to the side length of $B^k$.

This section is devoted to proving that any semi-algebraic set can be covered by a small neighborhood of several regular graphs. Here is the precise statement.

\begin{theorem}[Covering lemma]\label{algthm}
     Let $0 \leq k \leq d$ be an integer. Let $\Lambda \geq 2$. For $E, A>100$, there exist  positive numbers
     \begin{equation}
       M(d,E,\Lambda),\, C(d,E,\Lambda),\, \mu_0(d,E,A,\Lambda),\, c(d,E,A,\Lambda)  
     \end{equation}
     such that the followings are true.

Let $Z\subset [0,1]^d$ be a semi-algebraic set with complexity $\le E$ and dimension $k$. For any $0<\mu<\mu_0 (d,E,A,\Lambda)$, there exist
\begin{enumerate}[label=(\roman*)]
\item Scales: $\mu_l\in (\mu^{1/(2A)},\mu^{c(d,E,A,\Lambda)})$, $l=1,\dots,M(d,E,\Lambda)$.

\item Boxes: $\{\ci_l\}_{l=1}^{M(d,E,\Lambda)}$ where each $\ci_l=\{I_l\}$ consists of finitely overlapping boxes.\footnote{The number of overlapping is $1000^d$. This overlapping does not make trouble in our application.} Each $I_{l}\in\ci_l$ (after permutation of coordinates) has the form $I_l=I_l^k\times I_l^{n-k}$ where $I_l^k$ is a $k$-dimensional cube of side length $\mu_l$ and $I_l^{n-k}$ is an $(n-k)$-dimensional cube of side length $C(d,E,\Lambda)d\mu_l$.

\item For each box $I_l \in\ci_l$, there is a $(k,\mu_l, \Lambda,C(d,E,\Lambda))$-regular graph $\Ga_{I_l}$ over $I_l^k$ and contained in $I_l$. The graph $\Gamma_I$ can be expressed as
\begin{equation}
    \Gamma_I=I \cap \{x \in \R^d:  Q_I(x)=0 \}
\end{equation}
for some polynomial $Q_I$ of degree $O_{d,E,\Lambda}(1)$.
\end{enumerate}
They satisfy
\begin{equation}\label{cover}
    N_\mu(Z)\cap [0,1]^d\subset \bigcup_{l=1}^{M(d,E,\Lambda)} \bigcup_{I\in\ci_l} N_{(\mu_l)^A}(\Ga_I)\cap I. 
\end{equation}
\end{theorem}

\begin{remark}
{\rm 
There is a classical result in real algebraic geometry saying that any semi-algebraic set $Z$ with complexity $E$ and dimension $k$ is a disjoint union of manifolds diffeomorphic to $(0,1)^l$ (see \cite[Proposition 2.9.10]{BochnakCosteRoy}). More precisely,
\[ Z=\bigcup_{l=0}^k \bigcup_{i=1}^{O_E(1)} \Ga_{l,i}, \]
where each $\Ga_{l,i}$ is a semi-algebraic set diffeomorphic to $(0,1)^l$. However, this is not enough since we need a quantitative version. In our application, we need to cover a small neighborhood of $Z$ by a controlled number of graphs, and we require each graph to have bounded curvature (see \eqref{der}). 
}
\end{remark}

\begin{figure}[ht]
\centering
\begin{minipage}[b]{0.85\linewidth}
\includegraphics[width=11cm]{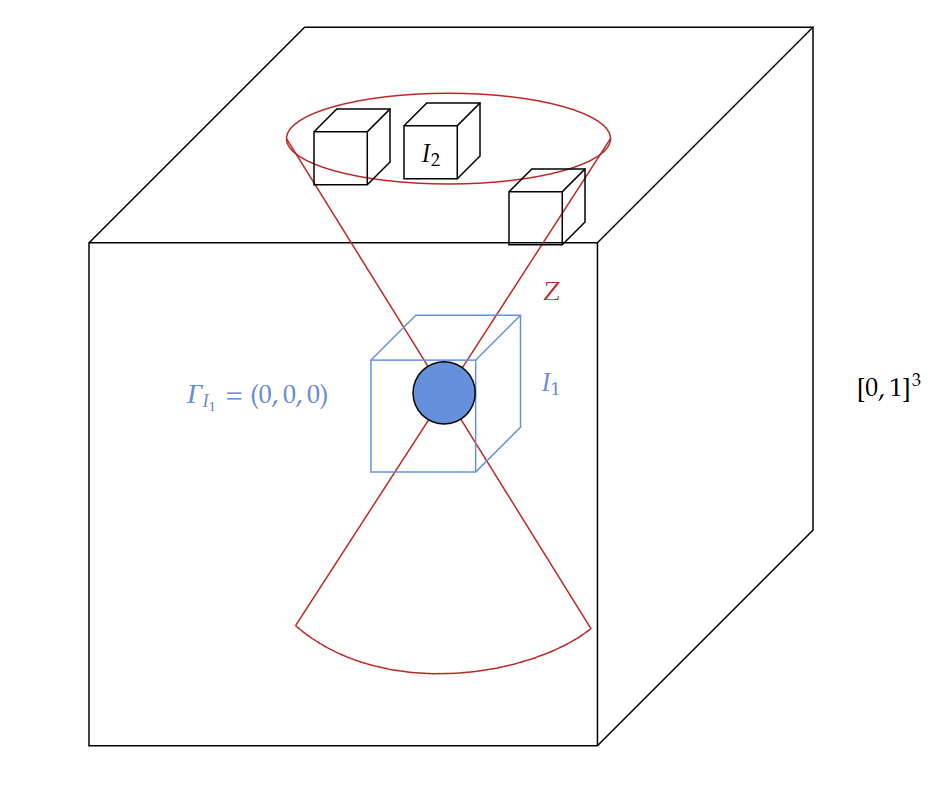}
\vspace{-5mm}
\caption{Cone}
\label{cone}
\end{minipage}
\end{figure}

\begin{figure}[ht]
\centering
\begin{minipage}[b]{0.85\linewidth}
\includegraphics[width=11cm]{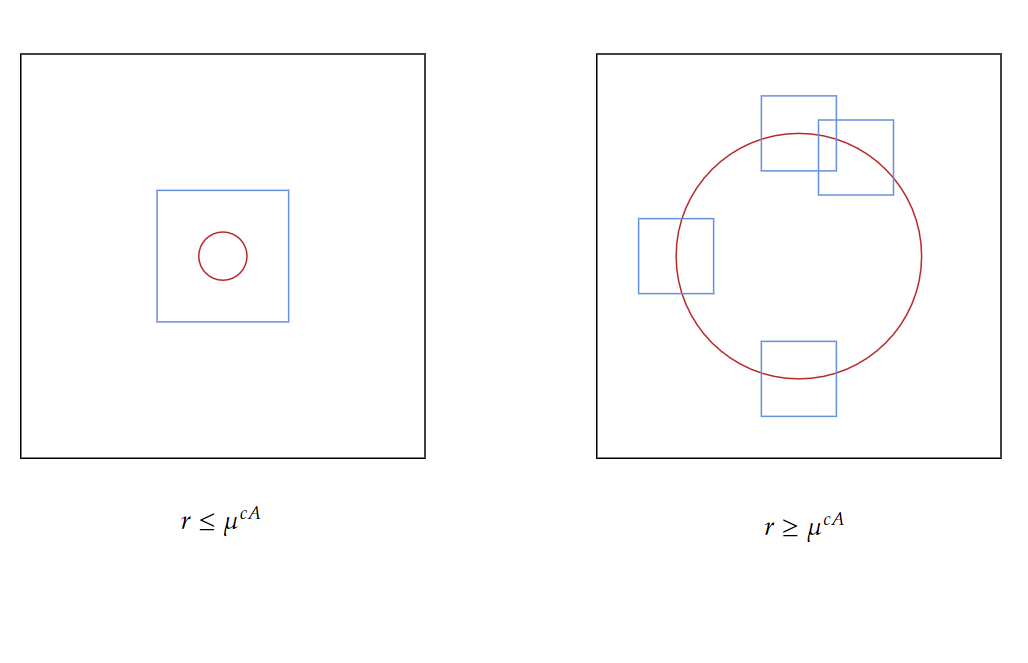} \vspace{-1.2cm}
\caption{Circle}
\label{circle}
\end{minipage}
\end{figure}

\begin{remark}
    {\rm
    Let us mention two examples in the proof of Theorem \ref{algthm}. 
    
    The first example is the cone $Z=\{x^2+y^2=z^2\}\cap [0,1]^3$ (see Figure \ref{cone}). Apparently, the origin  is a singular point. Away from the origin, the variety $Z$ behaves well. Here is how we choose the covering. We let $I_1$ be a box centered at the origin and let $\Ga_{I_1}$ be the origin which is a $0$-dimensional graph. For the remaining portion $Z\setminus N_{\mu_1^A}(\Ga_{I_1})$, we can easily find $\{I_2\}$ and regular graphs $\{\Ga_{I_2}\}$ to cover.

    The second example is the circle $Z_r=\{x^2+y^2=r^2\}\subset [0,1]^2$ (see Figure \ref{circle}). The tricky thing here is that the covering of $N_\mu(Z_r)$ depends on how large $r$ is, though $Z_r$ are topologically the same for all $r$. When $r$ is very small compared with $\mu$ (for example $r=\mu^{1000}$),  $Z_r$ is small and has large curvature. It is impossible to find  $(1,\mu_1,C)$-regular graphs with $\mu_1\in (\mu^{1/(2A)},\mu^c)$ to cover $Z_r$. In this case, we need to view $Z_r$ as a $0$-dimensional set. Then we just simply choose $\mu_1=\mu^{1/(4A)}$, $\Ga_{I_1}=\{(0,0,0) \}$, so we have $N_\mu(Z_r)\subset N_{\mu_1^A}(\Ga_{I_1})$ (see the left picture in Figure \ref{circle}). When $r$ is large compared with $\mu$,  we can cover $Z_r$ by $1$-dimensional regular graphs (see the right picture in Figure \ref{circle}).

    We suggest readers to keep these two examples in mind while reading the proof of Theorem \ref{algthm}. We also encourage interested readers to figure out how these examples are handled in the proof.
    }
\end{remark}

The special case $k=d-1$ of Theorem \ref{algthm} is proved in \cite[Lemma 2.14]{MR4143735}. Our theorem generalizes their theorem to all $k$. The new ingredient is the celebrated result of Tarski. To prove Theorem \ref{algthm}, we use an induction on the dimension $d$. Tarski's theorem enables us to reduce to a lower dimensional problem. As a side note, we  mention that a uniform $l^2$-decoupling inequality for a real variety in $\R^2$ is proved by \cite{li2021decoupling}, which is used as an ingredient of the proof of their main theorem.

\begin{lemma}[Tarski]\label{tarski}
    For any number $E$, there exists a number $\lambda_{d,E}$ so that the following is true. For any semi-algebraic set $Z\subset \R^d$ of complexity $\le E$ and any orthogonal projection $\pi$, the set $\pi(Z)$ is a semi-algebraic set of complexity $\le \lambda_{d,E}$.
\end{lemma}

In Subsection \ref{0707.sub41}, we give a proof of Theorem \ref{0619.thm32} under the assumption of the covering lemma (Theorem \ref{algthm}). In Subsection \ref{0707.sub42}, we give a proof of Theorem \ref{algthm}.

\subsection{Proof of Theorem \ref{0619.thm32} assuming Theorem \ref{algthm}}\label{0707.sub41}
\hfill

We need the following lemma. The case $k=d-1$ is proved in \cite{MR4143735} (see also Lemma 4.4 of \cite{MR3848437}). The general case can be proved identically and we omit the proof.

\begin{lemma}[cf. Theorem 2.5 of \cite{MR4143735}]\label{0713.lem47} Let  $B,k>0$. Suppose that $\mathcal{H}$ is a $k$-dimensional manifold, which is written as a graph of form
\begin{equation}
    (t_{k+1},\ldots,t_d)=F(t') \;\; \text{ with } \;\; \|F\|_{C^2} \leq B.
\end{equation}
Then for any $K$ and $\epsilon'>0$, we have
    \begin{equation}
        \Big\|\sum_{\tau \cap \mathcal{H} \neq \emptyset: \tau \in \mathcal{P}(K^{-1}) } E^{\mathcal{M}}f_{\tau} \Big\|_{L^p} \leq C_{\epsilon',B}K^{\epsilon'+\mathrm{Dec}_p(\mathcal{M}|_{L_k})
 }  \Big( \sum_{\tau \in \mathcal{P}(K^{-1}) } \|E^{\mathcal{M}}f_{\tau}\|_{L^p}^p \Big)^{\frac1p}.
    \end{equation}
The constant $C_{\epsilon',B}$ depends on $\epsilon'$ and $B$ but independent of the choice of $F$.
\end{lemma}

Let $Z \subset \R^d$ be a $k$-dimensional semi-algebraic set with complexity $E$.
We apply Theorem \ref{algthm} to $Z$ with $\mu=K^{-1}$, $A=(100\epsilon)^{-1}$, and $\Lambda=2$. Then we have
\begin{equation}
    N_\mu(Z)\cap [0,1]^n\subset \bigcup_{l=1}^M \bigcup_{I\in\ci_l} {N}_{(\mu_l)^A}(\Ga_I)\cap I. 
\end{equation}
Without loss of generality, we may assume that
\begin{equation}
    \mu_1 \geq \mu_2 \geq \cdots \geq \mu_M.
\end{equation}
Take $K_i:=(\mu_i)^{-A}$. By the bounds of $\mu_l$, we have \eqref{0729.31}. Cover $\bigcup_{I \in \ci_I} (N_{\mu_l^A}(\Ga_I)\cap I)$ by dyadic cubes $\tau_i$ with side length $K_i^{-1}$, and denote by $\mathcal{P}_i(K_i^{-1})$ the collection of $\tau_i$'s. Note that $\mathcal{P}_i(K_i^{-1}) \subset \mathcal{P}(K_i^{-1})$. We have
\begin{equation}
   {N}_{K^{-1}}(Z) \cap [0,1]^d \subset \bigcup_{i=1}^{M} \bigcup_{\tau_i \in \mathcal{P}_i(K_i^{-1}) } \tau_i. 
\end{equation}
It remains to prove
\begin{equation}\label{0713.411}
\begin{split}
    \Big\|\sum_{\substack{ \tau_i \in \mathcal{P}_i(K_i^{-1}) }} E^{\mathcal{M}}_{\tau_i}f \Big\|_{L^p} \leq
    C_{\epsilon,p,E} K_i^{\epsilon+\mathrm{Dec}_p(\mathcal{M}|_{L_{k}})}  \Big( \sum_{\tau_i \in \mathcal{P}(K_i^{-1}) } \|E^{\mathcal{M}}_{\tau_i}f\|_{L^p}^p \Big)^{\frac1p}.
    \end{split}
\end{equation}
Note that the cardinality of $\ci_I$ is  $O(K_i^{\frac1A})$. Hence, by the triangle inequality, for some $I$, we have
\begin{equation}
\begin{split}
    \Big\|\sum_{\substack{ \tau_i \in \mathcal{P}_i(K_i^{-1}) }} E^{\mathcal{M}}_{\tau_i}f \Big\|_{L^p}  \lesssim (K_i)^{\frac1A}
     \Big\|\sum_{\substack{ \tau_i \in \mathcal{P}(K_i^{-1}): \\ \tau_i \cap (\Gamma_I \cap I) \neq \emptyset } } E^{\mathcal{M}}_{\tau_i}f \Big\|_{L^p}. 
    \end{split}
\end{equation}
Recall that by Theorem \ref{algthm} the graph $\Ga_{I}$ is $(k,K_i^{-\frac1A},C(d,E))$-regular  over $I^k$ and contained in $I=I^k \times I^{d-k}$. The side length of $I$ is $O(K_i^{-1/A})$. For convenience, assume that $I=[0,K_i^{-\frac{1}{A}}]^d$. We do rescaling; take a linear transform $L$ so that $L(I)=[0,1]^d$. Define the function $g$ so that $g(\xi):=f(L^{-1}\xi)$. Then we have
\begin{equation}
    \|E_{I}^{\mathcal{M}}f\|_p \sim ({K_i}^{\frac{1}{A}})^{-d+\frac{d+2n}{p}} \|E^{\mathcal{M}}_{[0,1]^d}g\|_{p}, \;\;\; \|f\|_{L^p(\tau)} \sim ({K_i}^{\frac{1}{A}})^{-\frac{d}{p}}\|g\|_p.
\end{equation} 
After this rescaling, $\Gamma_I$ becomes
\begin{equation}
    \{(x_{(k)},x_{k+1},\ldots,x_d) \in \R^d: x_j=\widetilde{\phi}_j(x_{(k)}),\, j=k+1,\dots,d \}
\end{equation}
where $\widetilde{\phi}_j$ satisfies
\begin{equation}
     |\nabla\phi_j(x_{(k)})|\le C,\ \ |\nabla^2\phi_j(x_{(k)})|\le C. 
\end{equation}
Apply Lemma \ref{0713.lem47} to $Eg$ and rescale back, and we obtain
\begin{equation}
    \Big\|\sum_{\substack{ \tau_i \in \mathcal{P}(K_i^{-1}): \\ \tau_i \cap (\Gamma_I \cap I) \neq \emptyset } } E^{\mathcal{M}}_{\tau_i}f \Big\|_{L^p} \lesssim 
K_i^{\frac{\epsilon}{100}+\mathrm{Dec}_p(\mathcal{M}|_{L_{k}})}  \Big( \sum_{\tau_i \in \mathcal{P}(K_i^{-1}) } \|E^{\mathcal{M}}_{\tau_i}f\|_{L^p}^p \Big)^{\frac1p}.
\end{equation}
This gives \eqref{0713.411}, and completes the proof of Theorem \ref{0619.thm32}.

\subsection{Proof of Theorem \ref{algthm}}\label{0707.sub42}

Fix $\Lambda \geq 2$. Throughout the proof, $\Lambda$ will not change. So we do not keep track of dependence on $\Lambda$ and skip mentioning $\Lambda$. For example, $(k,\mu_l,\Lambda,C(d,E,\Lambda))$-regular graph will be called $(k,\mu_l,C(d,E))$-regular graph.

Let $Z \subset [0,1]^d$ be a semi-algebraic set with complexity $\leq E$ and dimension $k$. By the definition of a semi-algebraic set, $Z$ is a union of the form
\[ \{x\in \R^d: P_1(x)=0,\dots, P_l(x)=0, P_{l+1}(x)>0,\dots, P_{l+l'}(x)>0  \}. \]
It is convenient to assume that each polynomial in the representation of $Z$ is normalized in the following sense. 

\begin{definition}
    For a polynomial $P$, we define $\|P\|$ to be the $\ell^1$ sum of the coefficients of $P$. We say that $P$ is nomalized if $\|P\|=1$.
\end{definition}

Let us explain an outline of the proof of Theorem \ref{algthm}.
We will use an induction on the dimension $d$. The case $d=1$ is the base case. In Step 1, we verify Theorem \ref{algthm} for the base case. In Step 2, we prove a covering lemma for the sublevel set of a single polynomial (Proposition \ref{0630.prop49}). In Step 3, we use Proposition \ref{0630.prop49} to obtain some type of variety covering lemma, which has a certain inductive structure. In Step 4, we apply the induction hypothesis on the dimension $d$. In Step 5, we choose all the parameters and close the induction. 

\subsection*{Step 1. Base case of the induction on \texorpdfstring{$d$}{}}\label{0707.step1}
\hfill

The case $d=1$ will be the base case of the induction. Let us assume that $Z$ is nonempty. 
If $k=1$, then we just choose $\mu_1=\mu^{\frac{1}{3A}}$ and $\ci_1=\{I_1\}$ to be a collection of $\mu_1$-intervals that cover $[0,1]$. Since $d=1$, each $\Ga_{I_1}=I_1$ is trivially a $(1,\mu_1,C)$-regular graph. 


If $k=0$, then there must exist a polynomial $P_1$ so that
\begin{equation}
    Z \subset \{x \in \R: P_1(x)=0 \}.
\end{equation}
The degree of $P_1$ is controlled by a constant depending only on the complexity of $Z$.
By the fundamental theorem of algebra, the right hand side is a union of at most $C(E)$   many points for some constant $C(E)$ depending on $E$. So one can check that the theorem is true. We have verified that Theorem \ref{algthm} is true for the case $d=1$.
\\

\subsection*{Step 2. Covering lemma for the set \texorpdfstring{$\{|P_1| < 1/K \}$}{}}
\hfill

We now assume that $d \geq 2$ and Theorem \ref{algthm} is true for $\leq d-1$. If $k=d$, then $Z$ is a union of sets of the form
\begin{equation}
 Z=\{x\in \R^d:   P_{1}(x)>0,\dots, P_{s}(x)>0  \}   
\end{equation}
and this set is an open set. We deal with it as in the base case $d=k=1$. We choose $\ci_1$ to be a collection of $\mu_1$-boxes which form a finitely overlapping covering of $[0,1]^d$. Each $I_1$ is also $\Ga_{I_1}$ since $d=k$. Hence, we find a covering of $Z\cap [0,1]^d$. Next, we consider $k \leq d-1$. By the definition, the Hausdorff dimension of $Z$ is smaller or equal to $d-1$, so there must exist a polynomial $P_1$ such that
\begin{equation}\label{0630.414}
    Z \subset \{x \in \R^d: P_1(x)=0 \}.
\end{equation}
Since $P_1$ is normalized, for any $K$, one can see that 
\begin{equation}\label{0630.415}
    N_{K^{-1}}(Z) \subset \{x \in \R^d: |P_1(x)| \leq K^{-1} \}.
\end{equation}
    Take $K$ so that $\mu=:K^{-1}$. Then the left hand side of \eqref{0630.415} is that for \eqref{cover}.
For simplicity, we use the notation $Z_{P}:=\{x: P(x)=0\}$ for any polynomial $P$.
 We will first prove a structure lemma for $Z_{P_1}$ (Proposition \ref{0630.prop49}). Later, this will be used to find a structure lemma for $Z$ by using \eqref{0630.415}.

\begin{proposition}\label{0630.prop49}
For every natural numbers $d, A'\ge 1$, every $D\ge 1$, and every sufficiently large $K>1$ (depending on $d, D$ and $A'$), there exist dyadic numbers
    \[ K_1\le K_2\le \cdots\le K_{D+1}=K\textup{~with~}K_{j+1}\sim_{d,j}K_j^{A'+1}, \]
    and the numbers $C_{d,D}$ and $C_{d,D,\Lambda}$
    satisfying the following: for every normalized polynomial $P$ of degree $D$ in $d$ variables with real coefficients,
    \begin{enumerate} \item
     we have
    \begin{equation}
    \big\{ x \in \R^d: |P(x)| < \frac{1}{K}\big\}\cap [0,1]^d\subset \bigcup_{j=1}^D \bigcup_{B\in\cB_j} \Big(N_{K_j^{-A'}}(\Psi_B)\cap B\Big). 
\end{equation} 
\item
 $\cB_j=\{B\}$ is a collection of finitely overlapping boxes $B$.\footnote{{The number of overlapping is at most $1000^d$.} } Each box $B$ has dimensions
\begin{equation}
    (C_{d,D}K_j)^{-1} \times \cdots \times (C_{d,D}K_j)^{-1} \times dC_{d,D,\Lambda}(C_{d,D}K_j)^{-1}.
\end{equation}
\item
 $\Psi_B$ is an $(d-1,(C_{d,D}K_j)^{-1},C_{d,D,\Lambda})$-regular graph. Moreover, we have
\begin{equation}
    \Psi_B=B \cap Z_{Q_B}
\end{equation}
 for some polynomial $Q_B$ of degree $\le D$.
\end{enumerate}
\end{proposition}

To prove Proposition \ref{0630.prop49}, we use the following lemma from \cite[Lemma 2.14]{MR4143735} as a black box.

\begin{lemma}\label{gzlem}
    For every natural numbers $d, A', D \ge 1$, and every sufficiently large $K \geq 1$ (depending on $d, D$ and $A'$), there exist dyadic numbers
    \[ K_1\le K_2\le \cdots\le K_{D+1}:=K\textup{~with~}K_{j+1}\sim_{d,j}K_j^{A'+1} \]
    such that, for every normalized polynomial $P$ of degree $D$ in $d$ variables with real coefficients, there exists an increasing sequence of multi-indices $\alpha_D<\alpha_{D-1}\cdots<\alpha_1$ with $|\alpha_j|=D-j$ such that
    \begin{equation*}
    \big\{x \in \R^d :|P(x)|< \frac1K \big\}\cap [0,1]^d\subset \bigcup_{j=1}^D N_{K_j^{-A'}}\Big( Z_{\partial^{\alpha_j}P}\cap \{x \in \R^d:|\nabla\partial^{\alpha_j} P(x)|\ge \frac{1}{K_j}\} \Big).
    \end{equation*}  
\end{lemma}

\begin{proof}[Proof of Proposition \ref{0630.prop49}]
We apply Lemma \ref{gzlem} to $P$, and obtain
\begin{equation}\label{coverzp}
    \begin{split}
\{|P|<\frac1K\}\cap [0,1]^d \subset \bigcup_{j=1}^D N_{K_j^{-A'}}\Big( Z_{\partial^{\alpha_j}P}\cap \{|\nabla\partial^{\alpha_j} P|\ge \frac{1}{K_j}\} \Big).
    \end{split}
\end{equation}  
By pigeonholing, the right hand side of \eqref{coverzp} can be rewritten as
\begin{equation}\label{0703.421}
    \bigcup_{j=1}^D \bigcup_{i=1}^d N_{K_j^{-A'}}\Big( Z_{\partial^{\alpha_j}P}\cap \{|\nabla\partial^{\alpha_j} P|\ge \frac{1}{K_j}\} \cap \big\{ |\partial_{i}\partial^{\alpha_j} P|\ge\frac{1}{\sqrt{d}}|\nabla\partial^{\alpha_j} P| \big\}\Big).
\end{equation}
For simplicity, we write $Q:=\partial^{\alpha_j}P$. Next, we study the set
\[ Z_{Q}\cap \big\{|\nabla Q|\ge \frac{1}{K_j}\big\} \cap \big\{ |\partial_{i}Q|\ge\frac{1}{\sqrt{d}}|\nabla Q| \big\}. \]
By permutation of the coordinates, we may assume $i=d$, and this set becomes
\begin{equation}\label{thisset}
   \mc{Z}_Q:= Z_{Q}\cap \big\{|\nabla Q|\ge \frac{1}{K_j}\big\} \cap \big\{ |\partial_{d}Q|\ge\frac{1}{\sqrt{d}}|\nabla Q| \big\}. 
\end{equation} 
For $x\in\R^d$, we write $x=(y,x_d)$ where $y\in\R^{n-1}$ and $x_d\in\R$.
If $x^*=(y^*,x^*_d)$ lies in the set \eqref{thisset}, then
\[ Q(y^*,x^*_d)=0. \]
Since $|\partial_{d}Q(y^*,x^*_d)|\ge \frac{1}{\sqrt{d}}|\nabla Q(y^*,x^*_d)|\ge \frac{1}{\sqrt{d}K_j}$, by the quantitative version of the implicit function theorem, there exist a number $C_{d,D}$ and a function $\psi(y)$ defined in the $\frac{1}{C_{d,D}K_j}$-neighborhood of $y^*$ such that $\psi(y^*)=x_n^*$ and 
\begin{equation}\label{takeder}
    Q(y,\psi(y))=0, \textup{~for~} |y-y^*|\le \frac{1}{C_{d,D}K_j}.
\end{equation}
Denote by $B^{d-1} \subset \R^{d-1}$ the cube centered at $y^*$ of  side length $\frac{1}{C_{d,D}K_j}$.\footnote{
We remark that the existence of such $\frac{1}{C_{d,D}K_j}$-neighborhood of $y^*$ is by the $C^2$ boundedness of $Q$: $\|Q\|_{C^2([0,1]^d)}\lesssim_{d,D} 1$.} Note that the constant $C_{d,D}$ is independent of the choice of the point $x^{*}$. Let us  assume that $C_{d,D}$ is sufficiently large. From the lower bound of  $|\nabla Q|$, we can prove a quantitative $C^2$ bound of the function $\psi$. Let us state it as a lemma.

\begin{lemma}\label{0630.lem411} If $C_{d,D}$ is sufficeintly large, then
   for any $y \in B^{d-1}$ we have
    \begin{equation}
        \begin{split}
            &\label{mainest}
            |\nabla\psi(y)|\le 2d,
            \\&|\partial_d Q(y,\psi(y))|\ge \frac{1}{2\sqrt{d}K_j}, 
            \\&|\p_{il}\psi(y)| \leq 1000^{d+D^2} K_j.
        \end{split}
    \end{equation}
    Moreover, for a multi-index $\alpha$ with $|\alpha| \leq \Lambda$, we have
    \begin{equation}\label{0717.430}
        |\partial^{\al}\psi(y)| \leq C_{d,D,\Lambda}(K_j)^{|\alpha|-1}
    \end{equation}
    for some constant $C_{d,D,\Lambda}$.
\end{lemma}

The bounds \eqref{mainest} may not be sharp, but they are good enough for our purpose.

\begin{figure}[ht]
\centering
\begin{minipage}[b]{0.85\linewidth}
\includegraphics[width=11cm]{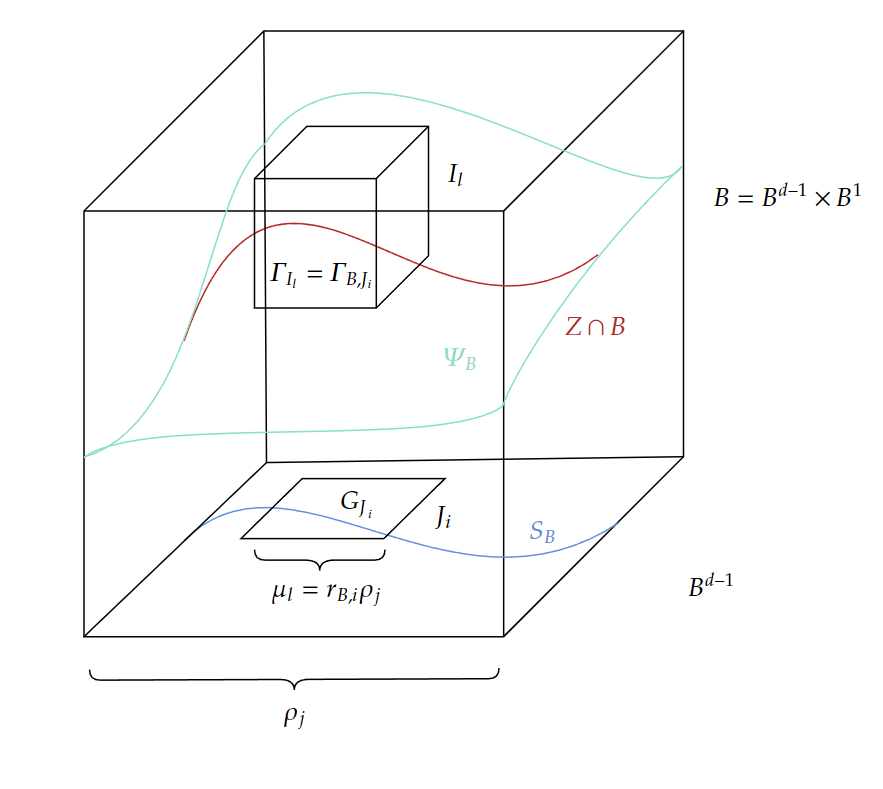} \vspace{-0.7cm}
\caption{Regular graph}
\label{fig2}
\end{minipage}
\end{figure}

Let us assume the lemma and continue proving Proposition \ref{0630.prop49}.  We claim the following: 
for each $x\in 
\mc{Z}_Q$ (see \eqref{thisset}), there exist 
\begin{enumerate}
    \item a box $B=B^{d-1}\times B^1$ centered at $x$ of dimensions 
    \begin{equation*}
       (C_{d,D}K_j)^{-1} \times \cdots \times (C_{d,D}K_j)^{-1} \times dC_{d,D,\Lambda} (C_{d,D}K_j)^{-1}.
    \end{equation*}
    \item an $(d-1,(C_{d,D}K_j)^{-1},C_{d,D,\Lambda})$-regular graph $\Psi_B$ over $B^{d-1}$ satisfying
    \begin{equation*}
        \Psi_B \subset B \;\; \mathrm{ and } \;\; B\cap Z_{Q}=\Psi_B.
    \end{equation*}
\end{enumerate}
Let us give a proof of the claim.
Define the graph
\[\Psi_B:=\{ (y,\psi(y)): y\in B^{d-1} \}\]
associated with the cube $B^{d-1}$.
By Lemma \ref{0630.lem411},  $\Psi_B$ is a $(d-1,\frac{1}{C_{d,D}K_j},C_{d,D,\Lambda})$-regular graph over $B^{d-1}$.    By lemma \ref{inbox}, there exists an interval $B^1\subset \R$ of length $\frac{dC_{d,D,\Lambda}}{C_{d,D}K_j}$ so that $\Psi_B$ is contained in the box $B^{d-1}\times B^1$. Define $B:=B^{d-1} \times B^1$.
We next show that 
\begin{equation}
    B \cap Z_Q = \Psi_B
\end{equation}
provided that $C_{d,D}$ is sufficiently large. Note that 
\begin{equation}
 |\partial_d Q(y^*,x_d^*)|\ge \frac{1}{{d}K_j}, \;\;\;    |\nabla^2 Q(x)| \leq 1000^{d+D^2}
\end{equation}
for $x\in [0,1]^d$. Hence, by Taylor's expansion and taking $C_{d,D}$ sufficiently large, we have $|\partial_d Q(x)|>0$ for any $x \in B$. This means that for any $y\in B^{d-1}$ there is at most one $x_d$ such that $Q(y,x_d)=0$. The above discussion gives that $x_d=\psi(y)$. This completes the proof of the claim. (See Figure \ref{fig2}, the green graph $\Psi_B$ is over $B^{d-1}$ inside the box $B$.)
\medskip

We are ready to finish the proof of Proposition \ref{0630.prop49}. By choosing $\{x\}$ to be a maximal $\frac{1}{10}(C_{d,D}K_j)^{-1}$-seperated subset of $\mc{Z}_Q$ (see \eqref{thisset}), and by the claim, we can find the corresponding set of finitely overlapping boxes $\cB_j=\{B\}$  and the associated graphs $\{\Psi_B\}_{B\in\cB_j}$ satisfying
\[ N_{K_j^{-A'}}(\mc{Z}_Q ) \subset \bigcup_{B\in\cB_j} N_{K_j^{-A'}}(\Psi_B)\cap B. \]
Plugging back to  \eqref{coverzp} (recall \eqref{thisset} and \eqref{0703.421}) yields 
\begin{equation}\label{structure1}
    \{|P_1| < 1/K\}\cap [0,1]^d\subset \bigcup_{j=1}^D \bigcup_{B\in\cB_j} N_{K_j^{-A'}}(\Psi_B)\cap B. 
\end{equation} 
This finishes the proof of Proposition \ref{0630.prop49}.
\end{proof}

\begin{proof}[Proof of Lemma \ref{0630.lem411}]

Let us first show that
\begin{equation}\label{0703.428}
    |\nabla\psi(y)|\le 2d \textup{~and~} |\partial_d Q(y,\psi(y))|\ge \frac{1}{2\sqrt{d}K_j}
\end{equation}
for $ |y-y^*|\le \frac{1}{C_{d,D}K_j}$.
Taking $i$-th derivative $(1\le i\le d-1)$ to \eqref{takeder}, we get
\[\partial_i Q+ \partial_d Q\cdot \partial_i\psi(y^*)=0.\]
So we have
\begin{equation}\label{boundpsi}
   |\nabla\psi(y^*)|\le \sqrt{d}\frac{|\nabla Q(y^*,x_d^*)|}{|\partial_d Q(y^*,x_d^*)|}\le d.  
\end{equation}  
Fix $1\le i\le d-1$. Let $e_i=(0,\dots,0,1,0,\dots,0)$ be a vector of $\R^{d-1}$ whose $i$-th coordinate is one.

We let $r$ be the largest number such that 
\[|\partial_i\psi(y^*+te_i)|\le 2\sqrt{d}, \textup{~for~} t\in [0,r].\]
We will show that $r\ge \frac{1}{C_{d,D} K_j}$. Note that this lower bound of $r$ gives the first inequality of \eqref{0703.428}.
By the definition of $r$, 
\begin{equation}\label{contr}
    |\partial_i\psi(y^*+re_i)|=2\sqrt{d}.
\end{equation} 
Similar to \eqref{boundpsi}, we have
\begin{equation}\label{plug1}
    |\partial_i\psi(y^*+re_i)|\le \frac{|\nabla Q(y^*+re_i,\psi(y^*+re_i))|}{|\partial_d Q(y^*+re_i,\psi(y^*+re_i))|}. 
\end{equation} 
Note that 
\begin{align*}
    &|\nabla Q(y^*+re_i,\psi(y^*+re_i))-\nabla Q(y^*,\psi(y^*))|\\
    &\lesssim \|Q\|_{C^2([0,1]^d)}\bigg( r+r\sup_{t\in[0,r]}|\partial_i\psi(y^*+te_i)| \bigg)\\
    &\leq 1000^{d+D} r.
\end{align*}
Hence, if $r\le \frac{1}{C_{d,D}K_j}$ where $C_{d,D}$ is large enough, then
\[ |\nabla Q(y^*+re_i,\psi(y^*+re_i))|\le |\nabla Q(y^*,\psi(y^*))|+\frac{1}{100dK_j}. \]
By the same reasoning, if $r\le \frac{1}{C_{d,D} K_j}$, then
\begin{equation}\label{lowqn}
    |\partial_d Q(y^*+re_i,\psi(y^*+re_i))|\ge |\partial_d Q(y^*,\psi(y^*))|-\frac{1}{100dK_j}. 
\end{equation} 
Plugging into \eqref{plug1}, we obtain 
\[ |\partial_i\psi(y^*+re_i)|\le \frac{|\nabla Q(y^*,\psi(y^*))|+\frac{1}{100dK_j}}{|\partial_d Q(y^*,\psi(y^*))|-\frac{1}{100dK_j}}.  \]
Recalling that $|\nabla Q(y^*,\psi(y^*))|\ge (K_j)^{-1}$ and $|\partial_d Q(y^*,\psi(y^*))|\ge \frac{1}{\sqrt{d}}|\nabla Q(y^*,\psi(y^*))|$, we have
\begin{equation}
    |\partial_d Q(y^*,\psi(y^*))|\ge \frac{1}{\sqrt{d} K_j}.
\end{equation}
Hence, we finally
obtain 
\[ |\partial_i\psi(y^*+re_i)|\le \frac{\sqrt{d}+\frac{\sqrt{d}K_j}{100dK_j}}{1-\frac{\sqrt{d}K_j }{100dK_j}}< 2\sqrt{d},  \]
which contradicts \eqref{contr}. Also, \eqref{lowqn} give a lower bound on $|\partial_d Q(y,\psi(y))|$.
We have finished the proof of \eqref{0703.428}.

Next, we will derive a $C^2$ bound for $\psi$. Taking $\partial_i\partial_l$ to \eqref{takeder}, we obtain
\[ \p_i\p_l Q+\p_d\p_l Q\cdot \p_i\psi+\p_i\p_d Q\cdot \p_l\psi+\p_d^2 Q\cdot \p_i\psi\cdot \p_l\psi+\p_d Q\cdot \p_i\p_l\psi=0. \]
The $C^2$ boundedness of $Q$ and \eqref{0703.428}
imply
\begin{equation}
    |\p_i\p_l\psi(y)| \leq 1000^{d+D^2} K_j, \textup{~for~} |y-y^*|\le \frac{1}{C_{d,D}K_j}.
\end{equation}

Lastly, let us show \eqref{0717.430}. We use an induction on $\Lambda=|\alpha|$. The base case $\Lambda \leq 2$ is proved by \eqref{mainest}. We may assume that \eqref{0717.430} is true for $|\alpha|<\Lambda$ and let us prove \eqref{0717.430} for $|\alpha|=\Lambda$. Write $\p^\al=\p_{i_\Lambda}\dots\p_{i_1}$, where $1\le i_j\le d$ for $j=1,\dots,\Lambda$. We will apply $\p^\al$ to the equation $Q(y,\psi(y))=0$.
First applying $\p_{i_1}$ to the equation $Q(y,\psi(y))=0$, we get
\[ \p_{i_1}Q+\p_d Q\cdot \p_{i_1}\psi=0. \]
We iteratively apply $\p_{i_j}$ to the equation. Then $\p_{i_j}$ will hit either on $\p^\beta Q$ or $\p^\ga\psi$.
By the chain rule, it is not hard to see that we obtain
\begin{equation}\label{440}
    \partial_d Q \cdot \partial^{\alpha} \psi + \sum_{m}\sum_{ \gamma_1,\dots,\ga_m}\sum_{|\beta| \leq \Lambda} a_{\beta,\{\gamma_i\}} \partial^{\beta}Q \cdot \prod_{i=1}^m (\partial^{\gamma_i} \psi) =0.
\end{equation}
The first term on the left hand side of \eqref{440} is the main term where all the $\p_{i_j}$ hit on $\psi$. For the second term, $a_{\beta,\{\ga_i\}}$ are constants, and each $\ga_i\in \N^d$. The summation rule is: when $m=1$, then $|\ga_1|<\Lambda$; when $m\ge 2$, then $\sum_{i=1}^m |\ga_i|\le \Lambda$.  This gives the inequality
\begin{equation}
    |\partial^{\alpha}\psi| \leq \frac{1}{|\partial_d Q|} \Big(\sum_{m}\sum_{ \gamma_1,\dots,\ga_m}\sum_{|\beta| \leq \Lambda} a_{\beta,\{\gamma_i\}} |\partial^{\beta}Q| \cdot \prod_{i=1}^m |\partial^{\gamma_i} \psi| \Big)
    .
\end{equation}
It suffices to use the lower bound of $\partial_d Q$ (see \eqref{mainest}) and induction hypothesis.
This finishes the proof of Lemma \ref{0630.lem411}.
\end{proof}

We have proved the structure lemma for $\{|P_1| < K^{-1} \}$. Let $A'$ be a large number to be determined later\footnote{We will later see $A'\sim_{d,E}  c(d-1,\lambda_{d,E},A^2)^{-1}A^{2}$}. Let $D:=\deg P_1$. Take $K:=\mu^{-1}$. Recall \eqref{0630.415}:
\begin{equation}\label{0702.433}
    N_{\mu}(Z) \subset \{x \in \R^d: |P_1(x)| \leq K^{-1} \}.
\end{equation}
In the next step, we apply Proposition \ref{0630.prop49} and find a covering of $Z$.
\\

\subsection*{Step 3. Inductive structure}
\hfill

We apply Proposition \ref{0630.prop49} to the right hand side of \eqref{0702.433}. Then we have 
\begin{equation}
    N_{\mu}(Z)\cap [0,1]^d\subset \bigcup_{j=1}^D \bigcup_{B\in\cB_j} \Big(N_{K_j^{-A'}}(Z_{Q_B} \cap B)\cap B\Big). 
\end{equation} 
This implies 
\begin{align*}
    N_{\mu}(Z)\cap [0,1]^d &\subset \bigcup_{j=1}^D \bigcup_{B\in\cB_j} \Big( N_{K_j^{-A'}}\big((N_{K_j^{-A'}+\mu}Z) \cap Z_{Q_B}\cap B \big)\cap B \Big)\\
    &\subset \bigcup_{j=1}^D \bigcup_{B\in\cB_j} \Big( N_{K_j^{-A'}}\big((N_{2K_j^{-A'}}Z) \cap Z_{Q_B}\cap B \big)\cap B \Big).
\end{align*}
The last inclusion is because of $\mu=K^{-1}\le K_j^{-A'}$. For simplicity, we may slightly modify $K_j$ to get rid of the factor $2$ in $2K_j^{-A'}$, which is harmless. We have  
\begin{equation}\label{contain}
    N_{\mu}(Z)\cap [0,1]^d\subset \bigcup_{j=1}^D \bigcup_{B\in\cB_j} \Big( N_{K_j^{-A'}}\big((N_{K_j^{-A'}}Z) \cap Z_{Q_B}\cap B \big)\cap B \Big).
\end{equation}

We will observe an inductive structure.
Fix $B\in\cB_j$. After permutation of coordinates, we may write $B=B^{d-1}\times B^1$ and $\Psi_B=\{ (y, \psi_B(y)): y\in B^{d-1} \}$. Recall that
\begin{equation}\label{0713.439}
    \Psi_B=Z_{Q_B} \cap B.
\end{equation}
Let $\pi_{B^{d-1}}$ be the orthogonal projection to the hyperplane spanned by $B^{d-1}$. 
Note that 
\begin{equation}
    (N_{K_j^{-A'}}Z) \cap Z_{Q_B}\cap B \subset \big\{ (y,\psi_B(y)): y\in \pi_{B^{d-1}}((N_{K_j^{-A'}}Z) \cap B) \big\}.
\end{equation}
Define the set (see Figure \ref{fig2}.)
\begin{equation}
    S_B:= \pi_{B^{d-1}}(Z \cap B).
\end{equation}
Notice that
\begin{equation}\label{0704.439}
    (N_{K_j^{-A'}}Z) \cap Z_{Q_B}\cap B \subset \{ (y,\psi_B(y)): y\in N_{K_j^{-A'}}(S_B) \}.
\end{equation}
By Tarski's theorem (Lemma \ref{tarski}), $S_B$ is a semi-algebraic set of complexity $\le \lambda_{d,E}$. This set is contained in $B^{d-1}$. We remark that the Hausdorff dimension of the set $S_B$ is less or equal to $k$. Noting that $B^{d-1}$ is an $(d-1)$-dimensional box lying on an $(d-1)$-dimensional affine subspace of $\R^d$, we can use induction to find a covering of $N_{K_j^{-A'}}(S_B)$ using regular graphs.
\\

\subsection*{Step 4. Applying the induction hypothesis}
\hfill

We will use induction hypothesis to cover the right hand side of $\eqref{contain}$. Fix $B\in\cB_j$. For simplicity, we set $\rho_j=\frac{1}{C_{d,D}K_j}$. Recall that the dimensions of $B$ is  
\begin{equation}
    \rho_j \times \cdots \times \rho_j \times dC_{d,D,\Lambda}\rho_j.
\end{equation}
By translation, we may assume $B=[0,\rho_j]^{d-1}\times [0,dC_{d,D,\Lambda}\rho_j]$. We also recall that 
\begin{equation}
    (N_{K_j^{-A'}}Z) \cap Z_{Q_B}\cap B \subset \{ (y,\psi_B(y)): y\in N_{K_j^{-A'}}(S_B) \}.
\end{equation}
By the definition of a regular graph,  $\psi_B:[0,\rho_j]^{d-1}\rightarrow [0,dC_{d,D,\Lambda}\rho_j]$ satisfies 
\begin{equation}\label{derbd}
    |\nabla^j \psi_B|\le C_{d,D,\Lambda} \rho_j^{-j+1}, \;\;\; j=1,\ldots, \Lambda. 
\end{equation} 
Note that after the rescaling $y\mapsto y\rho_j^{-1}$, $[0,\rho_j]^{d-1}$ becomes $[0,1]^{d-1}$, and $S_B$  becomes another semi-algebraic set $\wt S_B$ of complexity $\le \lambda_{d,E}$ inside $[0,1]^{d-1}$. The Hausdorff dimension of $\widetilde{S}_B$ is less or equal to $k$. Hence, we can use the induction hypothesis of Theorem \ref{algthm} in $[0,1]^{d-1}$ to $\wt S_B$ with $A$ replaced by $A^2$ and then rescale back. As a consequence, we obtain the following result. 
\medskip

There exist constants\footnote{
$\wt c=c(d-1,\lambda_{d,E},A^2), \wt \mu_0=\mu_0(d-1,\lambda_{d,E},A^2), \wt M= M(d-1,\lambda_{d,E})$ and $\wt C=C(d-1,\lambda_{d,E})$} $\wt c, \wt \mu_0, \wt M$, and $\wt C$ so that the following is true: for any $0<r_B<\wt \mu_0$, there exist
\begin{enumerate}[label=(\roman*)]
\item Scales: $r_{B,i}\in (r_B^{1/(2A^2)},r_B^{\wt c}), i=1,\dots, \wt M.$
\item Boxes: $\{\cj_{B,i}\}_{i=1}^{\wt M}$ where each $\cj_{B,i}=\{J_i\}$ consists of finitely overlapping boxes.\footnote{The number of overlapping is $1000^{d-1}$} Each $J_i\in \cj_{B,i}$ (after permutation of coordinates) has the form $J_i=J_i^k\times J_i^{d-1-k}$ where $J_i^k$ is a $k$-dimensional cube of side length $r_{B,i}\rho_j $ and $J_i^{d-1-k}$ is an $(d-1-k)$-dimensional cube of side length $ 
(d-1) \wt C r_{B,i}\rho_j $.
\item For each box $J_i \in\cj_{B,i}$, there is a $(k,r_{B,i}\rho_j,\wt C)$-regular graph $G_{J_i}$ over $J_i^k$ and contained in $J_i$. In other words, for $l=k+1,\dots,d-1$, there exists a function $\phi_{B,J_i,l}: J_i^k\rightarrow \R$ satisfying
\begin{equation}\label{lowdimder}
    |\nabla^j \phi_{B,J_i,l}(x_{(k)})|\le \wt C (r_{B,i}\rho_j)^{-j+1}, \;\;\; j=1,\ldots, \Lambda 
\end{equation} 
for $x_{(k)}\in J^k_i$. Also $G_{J_i}$ is of the following form
\begin{equation*}\begin{split}
    G_{J_i}=\{ (x_{(k)}, & x_{k+1},\dots,x_{d-1})\in \R^{d-1}: 
    \\&x_{(k)}\in J_i^k, x_l=\phi_{B,J_i,l}(x_{(k)}), l=k+1,\dots,d-1 \}.
\end{split} 
\end{equation*}
Moreover, there exists a polynomial $Q_{J_i}$ of degree  $O_{d-1,\Lambda,E}(1)$ 
such that
\begin{equation}\label{0713.447}
    G_{J_i}=Z_{Q_{J_i}} \cap J_i.
\end{equation}
\end{enumerate}
They satisfy
\begin{equation}\label{lowdimcover}
    N_{r_B\rho_j}( S_B)\cap B^{d-1}\subset \bigcup_{i=1}^{\wt M} \bigcup_{J\in\cj_i} N_{(r_{B,i})^{A^2}\rho_j}(G_J)\cap J. 
\end{equation}

We have applied the induction hypothesis.
Now, we lift \eqref{lowdimcover} to the graph $\Psi_B$. For simplicity, we introduce the notation \[\vec \phi_{B,J}=(\phi_{B,J,k+1},\dots, \phi_{B,J,d-1}).\]
Then we have the expression
\begin{equation}\label{0713.448}
   G_{J_i}=\{(x_{(k)},\vec \phi_{B,J_i}(x_{(k)})): x_{(k)}\in J_i^k\} \subset \R^{d-1}. 
\end{equation}
For $J_i\in\cj_i$, we define the graph
\begin{equation}\label{gaBJ}
    \Ga_{B,J_i}:=\Big\{ \big(x_{(k)}, \vec{\phi}_{B,J_i}(x_{(k)}), \psi_B\big(x_{(k)}, \vec{\phi}_{B,J_i}(x_{(k)})\big)\big): x_{(k)}\in J_i^{k} \Big\} \subset \R^d.
\end{equation}
Note that this set is a subset of $\Psi_B=\{(y,\psi_B(y)): y\in B^{d-1}\}$. For a visualization of $G_{J_i}, \Ga_{B,J_i}$, see Figure \ref{fig2}.
\\

We next prove the following proposition. 

\begin{proposition}\label{0703.prop412}
    Suppose that
    \begin{equation}\label{tocheck}
        r_B\rho_j\ge 2K_j^{-A'},\ \ \  10d \widetilde{C} C_{d,D}\max\{(r_{B,i})^{A^2}\rho_j, K_j^{-A'}\}\le  (r_{B,i}\rho_j)^A.
    \end{equation}
    Then we have
    \begin{equation}\label{rewrt}
        N_{\mu}(Z)\cap [0,1]^d \subset \bigcup_{j=1}^D\bigcup_{B\in\cB_j}\bigcup_{i=1}^{\wt M}\bigcup_{J_i\in\cj_{B,i}} \sigma_{B,J_i} \big( N_{(r_{B,i}\rho_j)^A}(\Ga_{B,J_i})\cap (J_i^k\times \R^{d-k}) \big),
    \end{equation}
    where $\sigma_{B,J_i}$ is a permutation of coordinates. 
\end{proposition}

\begin{proof}
By \eqref{contain} and \eqref{0704.439}, this follows from
\begin{equation}\label{0704.448}
\begin{split}
         &N_{K_j^{-A'}}\big(\{ (y,\psi_B(y)): y\in N_{K_j^{-A'}}(S_B) \} \big)\cap B 
\\& \subset \bigcup_{i=1}^{\wt M}\bigcup_{J_i\in\cj_{B,i}} N_{(r_{B,i}\rho_j)^A}(\Ga_{B,J_i})\cap (J_i^k\times \R^{d-k})
\end{split}
\end{equation}
under the assumption of \eqref{tocheck}.
This immediately follows from combining three lemmas stated below (Lemma \ref{0712.lem413}, \ref{0712.lem414}, and \ref{0712.lem415}).
\end{proof}

\begin{lemma}\label{0712.lem413}
Assume that $r_B\rho_j\ge 2K_j^{-A'}$. Then we have
\begin{equation}\label{0704.449}
    \begin{split}
       & \mathrm{LHS \; of \;} \eqref{0704.448} \\&
       \subset \{ (y,x_d): y\in N_{r_B\rho_j}(S_B)\cap B^{d-1}, |x_d-\psi_B(y)|\le 10dC_{d,D}K_j^{-A'} \}\cap B.
    \end{split}
\end{equation}
\end{lemma}
\begin{proof}
Let $x$ belong to the LHS of \eqref{0704.448}. Write $x=(y,x_d)$. We need to show that
\begin{equation}
    y\in N_{r_B\rho_j}(S_B)\cap B^{d-1} \;\; \mathrm{ and } \;\; |x_d-\psi_B(y)|\le 10d C_{d,D}K_j^{-A'}
\end{equation} 
Note that there exists $x'=(y',x_d') \in B$  such that
$y' \in N_{K_j^{-A'}}(S_B)$, $|x-x'| \leq K_j^{-A'}$, and $x_d'=\psi_B(y')$. By the inequality
\begin{equation}
    |y-y'|\le |x-x'|\le K_j^{-A'}\le r_B\rho_j/2
\end{equation}
we see that $y\in N_{r_B\rho_j}(S_B)\cap B^{d-1}$. On the other hand, 
\begin{equation}
 |x_d-\psi_B(y)|\le |x_d-x_d'|+|\psi_B(y')-\psi_B(y)|\le 10dC_{d,D}K_j^{-A'}   
\end{equation}
by the inequailty $|\nabla \psi_B|\le C_{d,D}$.
\end{proof}

\begin{lemma}\label{0712.lem414}
\begin{equation*}
    \begin{split}
        &\Big\{ (y,x_d): y\in N_{r_B\rho_j}(S_B)\cap B^{d-1}, |x_d-\psi_B(y)|\le 10dC_{d,D}K_j^{-A'} \Big\}\cap B\\&
    \subset \bigcup_{i=1}^{\wt M} \bigcup_{J_i\in\cj_{B,i}} \Big\{ (y,x_d): y\in N_{(r_{B,i})^{A^2}\rho_j}(G_{J_i})\cap J_i, |x_d-\psi_B(y)|\le 10dC_{d,D}K_j^{-A'} \Big\}\cap B.
    \end{split}
\end{equation*}
\end{lemma}
    The above lemma simply follows from
    \eqref{lowdimcover}.

\begin{lemma}\label{0712.lem415}
If $100d \widetilde{C} C_{d,D}\max\{(r_{B,i})^{A^2}\rho_j, K_j^{-A'}\}\le  (r_{B,i}\rho_j)^A$, then
    \begin{align*}
        &\Big\{ (y,x_d): y\in N_{(r_{B,i})^{A^2}\rho_j}(G_{J_i})\cap J_i, \,|x_d-\psi_B(y)|\le 10dC_{d,D}K_j^{-A'} \Big\}\cap B\\
        &\subset N_{(r_{B,i}\rho_j)^A}(\Ga_{B,J_i})\cap (J_i^k\times \R^{d-k}).
    \end{align*}
\end{lemma}
\begin{proof}
    Let $x=(y,x_d)=(x_{(k)},x_{(d-1-k)},x_d)$ be a point in the set on the left hand side, where $x_{(k)}\in J_i^k$. Recall that 
    \[ G_{J_i}=\{(x_{(k)},\vec \phi_{B,J_i}(x_{(k)}) ): x_{(k)}\in J_i^k\} \subset \R^{d-1} \]
    with $|\nabla \vec \phi_{B,J_i}|\le d \widetilde{C}$.  Hence, by some elementary geometry, we have
    \begin{equation}
        |\vec \phi_{B,J_i}(x_{(k)})-x_{(d-1-k)}| \leq 10d \widetilde{C} \dist(y, G_{J_i}) \leq 10d \widetilde{C} (r_{B,i})^{A^2}\rho_j.
    \end{equation}
    Further by $|\nabla\psi_B|\le C_{d,D}$, we see that 
    \begin{equation}
    \begin{split}
        &\mathrm{dist}\Bigg(  (y,x_d), \Big(x_{(k)},\vec \phi_{B,J_i}(x_{(k)}),\psi_B\big(x_{(k)},\vec \phi_{B,J_i}(x_{(k)})\big)\Big) \Bigg)
        \\&\leq |y-\big(x_{(k)},\vec \phi_{B,J_i}(x_{(k)})\big)|+|x_d-\psi_B(y)|+|\psi_B(y)-\psi_B\big(x_{(k)},\vec \phi_{B,J_i}(x_{(k)})\big)|
        \\&\leq 
        10d  C_{d,D}(r_{B,i})^{A^2}\rho_j+10d  C_{d,D}K_j^{-A'}+10d \widetilde{C} C_{d,D}(r_{B,i})^{A^2}\rho_j
        \\&\le (r_{B,i}\rho_j)^A.
    \end{split}
    \end{equation} 
    This completes the proof.
\end{proof}

We have proved Proposition \ref{0703.prop412}. Note that $\Ga_{B,J_i}$ is a $(k,\rho_jr_{B,i}, C(d,E))$-regular graph for some large constant $C(d,E)$, which will be the constant in the statement of Theorem \ref{algthm}. This follows from the fact that $G_{J_i}$ is $(k,r_{B,i}\rho_j,\widetilde{C})$-regular graph (see \eqref{0713.448}) and the inequalities
\begin{equation}
\begin{split}
    \nabla^j \Big(\psi_B\big(x_{(k)},\vec{\phi}_{B,J}(x_{(k)})\big) \Big) \lesssim_{d,E} (\rho_j r_{B,i})^{-j+1},\;\;\; j=1,\ldots, \Lambda.
\end{split}
\end{equation}

In the next step, we will choose parameters carefully so that \eqref{tocheck} holds true. Then we will show that \eqref{rewrt} gives the desired covering (the right hand side of \eqref{cover}).

\bigskip

\subsection*{Step 5. Choosing parameters and closing the induction}
\hfill

We are ready to finish the proof of Theorem \ref{algthm}. It remains to carefully choose all the parameters so that we can close the induction. Recall that $D:=\deg P_1$.
We choose  
\begin{enumerate}[label=(\roman*)]
    \item $M(d,E):=D\wt M=D M(d-1,\lambda_{d,E})\le EM(d-1,\lambda_{d,E})$.
    \item $\{\mu_l\}_{1\le l\le M(d,E)}:=\{r_{B,i}\rho_j\}_{1\le j\le D, 1\le i\le \wt M}$.
\end{enumerate}
For each pair $(B, J_i)$ with $J_i\in\cj_{B,i}$, we choose a box $I_l:=I_l^k\times I_l^{d-k}$ with $I_l^k=J_i^k$ and $I_l^{d-k}$ being an $(d-k)$-dimensional cube of side length $C(d,E)d\mu_l$. We denote
\[ \Ga_{I_l}:=\Ga_{B,J_i}. \]
Since $\Ga_{B,J_i}$ is a $(k,\mu_l,C(d,E))$-regular graph over $I_l^k$, by Lemma \ref{inbox}, we have
\[ N_{(r_{B,i}\rho_j)^A}(\Ga_{B,J_i})\cap (J_i^k\times \R^{d-k})= N_{(r_{B,i}\rho_j)^A}(\Ga_{I_l})\cap I_l. \]
By \eqref{0713.447} and \eqref{0713.439}, we have
\begin{equation}
    \Ga_{I_l}= I_l \cap Z_{Q_{J_i}} \cap Z_{Q_B}.
\end{equation}
We can rewrite \eqref{rewrt} as
\begin{equation}\label{0704.454}
    N_{\mu}(Z)\cap [0,1]^d \subset \bigcup_{l=1}^{M(d,E)}N_{\mu_l^A}(\Ga_{I_l})\cap I_l.
\end{equation}
This gives \eqref{cover}. We have shown that items (ii) and (iii) in Theorem \ref{algthm} are true. We have also determined $C(d,E)$ and $M(d,E)$.

Recall that we obtained \eqref{0704.454} under the assumption of \eqref{tocheck}.
It remains to choose the parameters 
\[A',\ \{r_B\}_{B\in\cup_j\cB_j},\ c(d,E,A),\ \mu_0(d,E,A)\] 
so that  \eqref{tocheck} holds true. We will also check that $c(d,E,A)$ will  satisfy the item (i) of Theorem \ref{algthm}. This will finish the proof of Theorem \ref{algthm}.
\medskip

First of all, we choose
\begin{equation}\label{0706.455'}
 A':= \wt C_{d,E} \wt c^{-1} A^{2}
\end{equation}
for some large constant $\wt C_{d,E}$. Since $\wt c:=c(d-1,\lambda_{d,E},A^2) \leq 1$, $A'$ is much larger than $A$.
For each $1\le j\le D$, we claim that 
\begin{equation}\label{0706.456}
    \rho_j\le O_{d,E}(1)\mu^{(A'+1)^{-D}}.
\end{equation} 
To prove the claim, we first
recall that $\mu=K^{-1}$ and $\rho_j=(C_{d,D}K_j)^{-1}$.
By the condition in Lemma \ref{gzlem}, we have 
\begin{equation}\label{estKj}
    K_j^{-1}\le K_1^{-1}\sim_{d,E} K^{-(A'+1)^{-D}} = \mu^{(A'+1)^{-D}}.
\end{equation}
 This gives the claim. 

We next choose $\mu_0(d,E,A)$ sufficiently small so that the following two inequalities hold true:
\begin{equation}\label{condm01}
\begin{split}
&\mu_0(d,E,A)\leq  \mu(d-1,\lambda_{d,E},A^2)^{C_{d,E}''(A'+1)^D c(d-1,\lambda_{d,E},A^2)^{-1}}=c_{d,E}''\wt \mu_0^{C_{d,E}''(A'+1)^D \wt c^{-1}}
\\&
\mu_0(d,E,A)\leq (c_{d,E}'' /C(d-1,\lambda_{d,E}))^{(A'+1)^{E+1}} = (c_{d,E}'' \widetilde{C}^{-1})^{(A'+1)^{E+1}}.
\end{split}
\end{equation}
for some large constant $C''_{d,E}$ and small constant $c''_{d,E}$.
Note the right hand sides are already fixed numbers by induction hypothesis. If $c''_{d,E}$ is sufficiently small, then for $\mu\le \mu_0(d,E,A)$, we can plug in the first line of \eqref{condm01} into \eqref{0706.456} to get 
\begin{equation}\label{0706.458}
 \rho_j \le O_{d,E}(1) \wt\mu_0^{C''_{d,E}} < \wt \mu_0,   
\end{equation}
if $C''_{d,E}$ is large.
For each $B\in\cB_j$, we choose 
\begin{equation}\label{defrB}
    r_B:=\rho_j^{\wt c^{-1}}.
\end{equation} 
In Step 4, $r_B$ is required to be bounded by $\wt {\mu}_0$. By \eqref{defrB} and \eqref{0706.458}, this is guaranteed. 

Let us determine the constant $c(d,E,A)$. Choose
\begin{equation}\label{0713.469}
    c(d,E,A):=\frac{1}{2} (A'+1)^{-E}
\end{equation}
We need to show the item $(i)$ of Theorem \ref{algthm}. By \eqref{0706.456}, we have 
\begin{equation}
 \mu_l:=r_{B,i}\rho_j\le \rho_j \le O_{d,E}(1) \mu^{(A'+1)^{-D}} \le O_{d,E}(1) \mu^{(A'+1)^{-E}}.   
\end{equation}
To show $\mu_l \leq \mu^{c(d,E,A)}$, we need to show $\mu^{-\frac12 (A'+1)^{-E}}\ge O_{d,E}(1)$. By the second line of \eqref{condm01}, we have
 \begin{equation}
     \mu^{-\frac12 (A'+1)^{-E}}\ge (c''_{d,E})^{-1/2}\ge O_{d,E}(1),
 \end{equation}
 if $c''_{d,E}$ is small.
 This gives the upper bound of $\mu_l$.

Next, we show the lower bound of $\mu_l$: $\mu_l\ge \mu^{\frac{1}{2A}}$. Note that
\begin{equation}
    \mu_l=r_{B,i}\rho_j\ge r_B^{\frac{1}{2A^2}}\rho_j=\rho_j^{\frac{1}{2A^2}\cdot \wt c^{-1}+1}.
\end{equation}
Recall that $\rho_j=(C_{d,D}K_j)^{-1}$. By the definition in Lemma \eqref{gzlem}, we have
\begin{equation}
    K_j^{-1}\ge K_D^{-1}\sim_{d,E} K^{-(A'+1)^{-1}}=\mu^{(A'+1)^{-1}}.
\end{equation}
Therefore, we get
\begin{equation}
    \mu_l\ge \mu^{(\frac{1}{2A^2}\cdot \wt c^{-1}+1)\frac{1}{A'+1}}/ O_{d,E}(1).
\end{equation}
Recalling the definition of $A'$ in \eqref{0706.455'}, if $\wt C_{d,E}$ is big enough, then we have
\[ \mu_l\ge \mu^{\frac{1}{4A}}/O_{d,E}(1). \]
To show $\mu_l\ge \mu^{\frac{1}{2A}}$, we just need to show $\mu^{-\frac{1}{4A}}\ge O_{d,E}(1)$. This is done by the second line of \eqref{condm01}.
This verifies item $(i)$ of Theorem \ref{algthm}.
\medskip

We have determined all the constants.
Let us now check \eqref{tocheck}. Recall the inequalities that we want to prove: 
\begin{equation}\label{0706.458'}
\begin{split}
    &r_B\rho_j\ge 2K_j^{-A'},  
    \\&10d \widetilde{C} C_{d,D}(r_{B,i})^{A^2}\rho_j\le  (r_{B,i}\rho_j)^A,
          \\&
          10d \widetilde{C} C_{d,D}
          K_j^{-A'} \le  (r_{B,i}\rho_j)^A.
        \end{split}
    \end{equation}
    
Let us check the first inequality of \eqref{0706.458'}. By the definition of $\rho_j:=(C_{d,D}K_j)^{-1}$, it suffices to prove
\begin{equation}
 r_B\rho_j \geq 2(C_{d,D})^{A'} \rho_j^{A'}.   
\end{equation}
Since $r_B=\rho_j^{1/\wt c}$, it boils down to $\rho_j^{1+\wt c^{-1}-A'} \geq 2(C_{d,D})^{A'}$. By the choice of $A'$, we have $-1-\wt c^{-1}+A'\ge \frac{1}{2}A'$. We just need $\rho_j^{-1}\ge 4 C_{d,D}^2$, or equivalently, 
\begin{equation}\label{KjlessC}
    K_j\ge 4 C_{d,D}.
\end{equation} To do this, we recall \eqref{estKj} and the second condition in \eqref{condm01}. We have $K_j\ge (c''_{d,E})^{-1}/O_{d,E}(1)$
Recalling that the number $C_{d,D}$ is defined in Proposition \ref{0630.prop49}, if $c''_{d,E}$ is sufficiently small, then \eqref{KjlessC} holds.

\medskip

We next check the second inequality of \eqref{0706.458'}.
 Since $r_{B,i}\le (r_B)^{\wt c}= \rho_j$, it suffices to show
\[ 10d \widetilde{C}C_{d,D}\le \rho_j^{-(A^2+1-2A)}, \]
or equivalently,
\begin{equation}\label{0706.462}
    (\rho_j)^{A^2+1-2A} \leq (10d\widetilde{C}C_{d,D})^{-1}.
\end{equation}
By \eqref{0706.456} and $\mu \leq \mu_0$ and the second condition in \eqref{condm01}, we have \[\rho_j\le O_{d,E}(1)\mu^{(A'+1)^{-D}}\leq O_{d,E}(1) c''_{d,E} (\widetilde{C})^{-1}.\]
Since $A>100$ and $c''_{d,E}$ is a small constant, this implies \eqref{0706.462}.

\medskip

Lastly, we check the last inequality of \eqref{0706.458'}. Plugging $K_j^{-1}=C_{d,D}\rho_j$ and $r_{B,i}\ge r_B^{\frac{1}{2A^2}}=\rho_j^{\frac{1}{2A^2 \widetilde{c}}}$, it suffices to show 
\begin{equation}
    10d \widetilde{C}(C_{d,D})^{1+A'}(\rho_j)^{A'} \leq (\rho_j)^{A(1+1/(2A^2\widetilde{c}))}.
\end{equation}
This inequality follows simply by noting that $A'$ is much larger than $A$ (see \eqref{0706.455'}), and $C_{d,D}^2\le (4\rho_j)^{-1}$ (see \eqref{KjlessC}).

\section{Examples: good manifolds}
In this section, we prove Theorem \ref{0409thm16}. To prove the theorem, we need to calculate the lower dimensional decouplings and the parameter $X(\mathcal{M},k,m)$.

Let us go back to the proof of Theorem \ref{0409thm16}. We are interested in the range of
\begin{equation}\label{0712.51}
   p> \min_{3 \leq k \leq d+1} \max (\frac{2(k+1)}{k-1}, \frac{2(2d-k+6)}{2d-k+2})=:p_c.
\end{equation}
Two numbers $\frac{2(k+1)}{k-1}$ and $\frac{2(2d-k+6)}{2d-k+2}$ are equal to each other when $k=2(d+2)/3$. This gives the asymptotics $p>2+\frac{6}{d}+O(\frac{1}{d^2})$.

For the rest of the section, we assume that $\mathcal{M}$ is a good manifold.
By Theorem \ref{restriction} and the epsilon removal lemma by \cite{MR1666558}, it suffices to prove that for $p>p_c$ there exists $3 \leq k \leq d+1$ such that
\begin{equation*}
\begin{split}
 \delta^{d-\frac{2d+2n}{p}}\delta^{-\mathrm{Dec}_p(\mathcal{M}|_{L_{k-2}} )}
+
       \sup_{0 \leq m \leq d+n}(\delta^{-\frac{m}{p}+\frac12{X(\mathcal{M},k,m)}})  \lesssim_{\epsilon} \delta^{-\epsilon}.
\end{split}
\end{equation*}
In other words, for $p>p_c$, it suffices to prove that for some $3 \leq k \leq d+1$
\begin{equation}\label{0705.52'}
    \begin{split}
         \delta^{-\mathrm{Dec}_p(\mathcal{M}|_{L_{k-2}} )} \lesssim_{\epsilon} \delta^{-\epsilon}\delta^{-d+\frac{2d+4}{p}},
    \end{split}
\end{equation}
and
\begin{equation}\label{2ndbd}
      \frac{2m}{p} \leq X(\mathcal{M},k,m)\ \ \   \textup{~holds~for~all~} 0 \leq m \leq d+2.
\end{equation}

\subsection{Decoupling inequality}
\eqref{0705.52'} follows from
\begin{equation}\label{0705.51}
    D_p(\mathcal{M}|_{L_{k-2}}) \lesssim_{\epsilon} \delta^{-\epsilon}\delta^{-d+\frac{2d+4}{p}} \;\;\;\; \mathrm{for} \;\; p>p_c.
\end{equation}
Write $\mathcal{M}=\{(\xi,P(\xi),Q(\xi)): \xi \in \R^{d} \}$ and ${\bf{Q}}:=(P,Q)$. Recall Definition \ref{0706.def15}.
\begin{lemma}\label{0705.lem51}
Suppose that $\mathcal{M}$ is a good manifold. Then
    \begin{equation}\label{0705.52}
        \fd_{d-m,2}({\bf Q}) \geq d-m, \;\;\; \fd_{d-m,1}({\bf Q}) \geq d-2m-1
    \end{equation}
    for all $0 \leq m \leq d$.\footnote{This lemma still holds true for manifolds satisfying the condition in Remark \ref{0707.remark111}}
\end{lemma}

\begin{proof}

Let us show that the first inequality of \eqref{0705.52}  under the condition mentioned in Remark \ref{0707.remark111}. For a contradiction, assume that $\fd_{d-m,2}({\bf{Q}})<d-m$. By Corollary \ref{0717.a10} and Lemma \ref{narrowrange}, there exists  $M \in \R^{d \times d}$ with rank $d-m$ such that 
\begin{equation}
    \Rrank(xM AM^T+yM B M^T) <d-m.
\end{equation}
Here 
\[ A=\textup{diag}\{a_1,\dots,a_d\}, B=\textup{diag}\{b_1,\dots,b_d\} \]
are the Hessian matrices of $P,Q$. We refer to Definition \ref{0707.defa3} for the definition of Row-rank.
By the definition of Row-rank, there exist linearly independent vectors ($1\times d$ matrices) $\{\bv_i\}_{i=1}^{m+1}$ such that
\begin{equation}
    \bv_i \cdot \big( xMAM^T+  yMBM^T \big) \equiv 0
\end{equation}
for $i=1,\ldots,m+1$. In other words,
\begin{equation}
    \bv_i \cdot MAM^T = 0, \;\;\; \bv_i \cdot  MBM^T = 0
\end{equation}
for $i=1,\ldots,m+1$. Since $M$ has rank $d-m$, by the rank-nullity theorem, there exists a nonzero vector ($1\times d$ matrix) ${\bf w}$ such that
\begin{equation}
    {\bf w}= \sum_{i=1}^{m+1}c_i \bv_iM
\end{equation}
for some constants $c_i \in \R$.
This implies that
\begin{equation}
    \bw A \bw^{T} =0, \;\;\; \bw B\bw^{T}=0.
\end{equation}
In other words, there exists a nonzero vector $(w_1,\ldots,w_d) \in \R^d$ such that
\begin{equation}
    \sum_{i=1}^d a_i (w_i)^2=0, \;\;\; \sum_{i=1}^d b_i (w_i)^2=0.
\end{equation}
This gives the contradiction.

Let us next show that the second inequality of \eqref{0705.52} under the condition mentioned in Remark \ref{0707.remark111}. Note that for any nonzero $t' \in \R$ the function $P+t'Q$ is of the form $H_{\widetilde{a}}:=\sum_{i=1}^{d}\widetilde{a}_i \xi_i^2$ and at most one of $\widetilde{a}_i$ is equal to zero. Observe that
\begin{equation}
    \fd_{d-m,1}({\bf Q}) \geq  \inf\big( \fd_{d-m,1}(H_{\widetilde{a}})\big)
\end{equation}
where the infimum runs over $H_{\widetilde{a}}$ such that at most one of $\widetilde{a}_i$ is zero. By Lemma 3.1 of \cite{MR4541334}, we  have $\fd_{d-m,1}(\bQ) \geq d-1-2m$.
\end{proof}

To prove \eqref{0705.51}, let us introduce some notations.
Let $Q$ be a quadratic form.
Let $H$ be a $(k-2)$-dimensional subspace in $\R^d$. Let $\mathrm{Rot}_H$ be a rotation on $\R^d$ that maps $\{\xi \in \R^d: \xi_{k-1}=\cdots=\xi_d=0 \}$ to $H$. Define
\begin{equation}
    Q|_H(\xi'):=Q((\xi',{\bf 0}) \cdot (\mathrm{Rot}_H)^{T})
\end{equation}
with $\xi' \in \R^{k-2}$ and ${\bf 0}=(0,\ldots,0) \in \R^{d-k+2}$. We denote by
\begin{equation}
    {\bf Q}|_H:=(P|_H,Q|_H).
\end{equation}
Define $\mathcal{M}|_H:=\{(\xi',{\bf{Q}}|_H(\xi')) \in [0,1]^{k-2}  \times \R^2\}$.

Let $L_{k-2}$ be a $(k-2)$-dimensional linear subspace of $\R^d$.
Note that 
\begin{equation}
    \fd_{d',n'}({\bf Q})
    \leq \fd_{d',n'}({\bf Q}|_{L_{k-2}})
\end{equation}
for $d' \leq k-2$ and $n' \leq 2$.
By this and Lemma \ref{0705.lem51},
\begin{equation}
    \begin{split}
        &\fd_{k-2-m,2}({\bf{Q}}|_{L_{k-2}}) \geq k-2-m, \\
        &\fd_{k-2-m,1}({\bf Q }|_{L_{k-2}}) \geq \fd_{d-(d-k+2+m),1}({\bf Q })\ge \max\{0, d-2(d-k+2+m)-1\}. 
    \end{split}
\end{equation}
We also trivially have
\[ \fd_{k-2-m,0}({\bf Q }|_{L_{k-2}})=0. \]
Note that $\mathcal{M}|_{L_{k-2}}$ is a $(k-2)$-dimensional manifold of codimension two.  By an uncertainty principle, the Fourier transform of $E^{\mathcal{M}}_{N_{\delta}(L|_{k-2})}f$ is contained in (after a linear transformation) the $\delta$-neighborhood of 
\begin{equation}
    \mathcal{M}|_{L_{k-2}} \times \R^{d-k+2}.
\end{equation}
By Corollary 1.2 of \cite{MR4541334} and the proof of Theorem 2.2 of \cite{MR4143735}, we have 
\begin{equation}\label{0705.56}
    D_p(\mathcal{M}|_{L_{k-2}})   \lesssim \max(\delta^{-(k-2)(\frac12-\frac1p) },\delta^{-\min\{2(k-2)(\frac12-\frac1p)-\frac2p,(d+1)(\frac12-\frac1p)-\frac2p\}},\delta^{-(k-2)+\frac{2(k-2)+4}{p} }).
\end{equation} 
 Recall that we are trying to prove \eqref{0705.51}. The last term on the right hand side of \eqref{0705.56} is always smaller than the right hand side of \eqref{0705.51} for $2\le k\le d+1$. So, we just need to find $2\le k\le d+1$ and $p \geq 2$ that satisfy
\begin{equation}
    \max(\delta^{-(k-2)(\frac12-\frac1p) },\delta^{-\min\{2(k-2)(\frac12-\frac1p)-\frac2p,(d+1)(\frac12-\frac1p)-\frac2p\}}  ) \lesssim \delta^{-d+\frac{2d+4}{p}}.
\end{equation}
This is equivalent to 
\begin{equation}\label{1stlowbd}
    p> \max {\big(\frac{2(2d-k+6)}{2d-k+2},\min \{ \frac{2(d-k)+6}{d-k+2} ,\frac{2(d+1)}{d-1}\} \big)}.
\end{equation}

\subsection{\texorpdfstring{$k$}{}-linear restriction estimate} 

In this subsection, we deal with the inequality \eqref{2ndbd}. It is done by proving the following proposition.

\begin{proposition}\label{0410prop53}
    Suppose that $\mathcal{M}$ is a good manifold. Let $3 \leq k \leq d+1$. Then
    \begin{equation}
    m \cdot \frac{k-1}{k+1} \leq X\big(\mathcal{M},k,m\big)
\end{equation}
for all $0 \leq m \leq d+2$.\footnote{This lemma still holds true for manifolds satisfying the condition in Remark \ref{0707.remark111}}
\end{proposition}

Let us postpone the proof of the proposition, and first prove our main goal: for $p>p_c$, there is some $3\le k\le d+1$ such that \eqref{0705.52'} and \eqref{2ndbd} hold.

By \eqref{1stlowbd} and Proposition \ref{0410prop53}, if
\begin{equation}
    p> \max (\frac{2(2d-k+6)}{2d-k+2},\min \{ \frac{2(d-k)+6}{d-k+2} ,\frac{2(d+1)}{d-1}\},\frac{2(k+1)}{k-1}),
\end{equation}
then \eqref{0705.52'} and \eqref{2ndbd} hold. Therefore, we just need to show
\begin{equation}
    p_c\ge \min_{3\le k\le d+1}\max (\frac{2(2d-k+6)}{2d-k+2},\min \{ \frac{2(d-k)+6}{d-k+2} ,\frac{2(d+1)}{d-1}\},\frac{2(k+1)}{k-1}).
\end{equation}
By the definition of $p_c$ in \eqref{0712.51}, we just need to show for any $3\le k\le d+1$,
\begin{equation}\label{ineqfinal}
    \max (\frac{2(2d-k+6)}{2d-k+2},\frac{2(k+1)}{k-1})\ge \min \{ \frac{2(d-k)+6}{d-k+2} ,\frac{2(d+1)}{d-1}\}.
\end{equation}
When $k\le d$, we have $ \frac{2(d+1)}{d-1}\le \frac{2(k+1)}{k-1}$, so LHS $\ge$ RHS. The remaining case is $k=d+1$, for which \eqref{ineqfinal} becomes
\begin{equation}
    \max (\frac{2(d+5)}{d+1},\frac{2(d+2)}{d})\ge \min \{ 4 ,\frac{2(d+1)}{d-1}\}.
\end{equation}
This is ture by noting: $\frac{2(d+5)}{d+1}\ge 4$ if $d\le 3$; $\frac{2(d+5)}{d+1}\ge \frac{2(d+1)}{d-1}$ if $d\ge 3$.



\bigskip

Let us  give a proof of Proposition \ref{0410prop53}. Recall that
for each $\xi \in [0,1]^d$, we denote by $V_{\xi}$ the tangent space of $\mathcal{M}$ at the point $(\xi,\mathbf{Q}(\xi))$, and 
 by $\pi_{V_{\xi}}$ the orthogonal projection of $\R^{d+n}$ to $V_{\xi}$. We note a simple lemma for a linear algebra.

\begin{lemma}\label{0408lem}
Let $2 \leq k \leq d+1$ and $0 \leq m \leq d+2$. Fix $i \in \mathbb{N}$. 
    The followings are equivalent.
    \begin{enumerate}
        \item $\sup_{\dim V=m} \dim \{\xi:\dim (\pi_{V_\xi}(V))<m-2+i \} \leq k-2 $.

        \item $\sup_{\dim W=d+2-m} \dim \{\xi: \dim(\pi_{V_{\xi}^{\bot}}(W) ) < i  \} \leq k-2$.
    \end{enumerate}
\end{lemma}

\begin{proof}
We claim that
\begin{equation}
    \dim \pi_{V_\xi}V = \dim V - \dim V_{\xi}^{\bot} + \dim \pi_{V_\xi^{\bot}}V^{\bot}.
\end{equation}
for any $V$ with dimension $m$. Note that 
\begin{equation}
    \dim V -\dim V_{\xi}^{\bot}=m-2.
\end{equation}
So the lemma follows from the claim by taking $W:=V^{\bot}$. 

It remains to prove the claim. We first note that for any subspace $A,B$ we have
\begin{equation}
    \dim \pi_A B= \dim B -\dim (A^{\bot} \cap B).
\end{equation}
This simply follows by considering the map $\pi_A:B \rightarrow \pi_A B$ and noting that the kernel of the map is $B \cap A^{\bot}$. By this fact,
\begin{equation}\label{0705.520}
    \dim \pi_{V_\xi}V = \dim V - \dim (V_{\xi}^{\bot} \cap V ).
\end{equation}
We write $V_{\xi}^{\bot} \cap V = V_{\xi}^{\bot} \cap (V^{\bot})^{\bot}$. Plugging with $A=V^{\bot}$ and $B=V_{\xi}^{\bot}$, we get
\begin{equation}\label{0705.521}
    \dim(V_{\xi}^{\bot}\cap V)=-\dim \pi_{V^{\bot}}V_{\xi}^{\bot}+\dim V_{\xi}^{\bot}.
\end{equation}
Note that
\begin{equation}\label{0705.522}
    \dim \pi_{V^{\bot}}V_{\xi}^{\bot}=\dim \pi_{V_{\xi}^{\bot}}V^{\bot}.
\end{equation}
The claim follows by combining \eqref{0705.520}, \eqref{0705.521}, and \eqref{0705.522}.
\end{proof}

\begin{proof}[Proof of Proposition \ref{0410prop53}]

For simplicity, define $p_k:=\frac{k+1}{k-1}$, and we will  prove
\begin{equation}
    \sup_{\dim V=m} \dim \Big\{\xi:\dim (\pi_{V_\xi}(V))< \frac{m}{p_k} \Big\} \leq k-2. 
\end{equation}
for all $m$.
By Lemma \ref{0408lem}, this is equivalent to
\begin{equation}
    \sup_{\dim W=d+2-m} \dim \Big\{\xi: \dim(\pi_{V_{\xi}^{\bot}}(W) ) < 2-m\big(1-\frac{1}{p_k}\big)  \Big\} \leq k-2.
\end{equation}
The above set is empty provided that $2-m(1-\frac{1}{p_k}) \leq 0$. So we consider $m$ satisfying
$
    m < \frac{2p_k}{p_k-1}.$ Note that the dimension of the projection of a subspace is an integer.

    \subsubsection*{{Case 1.} $\frac{p_k}{p_k-1} \leq m < \frac{2p_k}{p_k-1}$}

We need to prove
\begin{equation}\label{0408eq47}
    \sup_{\dim W=d+2-m}\dim \Big\{ \xi: \dim(\pi_{V_{\xi}^{\bot}}(W) )=0\Big\} \leq k-2.
\end{equation}
We may assume that $m=\frac{2p_k}{p_k-1}-1$ as it is the worst case. Recall that $p_k=\frac{k+1}{k-1}$, so $$\dim W=d+3-\frac{2p_k}{p_k-1}=d-k+2.$$
For simplicity, we use the notation $P_i:=\partial_{i}P=2a_i\xi_i$ and $Q_i:=\partial_iQ=2b_i\xi_i$. By abuse of notation, we use $V_{\xi}^{\bot}$ to denote the $2 \times (d+2)$ matrix:
\begin{equation}\label{0408eq47'}
    \begin{pmatrix}
        P_1 & P_2 & \cdots & P_d & -1 & 0 \\
        Q_1 & Q_2 & \cdots & Q_d & 0 & -1
    \end{pmatrix}.
\end{equation}
We also use $W$ to denote the $(d+2)\times (d+2-m)$ matrix (not unique) whose column vectors span the space $W$. We have
\begin{equation}\label{0804.536}
   \dim(\pi_{V_\xi^\perp}(W))=\rank(V_\xi^\perp W). 
\end{equation}

We would like to simplify the matrices $V_\xi^\perp$ and $W$.
Note that for a invertible matrix $B$, we have $\rank(V_\xi^\perp W)=\rank(V_\xi^\perp  W B)$. We may replace $W$ by $WB$. This is equivalent to do column operations to $W$.

After the matrix operations, we assume that $W$ has the column-echelon form, i.e., there is a $(\dim W\times \dim W)$-minor of $W$ which is the identity matrix. There are three scenarios on whether this identity matrix has some entries on the last two rows of $W$.

If this identity matrix has no entry on the last two rows of $W$, we are in the first scenario, in which we may assume $W$ is of form \eqref{0705.527}.

If this identity matrix has one row on the last two rows of $W$, we are in the second scenario, in which we may assume $W$ is of form \eqref{0705.528}.

If this identity matrix has two rows on the last two rows of $W$, we are in the third scenario, in which we may assume $W$ is of form \eqref{0705.529}.
\begin{enumerate}
    \item \begin{equation}\label{0705.527}
        \left[
\begin{array}{c}
    I_{\dim W \times \dim W}   \\
    *     
\end{array} \right]
    \end{equation}
    \item \begin{equation}\label{0705.528}
        \left[
\begin{array}{c|c}
    I_{(\dim W-1) \times (\dim W-1)} &   \\
    * & e_{d+2} \\
    * &  
\end{array} \right]
    \end{equation}
    \item \begin{equation}\label{0705.529}
        \left[
\begin{array}{c|c|c}
    I_{(\dim W-2) \times (\dim W-2)} &  &  \\
    * & e_{d+1} &e_{d+2} \\
    * & &
\end{array} \right]
    \end{equation}
\end{enumerate}
Here $I_{*}$ is an identity matrix and $e_i=(0,\ldots,0,1,0,\ldots,0)^T$ whose  $i$th component is one. Note that if we have $\dim (\pi_{V_{\xi}^{\bot}}(W))=0$ then the second and third forms cannot happen. So it suffices to consider only the first form.
By direct computations, 
\begin{equation}
    V_{\xi}^{\bot}W = \begin{pmatrix}
        P_1+F^{(1)} & P_2+F^{(2)} & \cdots & P_{\dim W} +F^{(\dim W)}
        \\
        Q_1+G^{(1)} & Q_2+G^{(2)} & \cdots &Q_{\dim W}+G^{(\dim W)}
    \end{pmatrix}.
\end{equation}
Here the functions $F^{(*)}$ and $G^{(*)}$ are linear functions depending only on $\xi_{\dim W+1}, \cdots , \xi_{d}$. Recall \eqref{0804.536}.
 The inequality \eqref{0408eq47} follows from the following lemma. 

\begin{lemma}\label{0409lem45}
    Under the settings above, we have
    \begin{equation}
        \bigcap_{i=1}^{\dim W} Z(P_i+F^{(i)},Q_i+G^{(i)})
    \end{equation}
    is contained in a $(k-2)$-dimensional plane in $\mathbb{R}^d$.
\end{lemma}

\begin{proof}
By the definition of $P$ and $Q$,
we have
\begin{equation}\label{0804.542}
    \bigcap_{i=1}^{\dim W} Z(P_i+F^{(i)},Q_i+G^{(i)})=\bigcap_{i=1}^{\dim W} Z(a_i\xi_i+ F^{(i)},b_i\xi_i+G^{(i)}).
\end{equation}
By the condition \eqref{0707.115}, we see that for any $i$, two numbers $a_i$ and $b_i$ cannot simultaneously be $0$. Hence, we can pick $\dim W=d-k+2$ equations from 
\begin{equation}
 \{a_i\xi_i+ F^{(i)}=0\},\;\;   \{b_i\xi_i+ G^{(i)}=0\},\;\; i=1,\ldots,d-k+2
\end{equation}
whose coefficients are linearly independent. We see that  the set \eqref{0804.542} is contained in the zero set of $d-k+2$ independent linear equations, so it is contained in a $(k-2)$-dimensional plane. 
\end{proof}


This finishes the proof of \textit{Case 1}. Next, we consider \textit{Case 2}.

\bigskip

   \subsubsection*{{Case 2.} $0 \leq m < \frac{p_k}{p_k-1}$}
  We need to prove
\begin{equation}\label{0408eq4a}
    \sup_{\dim W=d+2-m}\dim \Big\{ \xi: \dim(\pi_{V_{\xi}^{\bot}}(W) ) \leq 1\Big\} \leq k-2.
\end{equation}
We may assume that $m= \lceil \frac{k+1}{2} \rceil -1$  as it is the worst case. Recall that $W$ is of the forms \eqref{0705.527}, \eqref{0705.528}, and \eqref{0705.529}. If $\dim(\pi_{V_\xi^{\bot}}(W)) \leq 1$ then \eqref{0705.529} cannot happen. 

Let us consider  the case \eqref{0705.527}. By direct computations,
\begin{equation}
    V_{\xi}^{\bot}W = \begin{pmatrix}
        P_1+F^{(1)} & P_2+F^{(2)} & \cdots & P_{\dim W} +F^{(\dim W)}
        \\
        Q_1+G^{(1)} & Q_2+G^{(2)} & \cdots &Q_{\dim W}+G^{(\dim W)}
    \end{pmatrix}.
\end{equation}
Here the functions $F^{(*)}$ and $G^{(*)}$ are linear functions depending only on $\xi_{\dim W+1}, \cdots , \xi_{d}$. We write
\begin{equation}\label{0706.536}
\begin{split}
    &\dim \Big\{ \xi: \dim(\pi_{V_{\xi}^{\bot}}(W) ) \leq 1\Big\} 
    \\& = \bigcup_{\lambda \in \R} \bigcap_{i=1}^{\dim W} \{ \xi \in \R^d: P_i(\xi)+F^{(i)}(\xi)-\lambda(Q_i(\xi)+G^{(i)}(\xi) )=0 \}.
\end{split}
\end{equation}
Consider the $\dim W \times d$ matrix
\begin{equation}\label{0706.537}
    \begin{pmatrix}
        \nabla (P_1 + F^{(1)}-\lambda (Q_1+G^{(1)} )) \\
        \vdots \\
         \nabla (P_{\dim W} + F^{(\dim W)}-\lambda (Q_{\dim W}+G^{(\dim W)} )).
    \end{pmatrix}
\end{equation}
This matrix is of the form
\begin{equation}
    \begin{pmatrix}
        2(a_1-\lambda b_1) & 0 & \cdots & 0 & * &\cdots &* \\
        0 & 2(a_2-\lambda b_2) & \cdots & 0 & * & \cdots & * \\
        \vdots & \vdots & \ddots & \vdots & \vdots & \ddots & \vdots \\
        0 & 0 & \cdots & 2(a_{\dim W}-\lambda b_{\dim W}) & * & \cdots & *
    \end{pmatrix}
\end{equation}
Note that this matrix is independent of variables $\xi_1,\ldots,\xi_d$ because $P_i, F^{(i)}, Q_i, G^{(i)}$ are linear functions. By calculating minors, one can see that this matrix has full rank except for $\leq \dim W$ many $\lambda$. This means that except for such $\lambda$ we have
\begin{equation}
    \dim \Big(\bigcap_{i=1}^{\dim W} \{ \xi \in \R^d: P_i(\xi)+F^{(i)}(\xi)-\lambda(Q_i(\xi)+G^{(i)}(\xi) )=0 \} \Big) = d-\dim W.
\end{equation}
We now suppose that the rank of \eqref{0706.537} is not full. By the definition of a good manifold, by considering the first $\dim W$ many columns, we see that $\lambda$ must satisfy
\begin{equation}\label{0706.539}
    (a_1-\lambda b_1) \cdots (a_{\dim W}-\lambda b_{\dim W})=0.
\end{equation}
By the definition of a good manifold, the determinant of
\begin{equation}
    \begin{pmatrix}
        a_i & a_j \\
        b_i & b_j
    \end{pmatrix}
\end{equation}
is not equal to zero for $i \neq j$. Thus, if $\lambda$ satisfies the equation \eqref{0706.539}, the rank of $\eqref{0706.537}$ is equal to $\dim W-1$.
By combining all the information we have obtained, we have
\begin{equation*}
\begin{split}
    \dim \Big( \bigcup_{\lambda \in \R} \bigcap_{i=1}^{\dim W} \{ \xi \in \R^d: P_i(\xi)+F^{(i)}(\xi)-\lambda(Q_i(\xi)+G^{(i)}(\xi) )=0 \} \Big)
    \leq 1+d-\dim W.
\end{split}
\end{equation*}
By \eqref{0706.536}, it remains to prove
\begin{equation}
    1+d-\dim W \leq k-2.
\end{equation}
Recall that $\dim W=d+2-m$ and $m=\lceil \frac{k+1}{2} \rceil -1$. Then the inequality becomes
\begin{equation}
    \lceil \frac{k+1}{2} \rceil \leq k.
\end{equation}
This is true for any $k \geq 1$.

Let us next consider  the case \eqref{0705.528}. By direct computations, 
\begin{equation}
    V_{\xi}^{\bot}W = \begin{pmatrix}
        P_1+F^{(1)} & P_2+F^{(2)} & \cdots & P_{\dim W-1} +F^{(\dim W-1)} & 0
        \\
        Q_1+G^{(1)} & Q_2+G^{(2)} & \cdots &Q_{\dim W-1}+G^{(\dim W-1)} & -1
    \end{pmatrix}.
\end{equation} 
Note that
\begin{equation}\label{0706.545}
\begin{split}
    \dim \Big\{ \xi: \dim(\pi_{V_{\xi}^{\bot}}(W) ) \leq 1\Big\} 
     =  \bigcap_{i=1}^{\dim W-1} \{ \xi \in \R^d: P_i(\xi)+F^{(i)}(\xi)=0 \}.
\end{split}
\end{equation}
The dimension of this set is greater or equal to $d-(\dim W-2)$. It remains to prove
\begin{equation}
    d+2-\dim W \leq k-2.
\end{equation}
After some computations, this amounts to proving
\begin{equation}
    \lceil \frac{k+1}{2} \rceil \leq k-1.
\end{equation}
This is true for any $k \geq 3$.  So this completes the proof. \end{proof}

\appendix
\section{Algorithm to compute \texorpdfstring{$\fd_{d',n'}({\bf Q})$}{}}

We provide an equivalent formulation of the minimal number of variables (Collary \ref{0717.a10}) and an algorithm to calculate the number.

\subsection{Equivalent formulation of \texorpdfstring{$\fd_{d',n'}({\bf Q})$}{}}

Any quadratic polynomial $Q$ in $d$ variables uniquely determines a symmetric $d\times d$ matrix $A$ whose entry is given by $a_{ij}=\frac{1}{2}\partial_i\partial_j Q$. On the other hand, any symmetric $d\times d$ matrix $A=(a_{ij})$ uniquely determines the polynomial $Q(\xi)=\sum_{i,j=0}^d a_{ij}\xi_i\xi_j$. This gives a correspondence between quadratic polynomials in $d$ variables and symmetric $d\times d$ matrices. Under the change of variables $\xi\rightarrow M\xi$, we see that the new polynomial $Q \circ M$ corresponds to the matrix $M^T A M$.

It is not hard to see:
\begin{lemma} Let $Q$ be a quadratic polynomial and $A$ be its corresponding matrix. Then $NV(P)$ equals to the number of nonzero rows of $A$ (or nonzero columns since $A$ is symmetric).
\end{lemma}

We also have:
\begin{lemma}
Let $Q$ be a quadratic polynomial (in $d$ variables) and $A$ be its corresponding matrix. Then we have
    \[ \inf_{\substack{ M\in \R^{d\times d}:\\ \rank(M)=d}} NV(Q\circ M)=\rank (A). \]
\end{lemma}
\begin{proof}
    We first show ``$\ge$". For any $d\times d$ matrix $M$ with full rank, we have $NV(Q\circ M)=$ number of nonzero rows of $M^T A M\ge$ $\rank(M^T A M)=\rank(A)$.

    Next we show ``$\le$". Since $A$ is symmetric, we can always find a matrix $M$ of full rank $d$, such that $M^T A M$ is diagonal. We see that $\rank(A)=\rank(M^T A M)=NV(Q\circ M)$.
\end{proof}

\begin{definition}\label{0707.defa3}
Let $R$ be a $\R$-vector space. (For example, $R=\R$ or $R=\R[x_1,\dots,x_n]$).
    Given a matrix $B\in R^{n\times n}$, we define the Row-rank of $B$ to be the rank of the $\R$-vector space spanned by the row of $B$, denoted by $\Rrank(B)$.
\end{definition}

\begin{remark}
{\rm
When $R$ is a ring, $\det B$ is well-defined for $B\in R^{n\times n}$. However, it is not true that $\Rrank(B)=n\Leftrightarrow \det B\neq 0$. The counterexample is 
$$B=\begin{pmatrix}
    x_1 & x_1\\
    x_2 & x_2
\end{pmatrix},$$
whose determinant is $0$ but has row-rank 2.
\medskip

If $B\in R^{n\times n}$ is an $n\times n$ matrix with entries in $R$, and $M\in GL_n(\R)$ is an invertible $n\times n$ matrix with entries in $\R$, then we have
\begin{equation}\label{rowrankeq}
    \Rrank(B)=\Rrank(M B)=\Rrank(B M). 
\end{equation} 
It is not hard to see $\Rrank(B)=\Rrank(M B)$ since the row spaces of $B$ and $MB$ are the same. To show $\Rrank(B)=\Rrank(B M)$, it suffices to show $\Rrank(B)\ge\Rrank(B M)$, since then we have 
\[\Rrank(B)\ge\Rrank(B M)\ge \Rrank(BMM^{-1})=\Rrank(B).\]
Let us denote $r=\Rrank(B)$, and without loss of generality assume the row space of $B$ is spanned by the first $r$ rows of $B$. Then we see that the row space of $BM$ is spanned by the first $r$ rows of $BM$, which means $\Rrank(BM)\le r$.

}
\end{remark}

\begin{definition}\label{defrank}
    For a finite subset $\cB=\{B_1,\dots,B_l\}\subset \Sym_n(\R)$, we define
    \[ \rank(\cB):=\Rrank(\sum_{i=1}^l x_iB_i). \]
Here, $\Rrank$ is defined as in Definition \ref{0707.defa3} where we use $R=\R[x_1,\dots,x_l]$.    
These $x_i$ are indeterminates.
\end{definition}

It is not hard to see that if $\cB,\cB'\subset \Sym_n(\R)$ satisfy $\spn(\cB)=\spn(\cB')$, then
\[ \rank(\cB)=\rank(\cB'). \]
Therefore, for any subset $\cB\subset \Sym_n(\R)$ (possibly infinite), we can also define $\rank(\cB)$ in the following way. We choose a finite subset $\cB_0\subset \cB$ such that $\spn(\cB)=\spn(\cB_0)$ and define
\[ \rank(\cB):=\rank(\cB_0). \]

\bigskip

Suppose ${\bQ}=(Q_1,\dots,Q_n)$ are quadratic polynomials in $d$ variables, and let $(A_1,\dots,A_n)$ be the corresponding $d\times d$ matrices. We will find a new way to express $\fd_{d',n'}(\bQ)$.
The next lemma connects the rank defined in Definition \ref{defrank} with number of variables.

\begin{lemma}\label{lemA6}
    \[ \inf_{\substack{ M\in \R^{d\times d}:\\ \rank(M)=d}} NV({\bQ}\circ M)=\rank (\{A_1,\dots,A_n\}).\]
\end{lemma}
\begin{proof}
    We first show
    \[ \inf_{\substack{ M\in \R^{d\times d}:\\ \rank(M)=d}} NV({\bQ}\circ M)\ge \rank (\{A_1,\dots,A_n\}).\]
We need to show for any $M$,
\[ NV(\bQ\circ M)\ge \Rrank(\sum_{i=1}^nx_iA_i). \]
    The left hand side equals the number of nonzero rows of $M^T A_1 M,\dots, M^T A_n M$, which is $\ge $ the number of nonzero rows of $M^T (\sum_{i=1}^n x_iA_i) M$, which is $\ge \Rrank(\sum_{i=1}^nx_iA_i)$.

Next, we show
    \[ \inf_{\substack{ M\in \R^{d\times d}:\\ \rank(M)=d}} NV({\bQ}\circ M)\le \Rrank (\sum_{i=1}^nx_iA_i).\]
Denote $r=\Rrank (\sum_{i=1}^nx_iA_i)$. This means that the row-space of $\sum_{i=1}^n x_iA_i$ is generated by $r$ rows. With out loss of generality, we assume they are the first $r$ rows. Therefore, there exists invertible $M\in \R^{d\times d}$ such that $M^T (\sum_{i=1}^nx_iA_i)M$ has nonzero entries only in the left-top $r\times r$ block. This implies that the nonzero rows of $M^T A_1 M,\dots, M^T A_n M$ is $\le$ $r$.
\end{proof}

\begin{lemma}
Let $\bQ$ and $A_i$ be given as above. Then,
    \begin{equation}
        \inf_{\substack{M\in\R^{d\times d} \\ \rank (M)=d}}
    \inf_{\substack{M'\in \R^{n\times n}\\ \rank(M')=n'}} NV(M'\cdot \bQ\circ M)=\inf_{\substack{M'\in \R^{n\times n}\\ \rank(M')=n'}}\rank (\{\sum_{j=1}^n m'_{1j}A_j,\dots,\sum_{j=1}^n m'_{nj}A_j\}).
    \end{equation}
Here on the right hand side, $(m'_{ij})$ are the entries of $M'$.    
\end{lemma}

\begin{proof}
Denote $\widetilde{\bQ}=M'\cdot \bQ$. In other words, if $\widetilde{\bQ}=(\wt{Q_1},\dots, \wt{Q_n})$ then $\wt{Q_i}= \sum_{j=1}^n m'_{ij} Q_j$. We see that the corresponding matrix for $\wt{Q_i}$ is $\sum_{j=1}^n m'_{ij} A_j$.
By Lemma \ref{lemA6}, we have
\begin{align*}
    \inf_{\substack{M\in\R^{d\times d} \\ \rank (M)=d}} NV(M'\cdot \bQ\circ M)=\inf_{\substack{M\in\R^{d\times d} \\ \rank (M)=d}} NV(\wt{\bQ}\circ M)=\rank(\{\sum_{j=1}^n m'_{1j}A_j,\dots,\sum_{j=1}^n m'_{nj}A_j\}).
\end{align*}
Taking $\inf_{\substack{M'\in \R^{n\times n}\\ \rank(M')=n'}}$ on both sides will finish the proof.
\end{proof}

\begin{lemma}\label{0707.lema8}
Let $\bQ$ and $A_i$ be given as above. Then, 
    \begin{align}
        \label{lemA8left}&\inf_{\substack{M\in\R^{d\times d}\\ \rank (M)=d'}}\inf_{\substack{M'\in\R^{n\times n}\\ \rank(M')=n'}} NV(M'\cdot \bQ\circ M)\\
        \label{lemA8right}=&\inf_{\substack{M\in\R^{d\times d}\\ \rank (M)=d'}}\inf_{\substack{M'\in\R^{n\times n}\\ \rank(M')=n'}}\rank (\{\sum_{j=1}^n m'_{1j}M^TA_jM,\dots,\sum_{j=1}^n m'_{nj}M^TA_jM\}).
    \end{align} 
Here on the right hand side, $(m'_{ij})$ are the entries of $M'$.
\end{lemma}

\begin{proof}
   We see that \eqref{lemA8left} equals
   \[ \inf_{\substack{M\in\R^{d\times d}\\ \rank (M)=d'}}\inf_{\substack{N\in\R^{d\times d}\\ \rank (N)=d}}\inf_{\substack{M'\in\R^{n\times n}\\ \rank(M')=n'}} NV(M'\cdot \bQ\circ M\circ N). \]
Denote $\wt{\bQ}=M'\cdot \bQ\circ M$. In other words, if $\wt{\bQ}=(\wt{Q_1},\dots,\wt{Q_n})$ then $\wt{Q_i}=\sum_{i=1}^nm'_{ij}Q_j\circ M$. We see that the corresponding matrix for $\wt{Q_i}$ is $\sum_{i=1}^n m'_{ij}M^T A_j M$.

If we fix $M'$ and $M$, then by Lemma \ref{lemA6} we have
\begin{align*}
 \inf_{\substack{N\in\R^{d\times d}\\ \rank(N)=d}}NV(M'\cdot \bQ\circ M\circ N)=
    \inf_{\substack{N\in\R^{d\times d}\\ \rank(N)=d}}NV(\wt{\bQ}\circ N)\\
    =\rank (\{ \sum_{j=1}^n m'_{1j}M^TA_jM,\dots,\sum_{j=1}^n m'_{nj}M^TA_jM \}). 
\end{align*} 
By further taking $\inf_{\substack{M\in\R^{d\times d}\\ \rank (M)=d'}}\inf_{\substack{M'\in\R^{n\times n}\\ \rank(M')=n'}}$ on both sides, we finish the proof of the lemma.
\end{proof}

\bigskip

By the definition of $\fd_{d',n'}$, Lemma \ref{0707.lema8} and Definition \ref{defrank}, we obtain
\begin{corollary}\label{0717.a10}
\begin{equation}\label{rhsrrank}
    \fd_{d',n'}(\bQ)=\inf_{\substack{M\in\R^{d\times d}\\ \rank (M)=d'}}\inf_{\substack{N\in\R^{n\times n}\\ \rank(N)=n'}}\Rrank (\sum_{i=1}^n\sum_{j=1}^n x_i n_{ij}MA_jM^T).
\end{equation}
Here, $(n_{ij})$ are the entries of $N$, $x_i$'s are indeterminates and ``$\Rrank$" is defined in Definition \ref{0707.defa3}.
\end{corollary}

\subsection{Finding an algorithm}

The remaining work of this section is to find an algorithm to calculate the right hand side of \eqref{rhsrrank}.
At some point in the algorithm, the problem is reduced to determining whether a semi-algebraic set is empty or not. This step will be handled by the classical cylindrical decomposition from real algebraic geometry. We refer to \cite[Theorem 2.3.1]{BochnakCosteRoy}, which determines whether a semi-algebraic set is empty or not. 

\begin{lemma}\label{emptyornot}
    Suppose $Z$ is a semi-algebraic set and we explicitly know the defining polynomials of $Z$ (see Definition \ref{defsemialg}). Then there is an algorithm to determine whether $Z$ is empty or not. 
\end{lemma}

On the right hand side of \eqref{rhsrrank}, the infimum is taken over all the rank $d'$ matrices $M$ and rank $n'$ matrices $N$. We will use the following lemma to narrow the range of the matrices taken under the ``inf".  

\begin{lemma}\label{narrowrange}
    Suppose $M,M'\in \R^{d\times d}$ satisfy $M'=CM$ for some $C\in GL_d(\R)$, then
    \[\Rrank (\sum_{i=1}^n\sum_{j=1}^n x_i n_{ij}MA_jM^T)=\Rrank (\sum_{i=1}^n\sum_{j=1}^n x_i n_{ij}M'A_jM'^T).\]
    Suppose $N=(n_{ij}),N'=(n'_{ij})\in \R^{n\times n}$ satisfy $N'=CN$ for some $C\in GL_n(\R)$, then
    \[\Rrank (\sum_{i=1}^n\sum_{j=1}^n x_i n_{ij}MA_jM^T)=\Rrank (\sum_{i=1}^n\sum_{j=1}^n x_i n'_{ij}MA_jM^T).\]
\end{lemma}

\begin{proof}
    For the first part, we note
    \begin{align*}
        \Rrank (\sum_{i=1}^n\sum_{j=1}^n x_i n_{ij}M'A_jM'^T)&=\Rrank (C\sum_{i=1}^n\sum_{j=1}^n x_i n_{ij}MA_jM^TC^T)\\
        &=\Rrank (\sum_{i=1}^n\sum_{j=1}^n x_i n_{ij}MA_jM^T).
    \end{align*}
The last equality is by \eqref{rowrankeq}.

\medskip

Next, we prove the second part. Writing $C=(c_{ij})$, we have $n'_{ij}=\sum_l c_{il}n_{lj}$.
We have
\begin{align*}
    \Rrank (\sum_{i=1}^n\sum_{j=1}^n x_i n'_{ij}MA_jM^T)=\Rrank (\sum_{i=1}^n\sum_{j=1}^n\sum_{l=1}^n x_i c_{il}n_{lj}MA_jM^T).
\end{align*}
We define the new indeterminates $y_l=\sum_{i=1}^n x_i c_{il}$, or in other words,
\[ (y_1 \dots y_n)= (x_1\dots x_n) C, \]
if written in matrix.
We have
\begin{align*}
    \Rrank (\sum_{i=1}^n\sum_{j=1}^n x_i n'_{ij}MA_jM^T)=\Rrank (\sum_{i=1}^n\sum_{j=1}^n y_i n_{ij}MA_jM^T).
\end{align*}

Since $\{x_i\}_{i=1}^n$ are algebraically independent indeterminates and $C$ is invertible, we see that $\{y_i\}_{i=1}^n$ are also algebraically independent indeterminates.
Therefore, the expression above equals
\[\Rrank (\sum_{i=1}^n\sum_{j=1}^n x_i n_{ij}MA_jM^T).\]
This completes the proof.
\end{proof}

Lemma \eqref{narrowrange} tells us: if we want to calculate the right hand side of \eqref{rhsrrank}, it suffices to test those $M$ and $N$ with the row-echelon form. We briefly introduce the row-echelon form. 
Fix $d'\le d$.
For a sequence of integers $\e=(\e_1,\dots,\e_{d'})$ with $1\le\e_1<\e_2<\dots<\e_{d'}\le d$ define the set of matrices $\cM_\e(d,d')$ which consists of the matrices $M=(m_{ij})$ with the following properties. 
For $i>d'$, $m_{ij}=0$. For $i=1,\dots,d'$, $m_{l\e_i}=1$ when $l=i$; and $=0$ when $l\neq i$. We call $\cM_\e(d,d')$ the $d\times d$ $\e$-type row-echelon form with rank $d'$. We can view $\cM_\e(d,d')$ as $\R^{d'(d-d')}$.

For example, if $\e=(1,2,3)$, then
\[ \cM_\e(4,3)=\left\{ \begin{pmatrix}
    1 & 0 & 0 & a\\
    0 & 1 & 0 & b\\
    0 & 0 & 1 & c\\
    0& 0 & 0 & 0
\end{pmatrix}: a,b,c\in\R  \right\}. \]

Note that there are $\binom{d}{d'}$ many $\e$ that satisfies $1\le\e_1<\e_2<\dots<\e_{d'}\le d$. Denote $L=\binom{d}{d'}$. For simplicity, we just change the subscript of $\cM_\e(d,d')$, and denote all these different types of the row-echelon forms by
\[ \cM_1(d,d'),\cM_2(d,d'),\dots,\cM_L(d,d'). \]

Consider the matrix space 
\[ \{M\in \R^{d\times d}: \rank(M)=d'\}. \]
By a fundamental fact in linear algebra, for each $M\in \{M\in \R^{d\times d}: \rank(M)=d'\}$,
there exists a matrix $C\in GL_d(\R)$ so that 
$CM \in \cM_i(d,d')$ for some $i$.
 
By \eqref{rhsrrank} and Lemma \ref{narrowrange}, we have
\begin{equation}
    \fd_{d',n'}(\bQ)=\inf_{\substack{\al=1,\dots, \binom{d}{d'}\\ \beta=1,\dots,\binom{n}{n'}}}  \inf_{M\in \cM_\al(d,d')}\inf_{N\in \cM_\beta(n,n')}\Rrank (\sum_{i=1}^n\sum_{j=1}^n x_i n_{ij}MA_jM^T).
\end{equation}

Since there are finitely many $\al$ and $\beta$, we fix $\al$ and $\beta$, and it suffices to find an algorithm to calculate
\begin{equation}\label{infinf}
    \inf_{M\in \cM_\al(d,d')}\inf_{N\in \cM_\beta(n,n')}\Rrank (\sum_{i=1}^n\sum_{j=1}^n x_i n_{ij}MA_jM^T). 
\end{equation} 

Let us rewrite it. For any $M\in\cM_\al(d,d')$, there are $d'$ entries of $M$ that are $1$ with fixed position, and there are $d'(d-d')$ entries ranging in $\R$. We can think of $M\in \cM_\al(d,d')$ as $M(u)$ with $u$ ranging from $\R^{d'(d-d')}$.
Therefore, we have proved

\begin{corollary}
\[ \eqref{infinf}=\inf_{u\in\R^{d'(d-d')}}\inf_{v\in\R^{n'(n-n')}}\Rrank (R(u,v,x)). \]
Here $R(u,v,x)=\big(R_{ij}(u,v,x)\big)$ is an $n\times n$-matrix whose entries are polynomials of variables $u\in \R^{d'(d-d')},v\in\R^{n'(n-n')}$  and $x\in \R^n$.\footnote{We remark that the matrix $R$ also depends on a choice of $A_j$ which comes from $\bQ$. But since $\bQ$ is fixed, we omit it.}
\end{corollary}

Finally, we reduce the algorithm to the following problem.

Suppose $R(u,v,x)=\big(R_{ij}(u,v,x)\big)$ is an $n\times n$ matrix whose entries are polynomials of $u\in \R^a, v\in \R^b$ and $x\in \R^n$. {\bf Find an algorithm} to calculate
\begin{equation}\label{0718.1.8}
r:=\inf_{u\in\R^a}\inf_{v\in\R^b}\Rrank (R(u,v,x)).   
\end{equation}
Here $\{x_i\}_{i=1}^n$ are indeterminates, and we remind the reader the definition of ``$\Rrank$" in Definition \ref{0707.defa3}.
\\

The idea is as follows. Starting with $m=0$, we check whether $r\le m$. If not, we check whether $r\le m+1$; If yes, we obtain $r=m$.
To check whether $r\le m$, we note that following fact: $\Rrank (R(u_0,v_0,x))\le m$ for some $u_0,v_0$ is equivalent to that there are $m$ rows of $R(u_0,v_0,x)$ that span the row space of $R(u_0,v_0,x)$. For convenience, let us first consider the case that these $m$ rows are the first $m$ rows of $R(u_0,v_0,x)$. Other cases can be dealt with similarly.
The first $m$ rows span the row space of $R(u_0,v_0,x)$ is equivalent to that
\[ C R(u_0,v_0,x)=0, \]
for some $n\times n$ matrix $C$ of form
\begin{equation}
    C:=\begin{pmatrix}
        O_{m \times m} & O_{m \times (n-m)} \\
        A_{(n-m) \times m} & I_{(n-m) \times (n-m)}
    \end{pmatrix}.
\end{equation}
Here,
 the matrix $O_{m \times m}$ is an $m \times m$ matrix with zero entries (define $O_{m \times (n-m)}$ similarly),  $I_{(n-m) \times (n-m)}$ is an identity matrix, and $A_{(n-m) \times m}$ is any $(n-m) \times m$ matrix. The remaining entries (there are $(n-m)m$ of them) could be any real numbers. Denoting $c:=(n-m)m$, we can write $C$ as $C(w)$ where $w\in \R^{c}$. We have shown that: the first $m$ rows of $R(u_0,v_0,x)$ span the row space of $R(u_0,v_0,x)$ $\Leftrightarrow$ there exists $w_0\in\R^c$ such that
\[ C(w_0) R(u_0,v_0,x). \]

We considered the case that the first $m$ rows of $R(u_0,v_0,x)$ span the row space of $R(u_0,v_0,x)$. In general, there are $\binom{n}{m}$ many choices of choosing $m$ rows out of $n$ rows. By following the same argument, there are $C_1(w),\dots, C_L(w)$ with $L=\binom{n}{m}$. Each $m$ rows corresponds to a $C_i(w)$ so that: these $m$ rows of $R(u_0,v_0,x)$ span the row space of $R(u_0,v_0,x)$ $\Leftrightarrow$ there exists $w_0\in\R^c$ such that
\[ C_i(w_0) R(u_0,v_0,x). \]

In summary, we proved: $r\le m$ (see \eqref{0718.1.8} for the definition of $r$) $\Leftrightarrow$
there exist $u_0\in\R^a,v_0\in\R^b, w_0\in \R^c, 1\le i\le L$ such that
\[ C_i(w_0) R(u_0,v_0,x)=0 \]
for all $x\in \R^n$.
Hence, checking whether $r\le m$ is equivalent to checking whether
\begin{equation}\label{returntosolution}
    C_i(w) R(u,v,\cdot)\equiv 0
\end{equation} 
has a solution $(w_0,u_0,v_0)$ for some $i$. A good thing here is that $C_i(w)R(u,v,x)$ is polynomially depending on $w,u,v,x$. We will  use the following lemma.

\begin{lemma}\label{polynomial0}
    Let $P(x)$ be a polynomial of degree $D$ with $n$-variables. If $P(x)=0$ for all $x\in \{1,2,\dots,D+1\}^n$, then $P(x)\equiv 0$.
\end{lemma}
\begin{proof}
    Write $x=(y,x_n) \in \R^{n-1} \times \R$. Write $P(x)=\sum_{i=0}^D a_i(y)x_n^i$. Note that $a_i$ is a polynomial. Given $y$, by the fundamental theorem of algebra, $P(y,x_n)=0$ for $x_n\in\{1,2,\dots,D+1\}$ implies $a_i(y)=0$ for all $i=0,\dots,D$. By the same process, we can reduce the number of variables by one at each time. Finally we can deduce that $P(x)\equiv 0$.
\end{proof}

Let us return to \eqref{returntosolution}. 
Denote the $(k,l)$-entry of $C_i(w) R(u,v,x)$ by $[C_i(w) R(u,v,x)]_{kl}$, which is a polynomial of $w,u,v,x$.
Define
\[ D:=\sup_i \sup_{k,l}\deg \big([C_i(w) R(u,v,x)]_{kl} \big). \]
 Consider the algebraic variety
\[ Z=\bigcup_i \{ (u,v,w): C_i(w) R(u,v,x)=0 \;\; \mathrm{for \; all} \; x\in \{1,2,\dots,D+1\}^n  \}. \]
Then, \eqref{returntosolution} has a solution $(w_0,u_0,v_0)$ for some $i$ is equivalent to $Z\neq \emptyset$. The reason is as follows. If \eqref{returntosolution} has a solution $(u_0,v_0,w_0)$ for some $i$, then $(u_0,v_0,w_0)\in Z$ which is not empty. If $(u_0,v_0,w_0)\in Z$, then we have
\[ C_i(w_0) R(u_0,v_0,x)=0 \]
for all $x\in \{1,2,\dots,D+1\}^n$.
By Lemma \ref{polynomial0}, we can deduce that
\[C_i(w_0) R(u_0,v_0,\cdot)\equiv 0.\]
Determining whether $Z$ is empty or not is by Lemma \ref{emptyornot}.
This finishes the algorithm.

\section{Algorithm to compute \texorpdfstring{$X(\cM,k,m)$}{}}
The goal of this section is to find an algorithm to calculate the number $X(\cM,k,m)$ defined in Definition \ref{0618.def12}.
We remind the reader that $\bQ=(Q_1,\dots,Q_n)$ is an $n$-tuple of quadratic polynomials, and $V_\xi$ is the tangent space of the manifold $\{(\xi,\bQ(\xi)):\xi\in\R^d\}$ at the point $(\xi, {\bf Q}(\xi))$.
 By the definition of $X$, it suffices to find an algorithm to calculate
\[ \sup_{\dim V=m}\dim(\{\xi\in\R^d:\dim(\pi_{V_\xi}(V))<X\}), \]
for all the integers $1\le m\le d+n, 0\le X\le d+1$.

We can represent $V_\xi$ using the following matrix. 
\begin{equation}\label{matrixVxi}
    \begin{pmatrix}
    1 & 0 &\cdots &0\\
    0 & 1 &\cdots & 0\\
    \vdots & \vdots &\ddots &\vdots\\
    0 & 0& \cdots & 1\\
    \p_1 Q_1 & \p_2 Q_1 &\cdots & \p_d Q_1\\
    \p_1 Q_2 & \p_2 Q_2 &\cdots & \p_d Q_2\\
    \vdots & \vdots & \ddots & \vdots\\
    \p_1 Q_n & \p_2 Q_n &\cdots & \p_d Q_n
\end{pmatrix}.
\end{equation} 
The $j$-th column of the matrix is $(\nabla_j \xi, \nabla_j \bQ(\xi))^T$. By abuse of  notation, we still use $V_\xi$ to denote this $(d+n)\times d$ matrix. We see that the tangent space $V_\xi$ is spanned by the column vectors of the matrix $V_\xi$.

Let $G(k, n)$ be the set of $k$-dimensional subspaces in $\R^n$.
We will represent $V\in G(m,d+n)$ using matrix. For $V\in G(m,d+n)$, we choose $m$ vectors $v_1,\dots,v_m$  spanning the space $V$. We write the $m\times (d+n)$ matrix
\begin{equation}\label{notunique}
    \begin{pmatrix}
    v_1\\
    v_2\\
    \vdots\\
    v_m
\end{pmatrix}. 
\end{equation} 
By abuse of notation, we denote this matrix by $V$. We actually find a correspondence
\begin{equation}\label{0718.23} G(m,d+n)\longrightarrow \{ V\in \R^{m\times (d+n)}: \rank(V)=m  \}=:\cM.
\end{equation}
We remark that this map is not unique because of the choice \eqref{notunique}. What matters here is that we have
\[ \dim(\pi_{V_\xi}(V))=\rank(VV_\xi). \]
On the left hand side, we view $V_\xi, V$ as subspaces; while on the right hand side, we view $V_\xi, V$ as matrices. Therefore, we obtain
\begin{equation}
    \sup_{\dim V=m}\dim(\{\xi:\dim(\pi_{V_\xi}(V)<X)\})=\sup_{V\in \cM}\dim(\{ \xi: \rank(VV_\xi)<X \}).
\end{equation}
Recall that $\mathcal{M}$ is defined in \eqref{0718.23}.
Similar to Appendix A, we can find $L:=\binom{d+n}{m}$ many types of the row-echelon forms $\cM_1, \cM_2,\dots \cM_L \subset \mathcal{M}$. For example, $\cM_1$ consists of matrices of form
\[ \begin{pmatrix}
    1 & 0 & \cdots & 0 & \cdots\\
    0 & 1 & \cdots & 0 & \cdots\\
    \vdots & \vdots & \ddots & \vdots &\cdots \\
    0 & 0 & \cdots & 1 & \cdots\\
\end{pmatrix}=
\begin{pmatrix}
    I_m & M
\end{pmatrix}, \] for some
$M\in \R^{m\times (d+n-m)}$. For every $V\in \cM$, there is some $C\in GL_m(\R)$ so that $CV \in \cM_i$ for some $i$.
Since the rank of a matrix is invariant under the left-multiplication by a invertible matrix, we have
\[\sup_{V\in \cM}\dim(\{ \xi: \rank(VV_\xi)<X \})=\sup_{1\le i\le L}\sup_{V\in\cM_i}\dim(\{ \xi: \rank(VV_\xi)<X \}).\]
Since there are finitely many $i$, we just need to find an algorithm to calculate 
\[\sup_{V\in\cM_i}\dim(\{ \xi: \rank(VV_\xi)<X \}).\]
For simplicity, we write the matrix in $\cM_i$ as a function $V_i(\eta)$ with variables $\eta\in \R^{m (d+n-m)}$. After this identification, we have
\begin{equation}\label{comeback1}
    \sup_{V\in\cM_i}\dim(\{ \xi: \rank(VV_\xi)<X \})=\sup_{\eta\in\R^{m(d+n-m)}}\dim(\{ \xi: \rank(V_i(\eta)V_\xi)<X \}).
\end{equation}

We will use the following lemma to characterize the rank. 
    The proof of the lemma is straightforward, so we omit the details.

\begin{lemma}\label{charrank}
    Let $V=(v_{ij})$ be an $n\times m$ matrix. For an integer $X$, let the polynomial $P(v_{ij})$ be the square sum of all the determinant of the $X\times X$-minor of $V$. More precisely,
    \[ P(v_{ij}):=\sum_{V' \textup{~is~a~}X\times X\textup{~minor~of~}V}|\det(V')|^2. \]
Then $P$ is a polynomial in terms of $\{v_{ij}\}$ and $\rank(V)<X\Leftrightarrow P(v_{ij})=0$. 
\end{lemma}

By Lemma \ref{charrank}, there is a polynomial $P(\eta,\xi)$ 
 with $\eta\in\R^{m(d+n-m)}$ and $\xi \in \R^d$ so that
\[ \eqref{comeback1}=\sup_{\eta\in\R^{m(d+n-m)}}\dim(\{\xi\in\R^d: P(\eta,\xi)=0\}).\]
We can also calculate the coefficients of $P(\eta,\xi)$ explicitly from the given  $\bQ$.

It boils down to the following slicing problem. Let $P(\eta,\xi)$ be a polynomial with variables $(\eta,\xi)\in \R^a\times \R^b$. Recall that $Z_P=\{(\eta,\xi): P(\eta,\xi)=0\}$. {\bf Find an algorithm} to calculate the maximal dimension of the $\eta$-slice, i.e., to find
\begin{equation}\label{supslice}
    \sup_{y\in\R^a} \dim(Z_P\cap\{\eta=y\}).
\end{equation} 
We  first prove the following result.
\begin{theorem}\label{thmslice1}
    For any integer $m$, the set
    \begin{equation}\label{etaslice}
        \{\eta\in\R^a: \dim(\{\xi\in\R^b:P(\eta,\xi)=0\})= m\} 
    \end{equation}
    is a semi-algebraic set, and its defining polynomials (see Definition \ref{defsemialg}) can be expressed explicitly using the coefficients of $P(\eta,\xi)$.
\end{theorem}

If the theorem is true, then for any $m$, we have that \eqref{etaslice} is semi-algebraic. By Lemma \ref{emptyornot}, we can determine whether \eqref{etaslice} is empty or not. Therefore, the quantity \eqref{supslice}, which is the maximal integer $m$ for which \eqref{etaslice} is non-empty, can be calculated.

\begin{proof}[Proof of Theorem \ref{thmslice1}]
    We  show that $\{\eta\in\R^a: \dim(\{\xi\in\R^b:P(\eta,\xi)=0\})\ge m\} $ is a semi-algebraic set and can be expressed explicitly using the coefficients of $P(\eta,\xi)$. If this is true, then the subtraction will give the same property for 
    \[\{\eta\in\R^a: \dim(\{\xi\in\R^b:P(\eta,\xi)=0\})= m\} .\]

Denote $L:=\binom{b}{m}$. Let $\{V_i\}_{i=1}^L$ be the set of $m$-dimensional subspaces in $\R^b$ so that each of them is spanned by $m$ vectors from $\{\vec e_1,\dots, \vec e_b\}$ which are the coordinates of $\R^b$. We also denote $\pi_{V_i}$ to be the orthogonal projection onto $V_i$. We have the following observation. For fixed $\eta \in \R^a$,
\[ \dim(\{\xi:P(\eta,\xi)=0\})\ge m \Leftrightarrow \exists i\ \textup{s.t.~} \pi_{V_i}(\{\xi:P(\eta,\xi)=0\}) \textup{~contains~an~interior~point}. \]
Here is a sketch proof. ``$\Leftarrow$" is trivial. For $``\Rightarrow"$, note that $\{\xi:P(\eta,\xi)=0\}\subset \R^b$ is a variety. If the variety has dimension $m$, then there is a small patch of the variety that is diffeomorphic to an $m$-dimensional manifold. Hence, there exists $V_i$ so that the projection of this small patch of the variety onto $V_i$ contains an interior point. 

As a result, we have
\begin{align}
    &\{\eta\in \R^a: \dim(\{\xi\in\R^b:P(\eta,\xi)=0\})\ge m\}\\
    &=\bigcup_{i=1}^L\{ \eta\in\R^a: \pi_{V_i}(\{\xi\in\R^b:P(\eta,\xi)=0\}) \textup{~contains~an~interior~point} \}.
\end{align}
We just need to deal with one of the $V_i$. For simplicity, we consider the $m$-dimensional subspace $V\subset \R^b$ which is spanned by the first $m$ coordinates $\{x_1,\dots,x_m\}$. We show that
\begin{equation}\label{theset}
    \{ \eta\in\R^a: \pi_{V}(\{\xi\in\R^b:P(\eta,\xi)=0\}) \textup{~contains~an~interior~point} \}
\end{equation}
is semi-algebraic and whose defining polynomials can be expressed explicitly using the coefficients of $P(\eta,\xi)$.

To do this, we will use a decent tool from real algebraic geometry known as the first-order formula. For the definition of the first-order formula, we refer to \cite[Definition 2.2.3]{BochnakCosteRoy}. We will express \eqref{theset} using the first order formula, and apply \cite[Proposition 2.2.4]{BochnakCosteRoy} to show that it is semi-algebraic.

Let us give the details.
It is not hard to see that the set \eqref{theset} equals to
\begin{align}\label{1order}
    \bigg\{\eta\in\R^a: \exists \e>0, \exists \xi_0\in\R^b, \forall 0\le \de_j\le \e\  (j=1,\dots,b),\\
    \nonumber \xi_0+(\de_1,\dots,\de_b)\in \pi_V(\{\xi\in \R^b: P(\eta,\xi)=0\}) \bigg\}
\end{align}
By Tarski's projection theorem, $\pi_V(\{\xi\in \R^b: P(\eta,\xi)=0\})$ is semi-algebraic and its defining polynomials can be expressed explicitly using the coefficients of $P(\eta,\xi)$. Therefore, \eqref{1order} is expressed using the first-order formula, and hence is a semi-algebraic set whose defining polynomials can be expressed explicitly using the coefficients of $P(\eta,\xi)$ (see \cite[Proposition 2.2.4]{BochnakCosteRoy}).
\end{proof}

\section{Proof of Theorem \ref{0705.thm16}}

 In this section, we prove Theorem \ref{0705.thm16}.

 Consider the paraboloid $\mathcal{M}$ in $\R^{d+1}$. We are interested in the range of  
 \begin{equation}
    p> \min_{2 \leq k \leq d+1} \max{(\frac{2(2d-k+4)}{2d-k+2},\frac{2k}{k-1}  )}=:p_c.
\end{equation}
By taking $k=\lfloor \frac{2(d+2)}{3} \rfloor$, we get the bound
\begin{equation}
  p>2+\frac{6}{d}+O(\frac{1}{d^2}).
\end{equation}
 By Theorem \ref{restriction} and the epsilon removal lemma by \cite{MR1666558}, it suffices to prove that   \begin{equation*}
\begin{split}
 \delta^{d-\frac{2d+2n}{p}}\delta^{-\mathrm{Dec}_p(\mathcal{M}|_{L_{k-2}} )}
+
       \sup_{0 \leq m \leq d+n}(\delta^{-\frac{m}{p}+\frac12{X(\mathcal{M},k,m)}}) \lesssim_{\epsilon} \delta^{-\epsilon}
\end{split}
\end{equation*}
for every $p>p_c$ and $2 \leq k \leq d+1$.
The above inequality is equivalent to
\begin{equation}\label{0706.22}
    \delta^{-\mathrm{Dec}_p(\mathcal{M}|_{L_{k-2}} )} \lesssim  \delta^{-d+\frac{2d+2}{p}}
\end{equation}
and
\begin{equation}\label{0706.23}
    \frac{2m}{p} \leq X(\mathcal{M},k,m)
\end{equation}
for all $0 \leq m \leq d+1$.
\\

Let us show that \eqref{0706.22} and \eqref{0706.23} are true for $p>p_c$.

Consider $E^{\mathcal{M}}_Lf$ where $L$ is the $\delta$-neighborhood of a $(k-2)$-dimensional linear subspace. Fix a ball $B_{\delta^{-2}}$. By an uncertainty principle, after a linear transformation, the Fourier transform of $E^{\mathcal{M}}_{L}f$ is contained in the $\delta^2$-neighborhood of 
\begin{equation}
    \big\{(\xi_1,\ldots,\xi_{k-1}): \xi_{k-1}=\sum_{i=1}^{k-2}\xi_i^2 \big\} \times \R^{d-k+2}.
\end{equation}
Hence, by a decoupling for the $(k-2)$-dimensional paraboloid \cite{MR3374964}, we have
\begin{equation}
    \delta^{-\gamma_p(\mathcal{M}|_{L_{k-2}} )} \lesssim_{\epsilon}\delta^{-\epsilon} \max (\delta^{-(k-2)(\frac12-\frac1p)} ,\delta^{-(k-2)+\frac{2(k-1)}{p}}).
\end{equation}
Using this inequality, one can see that $\eqref{0706.22}$ is true if and only if 
\begin{equation}
    -(k-2)(\frac12-\frac1p) \geq -d+\frac{2d+2}{p}.
\end{equation}
By routine computations, this is equivalent to
\begin{equation}
    p>\frac{2(2d-k+4)}{2d-k+2}.
\end{equation}
This gives the range of $p$ for which \eqref{0706.22} holds true.
\\

Let us next calculate 
$X(\mathcal{M},k,m)$.
We need to find the largest integer $X$ satisfying
\begin{equation}
\sup_{\dim V=m}
   \dim \{ \xi \in \mathbb{R}^d: \mathrm{dim} (\pi_{V_\xi}(V)) < X  \} \leq k-2.
\end{equation}
Denote the basis of $V$ by $\{v_1,\ldots,v_m \}$. Then we have
\begin{equation}
\pi_{V_\xi}(V)=
\begin{pmatrix}
    1 & 0 & \cdots & 0 & 2\xi_1 \\
    0 & 1 & \cdots & 0 & 2\xi_2 \\
    \vdots & \vdots & \ddots & \vdots & \vdots \\
    0 & 0 & \cdots & 1 & 2\xi_d
\end{pmatrix} \begin{pmatrix}
    v_1 & \cdots & v_m
\end{pmatrix}.
\end{equation}
The case $m=d+1$ is trivial. Note that $\dim(\pi_{V_\xi}(V))=d$ for all $\xi \in \R^d$. Since $2 \leq k \leq d+1$, we have 
\begin{equation}
X(\mathcal{M},k,d+1)=d.
\end{equation}

Consider the case  $m \leq d$.
Note that
\begin{equation}
    m-1 \leq \dim \pi_{V_\xi}(V) \leq m.
\end{equation}
Also, one can see that
\begin{equation}
\begin{split}
& \sup_{\dim V=m}
   \dim \{ \xi \in \mathbb{R}^d: \mathrm{dim} (\pi_{V_\xi}(V)) \geq  m+1  \} =0,
\\&\sup_{\dim V=m}
   \dim \{ \xi \in \mathbb{R}^d: \mathrm{dim} (\pi_{V_\xi}(V)) = m  \} =d,
   \\&
   \sup_{\dim V=m}
   \dim \{ \xi \in \mathbb{R}^d: \mathrm{dim} (\pi_{V_\xi}(V)) = m-1  \} =m-1,
   \\&
   \sup_{\dim V=m}
   \dim \{ \xi \in \mathbb{R}^d: \mathrm{dim} (\pi_{V_\xi}(V)) \leq m-2  \} = 0.
\end{split}
\end{equation}
This gives the value of $X$;
\begin{equation}
X(\mathcal{M},k,m)=
    \begin{cases}
        m \;\;\;\;\;\;\;\;\; \mathrm{ for }\; m \leq k-1 \\
 m-1 \;\;\; \mathrm{ for }\; m > k-1    \end{cases}
\end{equation}
Hence, \eqref{0706.23} holds true for $p \geq \frac{2k}{k-1}$.
\\

To summarize, \eqref{0706.22} and \eqref{0706.23} hold true for 
\begin{equation}
    p> \min_{2 \leq k \leq d+1} \max{\Big(\frac{2(2d-k+4)}{2d-k+2},\frac{2k}{k-1}  \Big)}.
\end{equation}
This finishes the proof of Theorem \ref{0705.thm16}.



\bibliographystyle{alpha}
\bibliography{reference}

\begin{thebibliography}{GOZZK23}

\bibitem[BCR98]{BochnakCosteRoy}
Jacek Bochnak, Michel Coste, and Marie-Fran\c{c}oise Roy.
\newblock {\em Real algebraic geometry}, volume~36 of {\em Ergebnisse der
  Mathematik und ihrer Grenzgebiete (3) [Results in Mathematics and Related
  Areas (3)]}.
\newblock Springer-Verlag, Berlin, 1998.
\newblock Translated from the 1987 French original, Revised by the authors.

\bibitem[BCT06]{MR2275834}
Jonathan Bennett, Anthony Carbery, and Terence Tao.
\newblock On the multilinear restriction and kakeya conjectures.
\newblock {\em Acta Math.}, 196(2):261--302, 2006.

\bibitem[BD15]{MR3374964}
Jean Bourgain and Ciprian Demeter.
\newblock The proof of the {$l^2$} decoupling conjecture.
\newblock {\em Ann. of Math. (2)}, 182(1):351--389, 2015.

\bibitem[BDG16]{MR3548534}
Jean Bourgain, Ciprian Demeter, and Larry Guth.
\newblock Proof of the main conjecture in {V}inogradov's mean value theorem for
  degrees higher than three.
\newblock {\em Ann. of Math. (2)}, 184(2):633--682, 2016.

\bibitem[BG11]{MR2860188}
Jean Bourgain and Larry Guth.
\newblock Bounds on oscillatory integral operators based on multilinear
  estimates.
\newblock {\em Geom. Funct. Anal.}, 21(6):1239--1295, 2011.

\bibitem[BGZZK21]{basu2021stationary}
Saugata Basu, Shaoming Guo, Ruixiang Zhang, and Pavel Zorin-Kranich.
\newblock A stationary set method for estimating oscillatory integrals, 2021.

\bibitem[BL04]{MR2064058}
Jong-Guk Bak and Sanghyuk Lee.
\newblock Restriction of the {F}ourier transform to a quadratic surface in
  {$\Bbb R^n$}.
\newblock {\em Math. Z.}, 247(2):409--422, 2004.

\bibitem[BLL17]{MR3694011}
Jong-Guk Bak, Jungjin Lee, and Sanghyuk Lee.
\newblock Bilinear restriction estimates for surfaces of codimension bigger
  than 1.
\newblock {\em Anal. PDE}, 10(8):1961--1985, 2017.

\bibitem[Bou91]{MR1097257}
J.~Bourgain.
\newblock Besicovitch type maximal operators and applications to {F}ourier
  analysis.
\newblock {\em Geom. Funct. Anal.}, 1(2):147--187, 1991.

\bibitem[Chr82]{christthesis}
Michael Christ.
\newblock Restriction of the fourier transform to submanifolds of low
  codimension.
\newblock {\em Ph.D. thesis, University of Chicago}, 1982.

\bibitem[Chr85]{MR766216}
Michael Christ.
\newblock On the restriction of the {F}ourier transform to curves: endpoint
  results and the degenerate case.
\newblock {\em Trans. Amer. Math. Soc.}, 287(1):223--238, 1985.

\bibitem[CL17]{MR3653943}
Chu-Hee Cho and Jungjin Lee.
\newblock Improved restriction estimate for hyperbolic surfaces in
  {$\mathbb{R}^3$}.
\newblock {\em J. Funct. Anal.}, 273(3):917--945, 2017.

\bibitem[Dem18]{MR3966819}
Ciprian Demeter.
\newblock Decouplings and applications.
\newblock In {\em Proceedings of the {I}nternational {C}ongress of
  {M}athematicians---{R}io de {J}aneiro 2018. {V}ol. {III}. {I}nvited
  lectures}, pages 1539--1560. World Sci. Publ., Hackensack, NJ, 2018.

\bibitem[DGL17]{MR3702674}
Xiumin Du, Larry Guth, and Xiaochun Li.
\newblock A sharp {S}chr\"{o}dinger maximal estimate in {$\Bbb R^2$}.
\newblock {\em Ann. of Math. (2)}, 186(2):607--640, 2017.

\bibitem[DZ19]{MR3961084}
Xiumin Du and Ruixiang Zhang.
\newblock Sharp {$L^2$} estimates of the {S}chr\"{o}dinger maximal function in
  higher dimensions.
\newblock {\em Ann. of Math. (2)}, 189(3):837--861, 2019.

\bibitem[GO22]{MR4525760}
Shaoming Guo and Changkeun Oh.
\newblock Fourier restriction estimates for surfaces of co-dimension two in
  {$\Bbb{R}^5$}.
\newblock {\em J. Anal. Math.}, 148(2):471--499, 2022.

\bibitem[GOZZK23]{MR4541334}
Shaoming Guo, Changkeun Oh, Ruixiang Zhang, and Pavel Zorin-Kranich.
\newblock Decoupling inequalities for quadratic forms.
\newblock {\em Duke Math. J.}, 172(2):387--445, 2023.

\bibitem[Gre22]{gressman2022testing}
Philip~T Gressman.
\newblock Testing conditions for multilinear radon-brascamp-lieb inequalities,
  2022.

\bibitem[Gre23]{gressman2023local}
Philip~T. Gressman.
\newblock Local curvature of maximally nondegenerate radon-like transforms,
  2023.

\bibitem[Gut15]{MR3300318}
Larry Guth.
\newblock A short proof of the multilinear {K}akeya inequality.
\newblock {\em Math. Proc. Cambridge Philos. Soc.}, 158(1):147--153, 2015.

\bibitem[Gut16]{MR3454378}
Larry Guth.
\newblock A restriction estimate using polynomial partitioning.
\newblock {\em J. Amer. Math. Soc.}, 29(2):371--413, 2016.

\bibitem[Gut18]{guth2018}
Larry Guth.
\newblock Restriction estimates using polynomial partitioning ii.
\newblock {\em Acta Math.}, 221(1):81--142, 09 2018.

\bibitem[Gut22]{guth2022decoupling}
Larry Guth.
\newblock Decoupling estimates in fourier analysis, 2022.

\bibitem[GWZ23]{guo2023dichotomy}
Shaoming Guo, Hong Wang, and Ruixiang Zhang.
\newblock A dichotomy for h\"ormander-type oscillatory integral operators,
  2023.

\bibitem[GZ18]{MR3830894}
Larry Guth and Joshua Zahl.
\newblock Polynomial {W}olff axioms and {K}akeya-type estimates in {$\Bbb
  R^4$}.
\newblock {\em Proc. Lond. Math. Soc. (3)}, 117(1):192--220, 2018.

\bibitem[GZK20]{MR4143735}
Shaoming Guo and Pavel Zorin-Kranich.
\newblock Decoupling for certain quadratic surfaces of low co-dimensions.
\newblock {\em J. Lond. Math. Soc. (2)}, 102(1):319--344, 2020.

\bibitem[HBP17]{MR3652248}
D.~R. Heath-Brown and L.~B. Pierce.
\newblock Simultaneous integer values of pairs of quadratic forms.
\newblock {\em J. Reine Angew. Math.}, 727:85--143, 2017.

\bibitem[HI22]{MR4405679}
Jonathan Hickman and Marina Iliopoulou.
\newblock Sharp {$L^p$} estimates for oscillatory integral operators of
  arbitrary signature.
\newblock {\em Math. Z.}, 301(1):1143--1189, 2022.

\bibitem[Hic23]{hickman2023pointwise}
Jonathan Hickman.
\newblock Pointwise convergence for the schr\"odinger equation [after xiumin du
  and ruixiang zhang], 2023.

\bibitem[HR19]{HR2019}
Jonathan Hickman and Keith~M. Rogers.
\newblock Improved fourier restriction estimates in higher dimensions.
\newblock {\em Camb. J. Math.}, 7(3):219--282, 2019.

\bibitem[HRZ22]{MR4521046}
Jonathan Hickman, Keith~M. Rogers, and Ruixiang Zhang.
\newblock Improved bounds for the {K}akeya maximal conjecture in higher
  dimensions.
\newblock {\em Amer. J. Math.}, 144(6):1511--1560, 2022.

\bibitem[HZ20]{hickman2020note}
Jonathan Hickman and Joshua Zahl.
\newblock A note on fourier restriction and nested polynomial wolff axioms,
  2020.

\bibitem[KR18]{MR3881832}
Nets~Hawk Katz and Keith~M. Rogers.
\newblock On the polynomial {W}olff axioms.
\newblock {\em Geom. Funct. Anal.}, 28(6):1706--1716, 2018.

\bibitem[LY21]{li2021decoupling}
Jianhui Li and Tongou Yang.
\newblock Decoupling for smooth surfaces in $\mathbb{R}^3$, 2021.

\bibitem[Mal22]{MR4395082}
Dominique Maldague.
\newblock Regularized {B}rascamp-{L}ieb inequalities and an application.
\newblock {\em Q. J. Math.}, 73(1):311--331, 2022.

\bibitem[Moc96]{Mockenhaupt}
Gerd Mockenhaupt.
\newblock Bounds in {L}ebesgue spaces of oscillatory integral operators.
\newblock {\em Habilitationsschrift, Univ.-GHS Siegen, Siegen}, 1996.

\bibitem[Obe02]{articleoberlin}
Daniel~M. Oberlin.
\newblock A restriction theorem for a $k$-surface in $\mathbb{R}^n$.
\newblock {\em Canadian Mathematical Bulletin}, 48, 10 2002.

\bibitem[Oh18]{MR3848437}
Changkeun Oh.
\newblock Decouplings for three-dimensional surfaces in {$\R^6$}.
\newblock {\em Math. Z.}, 290(1-2):389--419, 2018.

\bibitem[Pie19]{MR3939285}
Lillian~B. Pierce.
\newblock The {V}inogradov mean value theorem [after {W}ooley, and {B}ourgain,
  {D}emeter and {G}uth], 2019.
\newblock S\'{e}minaire Bourbaki. Vol. 2016/2017. Expos\'{e}s 1120--1135.

\bibitem[PS07]{pramanik2007p}
Malabika Pramanik and Andreas Seeger.
\newblock $l^p$ regularity of averages over curves and bounds for associated
  maximal operators.
\newblock {\em American journal of mathematics}, 129(1):61--103, 2007.

\bibitem[Tao99]{MR1666558}
Terence Tao.
\newblock The {B}ochner-{R}iesz conjecture implies the restriction conjecture.
\newblock {\em Duke Math. J.}, 96(2):363--375, 1999.

\bibitem[Tao03]{MR2033842}
Terence Tao.
\newblock A sharp bilinear restrictions estimate for paraboloids.
\newblock {\em Geom. Funct. Anal.}, 13(6):1359--1384, 2003.

\bibitem[TVV98]{MR1625056}
Terence Tao, Ana Vargas, and Luis Vega.
\newblock A bilinear approach to the restriction and {K}akeya conjectures.
\newblock {\em J. Amer. Math. Soc.}, 11(4):967--1000, 1998.

\bibitem[Wan22]{MR4484215}
Hong Wang.
\newblock A restriction estimate in {$\Bbb R^3$} using brooms.
\newblock {\em Duke Math. J.}, 171(8):1749--1822, 2022.

\bibitem[WW22]{wang2022improved}
Hong Wang and Shukun Wu.
\newblock An improved restriction estimate in $\mathbb{R}^3$, 2022.

\bibitem[Zah12]{MR3231483}
Joshua Zahl.
\newblock On the {W}olff circular maximal function.
\newblock {\em Illinois J. Math.}, 56(4):1281--1295, 2012.

\bibitem[Zah18]{MR3820441}
Joshua Zahl.
\newblock A discretized {S}everi-type theorem with applications to harmonic
  analysis.
\newblock {\em Geom. Funct. Anal.}, 28(4):1131--1181, 2018.

\bibitem[Zah21]{MR4205111}
Joshua Zahl.
\newblock New kakeya estimates using gromov's algebraic lemma.
\newblock {\em Adv. Math.}, 380:107596, 2021.

\bibitem[Zah23]{zahl2023maximal}
Joshua Zahl.
\newblock On maximal functions associated to families of curves in the plane,
  2023.

\end{thebibliography}

\end{document}